\RequirePackage{tcolorbox} 

\documentclass[sn-mathphys]{sn-jnl}

\usepackage{tcolorbox}  
\usepackage{accents}	
\usepackage{stmaryrd}   
\usepackage{subcaption} 
\usepackage{overpic} 
\usepackage{ulem} 

\jyear{2021}%

\theoremstyle{thmstyleone}%
\theoremstyle{thmstyletwo}%
\newtheorem{remark}{Remark}%

\theoremstyle{thmstylethree}%

\newcommand{\norm}[1]{\left\lVert#1\right\rVert}

\newcommand\stateL[1]{\mathbf{#1}}
\newcommand\stateG[1]{\boldsymbol #1}
\newcommand{\bigpartialderiv}[2]{ \frac{\partial {#1}}{\partial {#2} } }
\newcommand{\Nabla} {\vec{\nabla}}
\newcommand\acclrvec[1]{\accentset{\,\leftrightarrow}{#1}}	
\newcommand{\blocktensor}[1]{\acclrvec{{\mathbf #1}}}
\newcommand{\DG}{{\mathrm{DG}}}
\newcommand{\FV}{{\mathrm{FV}}}
\newcommand{\numfluxb}[1]{\hat{\mathbf{#1}} }
\newcommand{\numflux}[1]{\hat{{#1}} }
\newcommand{\mat}[1]{\underline{\mathbf{#1}}}
\newcommand{\ii}{j} 
\newcommand{\jump}[1]{\ensuremath{\left\llbracket #1 \right\rrbracket}}

\newcommand\eqLGL{\mathrel{\stackrel{\makebox[0pt]{\mbox{\normalfont\tiny LGL}}}{=}}}

\newcommand{\entVar}{{\mathbf{v}}}

\def\NN{\mathcal{N}} 

\newcommand{\numDims}{D}

\begin{document}

\title[MCL for LGL-DGSEM]{Monolithic Convex Limiting for Legendre--Gauss--Lobatto Discontinuous Galerkin Spectral Element Methods}


\author*[1]{\fnm{Andrés M} \sur{Rueda-Ramírez}}\email{aruedara@uni-koeln.de}

\author[1]{\fnm{Benjamin} \sur{Bolm}}\email{bbolm@smail.uni-koeln.de}

\author[2]{\fnm{Dmitri} \sur{Kuzmin}}\email{kuzmin@math.uni-dortmund.de}

\author[1,3]{\fnm{Gregor J} \sur{Gassner}}\email{ggassner@uni-koeln.de}

\affil*[1]{\orgdiv{Department of Mathematics and Computer Science}, \orgname{University of Cologne}, \orgaddress{\street{Weyertal 86-90}, \city{Cologne}, \postcode{50931}, 
\country{Germany}}}

\affil[2]{\orgdiv{Institute of Applied Mathematics (LS III)}, \orgname{TU Dortmund University}, \orgaddress{\street{Vogelpothsweg 87}, \city{Dortmund}, \postcode{D-44227}, 
\country{Germany}}}

\affil[3]{\orgdiv{Center for Data and Simulation Science}, \orgname{University of Cologne}, \orgaddress{\street{Weyertal 86-90}, \city{Cologne}, \postcode{50931}, 
\country{Germany}}}


\abstract{We extend the monolithic convex limiting (MCL) methodology
  to nodal discontinuous Galerkin spectral element methods (DGSEM). The use of  
  Legendre--Gauss--Lobatto (LGL) quadrature endows  collocated  DGSEM space
  discretizations of nonlinear hyperbolic problems with properties that greatly
  simplify the design of  invariant domain preserving high-resolution schemes.
Compared to many other continuous and discontinuous Galerkin method variants,
  a particular  advantage of the LGL spectral operator is the  availability of a
  natural decomposition into a compatible subcell flux discretization. Representing a high-order spatial
  semi-discretization in terms of intermediate states, we perform flux limiting
  in a manner that keeps these states and the results of Runge--Kutta stages
  in convex invariant domains. Additionally, local bounds may be imposed on
  scalar quantities of interest. In contrast to limiting approaches
  based on predictor-corrector algorithms, our MCL procedure for LGL-DGSEM
  yields nonlinear flux approximations that are independent of the time-step size
  and can be further modified to enforce entropy stability.
  To demonstrate the robustness of MCL/DGSEM schemes for
  the compressible Euler equations, we run simulations for
  challenging setups featuring strong shocks, steep density gradients and vortex dominated flows.
}




\maketitle

\section{Introduction}\label{sec:intro}

A wealth of advanced stabilization procedures can be found in the literature on 
second-order finite element discretizations of hyperbolic problems.
In traditional artificial viscosity methods, the amount of nonlinear stabilization
is determined using residual-based shock detectors 
\cite{guermond2011,johnson1990,lv2016,nazarov2013}.
Limiter-based alternatives
adjust numerical fluxes \cite{kuzmin2010a,kuzmin2020monolithic,lohmann2016,guermond2018} or derivatives of 
piecewise-polynomial approximations
 \cite{dobrev2018,moe2017,zhang2010,zhang2012}. The purpose of flux/slope limiting is to enforce sufficient conditions
for positivity preservation, validity of local discrete maximum principles,
and/or entropy stability.
Some limiters are formally applicable
to arbitrary-order finite elements but ``discretization-independent''
\cite{guermond2019} black-box
extensions are far less accurate than piecewise-linear approximations using
the same total number of degrees of freedom \cite{hajduk2020,lohmann2017}.
Accuracy-preserving limiting procedures for high-order finite elements usually rely on the
use of nonoscillatory (WENO) reconstructions \cite{zhang2010,zhang2012,zhu2009,vedral2023},
smoothness indicators \cite{diot2012,krivodonova2004,persson2006}, 
or subcell flux limiting\,/\,shock capturing techniques
\cite{dumbser2014,hajduk2021monolithic,hennemann2021,kuzmin2020a,lohmann2017,vilar2019posteriori,vilar2022}.

The algebraic flux correction (AFC) schemes that we review and modify in the present
paper are based on the methodology that is currently known as \textit{convex limiting}
\cite{guermond2018,hajduk2021monolithic,kuzmin2020monolithic}. The underlying design
philosophy traces its origins
to localized flux-corrected transport (FCT) algorithms for scalar conservation laws
\cite{cotter2016,guermond2017,lohmann2017}. The first extension to nonlinear hyperbolic 
systems was proposed by Guermond et al. \cite{guermond2018}. 
In contrast to Zalesak's multidimensional FCT limiter \cite{zalesak1979} and its edge-based
generalizations to continuous finite element methods for the Euler equations
\cite{kuzmin2010a,kuzmin2012,lohmann2017,lohner2008,selmin1987b}, convex limiting approaches
enforce preservation of local and global bounds by constraining individual fluxes
rather than sums of fluxes. In the explicit case, the local extremum diminishing (LED)
and/or invariant domain preserving
(IDP) properties of flux-limited approximations
are shown using representations in terms of intermediate states that stay
in convex admissible sets \cite{guermond2018,guermond2019}.

All of the aforementioned FCT algorithms belong to the family of AFC
schemes in which the computation of a property-preserving low-order predictor is
followed by an anti-diffusive correction stage.
The monolithic convex limiting (MCL) methodology developed in \cite{kuzmin2020monolithic} 
differs from such fractional-step approaches in that limited anti-diffusive fluxes are
incorporated into the residual of the semi-discrete scheme. The resulting nonlinear
system of ordinary differential equations has a well-defined steady state, and the use of
implicit time integrators is an option. The IDP property of the
explicit version can be shown following the analysis of the low-order (local
Lax--Friedrichs) method in \cite{Guermond2016}. Moreover, the validity of
(semi-)discrete entropy inequalities can be enforced using limiter-based or
dissipation-based fixes \cite{KuHaRu2021,kuzmin2020subcell}.

The first successful extensions of FCT and MCL to high-order finite elements
\cite{anderson2017,hajduk2020,kuzmin2020a,lohmann2017,hajduk2021monolithic} used Bernstein polynomials as local
basis functions. In this context, a key to achieving optimal
accuracy lies in the use of sparse discrete gradient/Laplacian operators
and subcell flux limiting techniques. The discontinuous Galerkin spectral element
 methods (DGSEM) proposed by Pazner \cite{Pazner2020}, Lin et al. \cite{lin2022}, and Rueda-Ram\'irez et al. 
\cite{RUEDARAMIREZ2022} extend subcell convex limiting of FCT type to
Legendre--Gauss--Lobatto (LGL) bases. The underlying low-order method
has the structure of the subcell finite volume scheme employed in
\cite{hennemann2021}. The high-order DGSEM discretization also admits a
natural sparse representation in terms of subcell fluxes between neighbor
nodes. Hence, there is no need for artificial flux reconstructions or decompositions.
Moreover, the mass matrices of collocated  LGL-DGSEM approximations
are diagonal and the discrete gradient/divergence
operators possess summation-by-parts (SBP) properties, which are
needed to achieve entropy stability \cite{gassner2013}. 

As an alternative to the LGL versions \cite{lin2022,Pazner2020,RUEDARAMIREZ2022,RR2021,wu2021}
of high-order FCT algorithms and sophisticated limiters
 for Bernstein finite elements \cite{kuzmin2020subcell,kuzmin2020g,hajduk2021monolithic},
 we introduce a tailor-made LGL-DGSEM counterpart of
 Hajduk's \cite{hajduk2021monolithic}
subcell MCL scheme for conservation laws. 
In fact, the proposed methodology is also applicable to any other spatial semi-discretization that produces sparse discrete gradient operators with SBP properties, such as Gauss-DGSEM discretizations \cite{mateo2022entropy} or general SBP discretizations of nonconservative systems of balance laws \cite{rueda2022flux,rueda2021entropy}.
The flux constraints of the MCL procedure and steady-state
solutions are independent of the time step.
In the context of subcell flux limiting for the
Euler equations of gas dynamics, the density, momentum, and total energy fluxes are
limited sequentially to enforce local bounds for the density, individual
velocity components, and specific total energy. If the pressure
becomes negative, a simple scaling limiter is applied. No
local bounds are imposed on the physical entropy because
the  limiter-based fixes proposed in \cite{KuHaRu2021,kuzmin2020g,kuzmin2020subcell}
guarantee entropy stability under less restrictive constraints. The above
limiting strategy enables us to achieve high resolution without sacrificing
any important properties or using impractically small time steps.

The remainder of this paper is organized as follows. 
In Section \ref{sec:num_methods}, we briefly present the LGL-DGSEM, derive the LGL-DGSEM subcell MCL method, and discuss some of its properties.
In Section \ref{sec:results}, we use the LGL-DGSEM/MCL method to perform challenging simulations of the compressible Euler equations, and present some comparisons with FCT/IDP strategies.
Finally, we draw our conclusions in Section \ref{sec:conclusions}.

\section{Numerical Methods}\label{sec:num_methods}

In this work, we deal with hyperbolic systems of conservation laws of the form
\begin{equation} \label{eq:consLaw}
\bigpartialderiv {\stateL{u}}{t}
+ \Nabla \cdot \blocktensor{f}(\stateL{u})
= \stateL{0}, \quad \mathrm{in} \ \Omega \times \mathbb{R}^+,
\end{equation}
where $\Omega \subseteq \mathbb{R}^{\numDims}$ is a computational domain. The number
of space dimensions is $\numDims\in\{1,2,3\}$.
The vector  $\stateL{u}(\vec x,t)\in \mathbb{R}^{n_{eq}}$ of conserved quantities depends on the space location $\vec x$ and the time instant $t$. 
The flux function $\blocktensor{f}:  \mathcal G\to\mathbb{R}^{\numDims \times n_{eq}}$ depends
on $\stateL{u}:\bar\Omega \times \mathbb{R}_0^+\to \mathcal G$. The set
$\mathcal G\subseteq\mathbb{R}^{n_{eq}}$ is called an invariant domain if
$\mathcal G$ is convex and
$\stateL{u}(\vec x,t)\in\mathcal G$ for all $(\vec x,t)\in \bar\Omega \times \mathbb{R}^+_0$.
System \eqref{eq:consLaw} is equipped with an initial condition, $\stateL{u}(\cdot,0)=\stateL{u}_0$, and suitable boundary conditions on $\partial \Omega$.

For brevity and better readability, we introduce the methods under investigation in the simple context of a one-dimensional ($\numDims=1$) conservation law or system.
All algorithms to be discussed admit straightforward tensor-product extensions to two and three space dimensions, and to curvilinear grids.

\subsection{The Discontinuous Galerkin Spectral Element Method}

Let $\mathcal{T} = \lbrace \Omega^1, \ldots, \Omega^K \rbrace$ be a tessellation of the domain $\Omega$ into $K$ non-overlapping elements. Within each element, we approximate the solution $\stateL{u}$ by a polynomial of degree $N$. A piecewise-polynomial DG approximation $\stateL{u^{\DG}}\approx\stateL{u}$ may be discontinuous at the element interfaces. We seek 
$\stateL{u^{\DG}}$ in the space
$$
\mathcal{V}^N = \lbrace \phi \in L^2(\Omega) \colon \phi \rvert_{\Omega^e} \in \mathcal{P}^N(\Omega^e) \, \forall \Omega^e \in \mathcal{T} \rbrace.
$$
Restricting our attention to a single element $\Omega^e$, we multiply \eqref{eq:consLaw} by an arbitrary polynomial test function $\stateG{\phi}\in (\mathcal{P}^N(\Omega^e))^{n_{eq}}$, integrate the weighted residual over $\Omega^e$, and perform integration by parts to obtain the weak form
\begin{equation} 
\int_{\Omega^e}
\bigpartialderiv {\stateL{u}^{\DG}}{t} \cdot \stateG{\phi} \mathrm{d} \vec{x}
-
\int_{\Omega^e} \blocktensor{f}(\stateL{u}^{\DG}) \cdot \Nabla \stateG{\phi} \mathrm{d} \vec{x}
+
\oint_{\partial \Omega^e} \stateG{\phi} \cdot \numfluxb{f}\mathrm{d} \vec{S}
= \stateL{0}
\end{equation}
of the local conservation law.
Since $\stateL{u}^{\DG}$ is generally
not uniquely defined at the element interfaces, we calculate
$\numfluxb{f}\approx
\blocktensor{f} \cdot \vec{n} \rvert_{\partial \Omega^e}$ using an
approximate Riemann solver that receives two one-sided limits
and returns a numerical flux.

The Legendre--Gauss--Lobatto (LGL) discontinuous Galerkin spectral element method (DGSEM) is a so-called nodal collocation variant of the DG method. It produces discrete gradient/divergence operators that possess summation-by-parts (SBP) properties \cite{gassner_skew_burgers}. The restriction of $\stateL{u}^{\DG}$ to $\Omega^e$ is represented using Lagrange basis functions that are associated with $(N+1)^{\numDims}$ LGL interpolation points. The quadrature rule for numerical integration on $\Omega^e$ uses the LGL collocation nodes on the reference element $\tilde{\Omega}=[-1,1]^{\numDims}$. A mapping $F^e: \tilde{\Omega}\to\Omega^e$ is used for transformations  from the reference space to the physical space ($\vec \xi\mapsto\vec x$ for  $\vec{\xi} \in \tilde{\Omega}$ and
$\vec x=F^e(\xi)\in\Omega^e$). After some manipulations, the evolution equation for the $i$th local degree of freedom
of a one-dimensional LGL-DGSEM discretization of \eqref{eq:consLaw} on $\Omega^e$ can be written as \cite{ranocha2021efficient,RUEDARAMIREZ2022}
\begin{equation} \label{eq:DGSEMstd}
    J \omega_i \dot{\stateL{u}}^{\DG}_i 
    + 
\sum_{k=0}^N \bar{S}_{ik} \stateL{f}_{k}
- \delta_{i0} \numfluxb{f}_{(0,L)}
+ \delta_{iN} \numfluxb{f}_{(N,R)}
= \stateL{0},
\end{equation}
where $J$ denotes the constant determinant of the Jacobian of the mapping from the reference element, $\omega_i$ denotes the reference-space quadrature weight, and $\delta_{ij}$ is the Kronecker delta of the node indices $i$ and $j$. The numerical fluxes $\numfluxb{f}_{(0,L)}$ and  $\numfluxb{f}_{(N,R)}$ are calculated using the inner and outer limits of $\stateL{u}^{\DG}$ on the boundaries of the element $\Omega^e=(x_0^e,x_N^e)$ containing the LGL nodal point $x_i^e$. 
The strong form derivative matrix $\bar{\mat{S}} =( \bar{S}_{ik})_{i,k=0}^N$ admits the representation
$$
\bar{\mat{S}} = \mat{Q} - \mat{B},
$$
where $\mat{B} := \text{diag} (-1, 0, \ldots, 0, 1)$ is the so-called boundary evaluation matrix. 
The entries $Q_{ij} := \omega_i \ell'_j(\xi_i)$ of the weak form derivative matrix $\mat{Q}=(Q_{ij})_{i,k=0}^N$ are defined using the derivatives of the Lagrange basis polynomials $\{ \ell_i \}_{i=0}^N$.

Using the skew-symmetric matrix ${\mat{S}} = 2\mat{Q} - \mat{B} = \mat{Q} - \mat{Q}^T$, whose entries we denote by $S_{ik}$, the discretized volume integral can be expressed in terms of two-point numerical fluxes $\stateL{f}^{*}_{(i,k)}$ \cite{Fisher2013a}. The semi-discrete scheme 
\begin{equation}\label{eq:DGSEMsplit}
    J \omega_i \dot{\stateL{u}}^{\DG}_i 
    + 
\sum_{k=0}^N S_{ik} \stateL{f}^{*}_{(i,k)}
- \delta_{i0} \numfluxb{f}_{(0,L)}
+ \delta_{iN} \numfluxb{f}_{(N,R)}
= \stateL{0}
\end{equation}
 is equivalent to \eqref{eq:DGSEMstd} if the standard average $\stateL{f}^*_{(i,k)} = (\stateL{f}_i+\stateL{f}_k)/2$ is used.
However, additional robustness can be achieved with other choices of the volumetric numerical flux $\stateL{f}^{*}_{(i,k)}$. 
For instance, some two-point approximations to fluxes of the Euler equations guarantee kinetic energy preservation \cite{gassner2016split}, entropy conservation/dissipation \cite{ismail2009affordable,Chandrashekar2013}, pressure equilibrium preservation \cite{shima2021preventing}, or all of these properties together \cite{ranocha2018generalised,ranocha2021preventing}.

All diagonal-norm SBP discretizations of conservation laws (and hence also the LGL-DGSEM considered here) can be written in the so-called \textit{flux-differencing} form \cite{Fisher2013a}
\begin{equation} \label{eq:DGSEM_fluxDiff}
J \dot{\stateL{u}}^{\DG}_{i} =
\frac{1}{\omega_i}
\left(
  \numfluxb{f}^{\DG}_{(i-1,i)}
- \numfluxb{f}^{\DG}_{(i,i+1)}
\right)
,
~~~~~~
\forall i=0, \ldots, N,
\end{equation}
where the indices $-1$ and $N+1$ refer to the outer states. The  symmetric and consistent fluxes $\numfluxb{f}^{\DG}_{(i,\ii)}=\numfluxb{f}^{\DG}_{(\ii,i)}$ are defined by \cite{Fisher2013a,RUEDARAMIREZ2022}
\begin{align}
\numfluxb{f}^{\DG}_{(-1,0)}  &= \numfluxb{f}_{(0,L)},
\label{eq:leftFlux}
\\
\numfluxb{f}^{\DG}_{(i,i+1)} &= \sum_{l=0}^i \sum_{k=0}^N S_{lk} \stateL{f}^{*}_{(l,k)}, & i=0, \ldots, N-1,
\label{eq:inteFlux}\\
\numfluxb{f}^{\DG}_{(N,N+1)} &= \numfluxb{f}_{(N,R)}. \label{eq:rightFlux}
\end{align}
Note that the flux $\numfluxb{f}^{\DG}_{(i,\ii)}$ is multiplied
 by the one-dimensional unit normal
$n_{(i,\ii)}\in\{-1,1\}$
 in \eqref{eq:DGSEM_fluxDiff}. The normal fluxes $n_{(i,\ii)}
 \numfluxb{f}^{\DG}_{(i,\ii)}$ are anti-symmetric, that is,
 $n_{(i,\ii)}\numfluxb{f}^{\DG}_{(i,\ii)}=-n_{(\ii , i)}
 \numfluxb{f}^{\DG}_{(\ii ,i)}$. Hence, \eqref{eq:DGSEMsplit} has local (subcell-level) conservation properties, as required by the Lax--Wendroff theorem~\cite{laxwendroff}.

\begin{remark}
  Let $\Delta x_i=J\omega_i$. Then
  \eqref{eq:DGSEM_fluxDiff} corresponds to the subcell finite volume scheme
  $$
\dot{\stateL{u}}^{\DG}_{i} =-\frac{\numfluxb{f}_{i+1/2}^{\DG}-\numfluxb{f}_{i-1/2}^{\DG}}{\Delta x_i},
~~~~~~
\forall i=0, \ldots, N,
$$
where $\numfluxb{f}_{i+1/2}^{\DG}=\numfluxb{f}^{\DG}_{(i,i+1)}$
and $\numfluxb{f}_{i-1/2}^{\DG}=\numfluxb{f}^{\DG}_{(i-1,i)}$. We adopt the
two-subscript notation because it is better suited for flux-based
finite element discretizations.
  \end{remark}

\begin{remark}\label{remark:fluxdiff_other_dg}
Since the LGL-DGSEM is a diagonal-norm SBP operator, its representation in the flux-differencing form  \eqref{eq:DGSEM_fluxDiff}-\eqref{eq:rightFlux} is readily available.
Other DG approximations with dense mass matrices need the application of a sparsification operator to recover the flux-differencing form.
Examples of decompositions into subcell fluxes can be found, e.g., in \cite{hajduk2021monolithic,kuzmin2020subcell,vilar2022}. For a DG method using Bernstein polynomials of degree $N>1$ as local basis functions, $n_{eq}$ sparse linear systems of size $(N+1) \times (N+1)$ need to be solved for each element in each Runge--Kutta stage \cite{hajduk2021monolithic,kuzmin2020subcell}.
Vilar \cite{vilar2019posteriori} showed that it is possible to obtain a flux-differencing formula for any (modal or nodal) representation  of the DG solution if one expresses the test function as a combination of so-called \textit{subresolution} basis functions and exploits existing relationships to the histopolation theory. An adaptation to unstructured triangular grids was proposed by Vilar and Abgrall \cite{vilar2022}, who parametrized $\stateL{u}^{\DG}$ in terms of subcell averages that satisfy a two-dimensional version of $\eqref{eq:DGSEM_fluxDiff}$. The calculation of subcell fluxes involves solving small linear systems with sparse graph Laplacians again.
\end{remark}

\begin{remark}
Mateo-Gabín et al.~\cite{mateo2022entropy} showed that (both standard and split-form versions of) the Legendre--Gauss DGSEM scheme can also be written in the flux-differencing form with explicit staggered fluxes and a diagonal mass matrix.
As a result, most of the algorithms to be presented in this paper are applicable to the Legendre--Gauss DGSEM. However, the treatment of inter-element fluxes and projection operators requires additional analysis and, possibly, appropriate modifications.
\end{remark}

\subsection{Monolithic Convex Limiting} \label{sec:mcl}

The monolithic convex limiting (MCL) methodology \cite{hajduk2021monolithic,kuzmin2020monolithic,kuzmin2020subcell} is a subcell flux correction procedure that combines a high-order baseline discretization with a compatible and invariant domain preserving low-order scheme. The validity of physical and numerical admissibility conditions is enforced using a representation in terms of intermediate states (similarly to the predictor-corrector approaches proposed in \cite{guermond2018,guermond2019,lin2022,Pazner2020,RUEDARAMIREZ2022}). Flux limiters for semi-discrete MCL schemes can be designed to enforce entropy stability conditions in addition to local and/or global maximum principles \cite{kuzmin2022limiter,kuzmin2020g}. To minimize the levels of low-order numerical dissipation, localized subcell limiting procedures are used for high-order finite elements
\cite{hajduk2021monolithic,kuzmin2020subcell,kuzmin2020g}. Moreover, sequential MCL algorithms for systems support the possibility of using individually chosen correction factors for different conserved or derived quantities \cite{hajduk2021monolithic,kuzmin2020monolithic}.

\subsubsection{Low-Order Invariant Domain Preserving Scheme}

As explained in \cite{hennemann2021}, one can obtain a low-order finite volume scheme that is compatible with the LGL-DGSEM discretization by interpreting the nodal values of the DGSEM scheme as mean values of the subcells.
Let
\begin{equation} \label{eq:FV_scheme}
J \dot{\stateL{u}}^{\FV}_{i} =
\frac{1}{\omega_i}
\left(
  \numfluxb{f}^{\FV}_{(i-1,i)}
- \numfluxb{f}^{\FV}_{(i,i+1)}
\right)
,
~~~~~~
\forall i=0, \ldots, N,
\end{equation}
where $\numfluxb{f}^{\FV}_{(i,\ii)}$ is a low order numerical approximation to the flux between nodes $i$ and $\ii$. Such a subcell FV scheme exhibits the same structure as \eqref{eq:DGSEM_fluxDiff}. Therefore, the two schemes are compatible and can be hybridized.

It has been shown that \eqref{eq:FV_scheme} is invariant domain preserving (IDP) for the first-order Rusanov (also known as local Lax-Friedrichs, LLF) fluxes \cite{Pazner2020,RUEDARAMIREZ2022}
\begin{equation}
\label{eq:rusonaov_llf}
    \numfluxb{f}^{\FV}_{(i,\ii)} = \frac{\stateL{f}_{i} + \stateL{f}_{\ii}}{2} - (i - \ii)\,\frac{\lambda^{\max}_{(i,\ii)}}{2}(\stateL{u}_i - \stateL{u}_{\ii}),
\end{equation}
where $\lambda^{\max}_{(i,\ii)} = \lambda^{\max}_{(\ii,i)}>0$ is an upper bound for the maximum wave speed of the Riemann problem with the initial states
$\stateL{u}_i$ and $\stateL{u}_{\ii}$. Estimation of this speed is
addressed, e.g., in \cite{guermond2016fast}.
In our notation, the presence of $(i - \ii)$ in the dissipative part of \eqref{eq:rusonaov_llf} ensures that the flux $\numfluxb{f}^{\FV}_{(i,\ii)}=\numfluxb{f}^{\FV}_{(\ii,i)}$ is symmetric. The multiplication by the unit normal $n_{(i,\ii)}=-n_{(\ii , i)}$ makes it anti-symmetric.

Inserting the Rusanov fluxes \eqref{eq:rusonaov_llf} into \eqref{eq:FV_scheme} and using the forward Euler method for time integration
yields a fully discrete version of the low-order scheme. Following
Guermond and Popov \cite{Guermond2016}, we write it in the form
\begin{align}
\stateL{u}_i^{\FV,n+1} =& \left(1 - \frac{\Delta t}{J\,\omega_i}\left(\lambda^{\max}_{(i,i-1)}+\lambda^{\max}_{(i,i+1)}\right)\right)\stateL{u}_i^n \nonumber
\\
& + \frac{\Delta t}{J\,\omega_i}\,\lambda^{\max}_{(i,i-1)}\, \overline{\stateL{u}}_{(i,i-1)}^n + \frac{\Delta t}{J\,\omega_i}\,\lambda^{\max}_{(i,i+1)}\, \overline{\stateL{u}}_{(i,i+1)}^n\label{barstateformLO}
\end{align}
using the auxiliary \textit{bar states}
\begin{equation}\label{def:barstates}
    \overline{\stateL{u}}_{(i,\ii)} := \frac{\stateL{u}_{i} + \stateL{u}_{\ii}}{2} - (i -\ii)\,\frac{\stateL{f}_{i} - \stateL{f}_{\ii}}{2\,\lambda^{\max}_{(i,\ii)}}.
\end{equation}
If the time-step size $\Delta t$ satisfies the CFL condition
\begin{equation}
\Delta t \leq \frac{J\,\omega_i}{\lambda_{(i,i-1)}^{\max}+\lambda_{(i,i+1)}^{\max}},
\end{equation}
then the result $\stateL{u}_i^{\FV,n+1}$ of the explicit update \eqref{barstateformLO}
is a property-preserving convex combination of the states $\stateL{u}_i^n$
and $\overline{\stateL{u}}_{(i,i\pm1)}^n$.
We discuss the time-step restriction in Section~\ref{sec:timestep}.

The most important property of the so-defined LLF bar states is that they preserve all \textit{convex invariants} of initial value problems for hyperbolic systems, as shown in \cite{Guermond2016} in the context of a continuous (multi-)linear finite element discretization. In fact,  $\overline{\stateL{u}}_{(i,\ii)}$ defined by \eqref{def:barstates} is the intermediate state of the HLL approximate Riemann solver \cite{hll1983}. Positivity preservation and the validity of entropy conditions can be deduced from this interpretation.

\begin{remark}
The bar states \eqref{def:barstates} of the (semi-discrete or fully discrete) low-order LLF scheme
  are symmetric in the sense that $\overline{\stateL{u}}_{(i,\ii)} = \overline{\stateL{u}}_{(\ii,i)}$ for any pair of adjacent nodes with local indices $i \in \{0,\dots,N\}$ and $\ii \in \lbrace i-1, i+1 \rbrace$.
\end{remark}

\begin{remark}
  To obtain compatible low-order IDP schemes for general high-order DG methods, it is necessary to first replace the discrete gradient/divergence operators with sparse approximations and then apply low-order dissipation (e.g., using a sparse graph Laplacian operator as in \cite{hajduk2020,hajduk2021monolithic}).
For our LGL-DGSEM scheme and, in general, for all diagonal norm SBP operators, the low-order IDP scheme \eqref{eq:FV_scheme} is readily available and compatible with the flux-differencing form \eqref{eq:DGSEM_fluxDiff} of \eqref{eq:DGSEMsplit}.
\end{remark}

\subsubsection{Limiting Procedure}

To enforce relevant inequality constraints, we replace 
\eqref{eq:DGSEM_fluxDiff} with (cf.~\cite{hajduk2021monolithic,kuzmin2020subcell})
\begin{equation} \label{eq:blendedScheme}
J \dot{\stateL{u}}_{i} =
\frac{1}{\omega_i}
\left(
  \numfluxb{f}_{(i-1,i)}
- \numfluxb{f}_{(i,i+1)}
\right)
,
~~~~~~
\forall i=0, \ldots, N.
\end{equation}
In the simplest case, the hybrid subcell fluxes $\numfluxb{f}_{(i,\ii)}\in\mathbb{R}^{n_{eq}}$ are given by
\begin{equation} \label{flux:alpha}
    \numfluxb{f}_{(i,\ii)} = \stateG{\alpha}_{(i,\ii)}  \circ \numfluxb{f}^{\DG}_{(i,\ii)} +  
    (\stateL{1}-\stateG{\alpha}_{(i,\ii)}) \circ \numfluxb{f}^{\FV}_{(i,\ii)},
\end{equation}
where $\circ$ denotes the Hadamard (or component-wise) product. The  scalar-valued components of $\stateG{\alpha_{(i,\ii)}}\in\mathbb{R}^{n_{eq}}$ are weights that attain values between 0 and 1. The high-order DG method \eqref{eq:DGSEM_fluxDiff} and the low-order FV scheme \eqref{eq:FV_scheme} can be recovered using $\stateG{\alpha}_{(i,\ii)}=\stateL{1}$ and  $\stateG{\alpha}_{(i,\ii)}=\stateL{0}$, respectively.

As detailed in the next section, the computation of $\stateG{\alpha_{(i,\ii)}}$
 might not be numerically well posed. Therefore, it should be avoided in practical implementations if there is a direct way to calculate the fluxes $\numfluxb{f}_{(i,\ii)}$.

\begin{remark}
  If  $\numfluxb{f}^{\DG}_{(i,\ii)}=\numfluxb{f}^{\FV}_{(i,\ii)}$ for
  $j\in\{-1,N+1\}$, then the DG and FV schemes use the same two-point
  flux approximation on the boundaries of $\Omega^e=(x_0^e,x_N^e)$. In this case, formula
  \eqref{flux:alpha} will produce $\numfluxb{f}_{(i,\ii)}=
\numfluxb{f}^{\DG}_{(i,\ii)}=\numfluxb{f}^{\FV}_{(i,\ii)} $
  for any choice of $\stateG{\alpha_{(i,\ii)}}$. This desirable property of boundary fluxes is specific to the LGL-DGSEM discretization because it includes the boundary nodes. It leads to a very local implementation of the limiting procedure.
\end{remark}

Since the low-order component of  $\numfluxb{f}_{(i,\ii)}$ is provably IDP,
the purpose of subcell limiting is to constrain the anti-diffusive components
\begin{equation}
\Delta \numfluxb{f}_{(i,\ii)} = (i -\ii)\,\left(\numfluxb{f}^{\DG}_{(i,\ii)} - \numfluxb{f}^{\FV}_{(i,\ii)} \right)    
\end{equation}
of the fluxes $\numfluxb{f}^{\DG}_{(i-1,i)} = \numfluxb{f}^{\FV}_{(i-1,i)} + \Delta \numfluxb{f}_{(i-1,i)}$ and 
$\numfluxb{f}^{\DG}_{(i,i+1)} = \numfluxb{f}^{\FV}_{(i,i+1)} - \Delta \numfluxb{f}_{(i,i+1)}$.

\begin{remark}
The so-called anti-diffusive flux $\Delta \numfluxb{f}_{(i,\ii)} = - \Delta \numfluxb{f}_{(\ii,i)}$ is anti-symmetric because the fluxes $\numfluxb{f}^{\DG}_{(i,\ii)}=\numfluxb{f}^{\DG}_{(\ii,i)}$ and $\numfluxb{f}^{\FV}_{(i,\ii)}=\numfluxb{f}^{\FV}_{(\ii,i)}$ are symmetric.
\end{remark}

With this notation, the high-order update can be written as 
\begin{equation} \label{eq:flux-diff2}
J\omega_i\frac{\stateL{u}_i^{\DG,n+1} - \stateL{u}_i^{n}}{\Delta t} = \numfluxb{f}^{\FV,n}_{(i-1,i)}
- \numfluxb{f}^{\FV,n}_{(i,i+1)} + \Delta \numfluxb{f}_{(i,i-1)}^n + \Delta \numfluxb{f}_{(i,i+1)}^n.
\end{equation}

Using the representation of the flux difference $\numfluxb{f}^{\FV}_{(i-1,i)}
- \numfluxb{f}^{\FV}_{(i,i+1)}$ in terms of the bar states defined by \eqref{def:barstates}, we find that
\begin{equation}\label{def:DGunlim}
\begin{split}
\stateL{u}_i^{\DG,n+1} &= \left(1 - \frac{\Delta t}{J\,\omega_i}\left(\lambda^{\max}_{(i,i-1)}+\lambda^{\max}_{(i,i+1)}\right)\right)\stateL{u}_i^n\\
&+ \frac{\Delta t}{J\,\omega_i}\,\lambda^{\max}_{(i,i-1)}\, \overline{\stateL{u}}_{(i,i-1)}^n + \frac{\Delta t}{J\,\omega_i}\,\lambda^{\max}_{(i,i+1)}\, \overline{\stateL{u}}_{(i,i+1)}^n\\
&+ \frac{\Delta t}{J\,\omega_i}\left(\Delta \numfluxb{f}_{(i,i-1)}^n + \Delta \numfluxb{f}_{(i,i+1)}^n\right).
\end{split}
\end{equation}

We can now define the bar state of the high-order method as 
\begin{equation}\label{def:barsttateHO}
    \overline{\stateL{u}}^{\DG}_{(i,\ii)} = \overline{\stateL{u}}_{(i,\ii)} + \frac{\Delta \numfluxb{f}_{(i,\ii)}}{\lambda^{\max}_{(i,\ii)}},
\end{equation}
and cast \eqref{def:DGunlim} into the bar state form 
\begin{equation} \label{eq:convexSumBarStates}
\begin{split}
\stateL{u}_i^{\DG,n+1} &=\left(1 - \frac{\Delta t}{J\,\omega_i}\left(\lambda^{\max}_{(i,i-1)}+\lambda^{\max}_{(i,i+1)}\right)\right) \stateL{u}_i^n
\\
&+ \frac{\Delta t}{J\,\omega_i}\,\lambda^{\max}_{(i,i-1)}\, \overline{\stateL{u}}_{(i,i-1)}^{\DG,n} + \frac{\Delta t}{J\,\omega_i}\,\lambda^{\max}_{(i,i+1)}\, \overline{\stateL{u}}_{(i,i+1)}^{\DG,n}
\end{split}
\end{equation}
which has the same structure as \eqref{barstateformLO}. 

Following the derivation of MCL schemes for Lagrange and Bernstein finite elements \cite{kuzmin2020monolithic,kuzmin2020subcell,hajduk2021monolithic}, we replace the DG bar states defined by
\eqref{def:barsttateHO} with 
\begin{equation}\label{def:barsttateMCL}
    \overline{\stateL{u}}^{Lim}_{(i,\ii)} = \overline{\stateL{u}}_{(i,\ii)} + \frac{\Delta \numfluxb{f}^{Lim}_{(i,\ii)}}{\lambda^{\max}_{(i,\ii)}},
\end{equation}
where $\Delta \numfluxb{f}^{Lim}_{(i,\ii)}=
(i -\ii)\,\left(\numfluxb{f}_{(i,\ii)} - \numfluxb{f}^{\FV}_{(i,\ii)} \right)$
is a limited approximation to $\Delta\numfluxb{f}_{(i,\ii)}$.
In contrast to the low-order component $\overline{\stateL{u}}_{(i,j)} = \overline{\stateL{u}}_{(j,i)}$, the limited bar state \eqref{def:barsttateMCL} is generally not symmetric due to the skew-symmetry of $\Delta \numfluxb{f}^{Lim}_{(i,\ii)}=-\Delta \numfluxb{f}^{Lim}_{(\ii,i)}$. 

The MCL 
bar states $\overline{\stateL{u}}^{Lim}_{(i,\ii)}$ should satisfy the
same inequality constraints as $\overline{\stateL{u}}_{(i,\ii)}$ and
stay as close as possible to the high-order target $\overline{\stateL{u}}^{\DG}_{(i,\ii)}$. The forward Euler time discretization should be replaced with
a high-order Runge--Kutta method. For SSP-RK schemes with forward Euler
stages, the IDP property can be shown in the same way as for
the low-order scheme  \cite{kuzmin2020monolithic,hajduk2021monolithic}.
A general Runge--Kutta method may require flux limiting
in time \cite{kuzmin2022timelim,quezada2022}.

The convex limiting techniques employed in 
\cite{guermond2018,guermond2019,lin2022,Pazner2020,RUEDARAMIREZ2022}
differ from MCL in that they split the computation of $\stateL{u}_i^{\DG,n+1}$
into a low-order IDP update and an anti-diffusive correction stage. This
predictor-corrector strategy is also used in older FCT-type
algorithms for finite element discretizations of hyperbolic systems
\cite{kuzmin2010a,kuzmin2012,lohmann2017,lohner2008,selmin1987b}.
In contrast to MCL, the resulting schemes have no semi-discrete
counterparts. Moreover, the bounds of the limiting constraints depend
on the time step. Depending on the application, this
peculiarity of FCT/IDP approaches may be an advantage or a disadvantage.

We will now describe the computation of the limited anti-diffusive fluxes
$\Delta \numfluxb{f}^{Lim}_{(i,\ii)}$
for the MCL version with generic bounds. Appropriate definitions of
the bounds are discussed in Section \ref{sec:bounds}.

\paragraph{Limiter for conservative quantities} \label{sec:consLimiter}
The simplest limiting strategy for systems is to treat each equation 
as a scalar conservation law and to limit the anti-diffusive fluxes of
each conserved variable individually. Let $\rho$ be a scalar component
of $\stateL{u}$. We denote by  $\overline{\rho}_{(i,j)}$ and $\Delta \numflux{f}^{\rho,Lim}_{(i,\ii)}$ the corresponding components of $\overline{\stateL{u}}_{(i,j)}$ and $\Delta \numfluxb{f}^{Lim}_{(i,\ii)}$, respectively. 
To keep
$\overline{\rho}_{(i,j)}^{Lim}$ and $\overline{\rho}_{(j,i)}^{Lim}$
in the range $[\rho^{\min}_i,\rho^{\max}_i]$, we impose the
inequality constraints
\begin{align} \label{eq:bounsCons}
    \rho^{\min}_i 
    \le \overline{\rho}_{(i,j)} + \frac{\Delta \numflux{f}^{\rho,Lim}_{(i,\ii)}}{\lambda^{\max}_{(i,\ii)}}
    \le
    \rho^{\max}_i,
    ~~~~
    \rho^{\min}_j 
    \le \overline{\rho}_{(i,j)} - \frac{\Delta \numflux{f}^{\rho,Lim}_{(i,\ii)}}{\lambda^{\max}_{(i,\ii)}}
    \le
    \rho^{\max}_j.
\end{align}
A positive/negative anti-diffusive flux
$\Delta \numflux{f}^{\rho,Lim}_{(i,\ii)}$
may violate the upper/lower bound for $\overline{\rho}_{(i,j)}^{Lim}$ or
the lower/upper bound for $\overline{\rho}_{(j,i)}^{Lim}$.  Introducing
\begin{align*}
   \Delta \numflux{f}^{\rho,+}_{(i,\ii)} &=
   \lambda^{\max}_{(i,\ii)} \min \lbrace \rho^{\max}_i - \overline{\rho}_{(i,j)}, \overline{\rho}_{(i,j)} - \rho^{\min}_j  \rbrace , \\
      \Delta \numflux{f}^{\rho,-}_{(i,\ii)} &=
   \lambda^{\max}_{(i,\ii)} \max \lbrace \rho^{\min}_i - \overline{\rho}_{(i,j)}, \overline{\rho}_{(i,j)} - \rho^{\max}_j  \rbrace,
\end{align*}
we define \cite{kuzmin2020monolithic,hajduk2021monolithic}
\begin{equation} \label{eq:antiDiffCons}
   \Delta \numflux{f}^{\rho,Lim}_{(i,\ii)} =
   \begin{cases}
   \min \lbrace \Delta \numflux{f}^{\rho}_{(i,\ii)},  \Delta \numflux{f}^{\rho,+}_{(i,\ii)}\rbrace  & \mathrm{if}~ \Delta \numflux{f}^{\rho}_{(i,\ii)} \ge 0, \\
   \max \lbrace \Delta \numflux{f}^{\rho}_{(i,\ii)},  \Delta \numflux{f}^{\rho,-}_{(i,\ii)} \rbrace
   & \mathrm{otherwise.}
   \end{cases}
\end{equation}
It is easy to verify that conditions \eqref{eq:bounsCons} are met
for this choice of $\Delta \numflux{f}^{\rho,Lim}_{(i,\ii)}$. Moreover, there
exists $\alpha^{\rho}_{(i,j)}\in[0,1]$ such that
$\Delta \numflux{f}^{\rho,Lim}_{(i,\ii)}=\alpha^{\rho}_{(i,j)}
\Delta \numflux{f}^{\rho}_{(i,\ii)}$ and
\begin{equation} 
  \numflux{f}^{\rho}_{(i,\ii)} =\numflux{f}^{\rho,\FV}_{(i,\ii)}
  + \frac{\Delta \numflux{f}^{\rho,Lim}_{(i,\ii)}}{i-j}=
          {\alpha}^{\rho}_{(i,\ii)}  \numflux{f}^{\rho,\DG}_{(i,\ii)} +  
    ({1}-{\alpha}^{\rho}_{(i,\ii)})  \numflux{f}^{\rho,\FV}_{(i,\ii)}
\end{equation}
is a convex combination of the FV and DG fluxes. Kuzmin \cite{kuzmin2020monolithic} noticed that the computation of  $\alpha^{\rho}_{(i,j)}=
\Delta \numflux{f}^{\rho,Lim}_{(i,\ii)}/\Delta \numflux{f}^{\rho}_{(i,\ii)}$
is unnecessary and
numerically ill posed in the case of a small nonvanishing
denominator. The direct computation of the limited anti-diffusive flux
\eqref{eq:antiDiffCons} is therefore preferable in practice.

\paragraph{Sequential limiter for ``primitive'' quantities}
In some situations, we are interested in imposing bounds on scalar quantities that are not included in the state vector $\stateL{u}$.
If the quantity of interest represents the ratio of two conservative variables, we can use the sequential limiting approach proposed in \cite{dobrev2018,kuzmin2020monolithic}.
For instance, components of the velocity field, $\vec{v} = (\rho \vec{v}) / \rho$, and the total specific energy, $E = (\rho E) / \rho$, of the Euler equations of gas dynamics (see Appendix \ref{app:euler}) might belong to the set of
control variables.

Let $\rho$ and $\rho\phi$ be generic conservative variables. To ensure that
$\rho_i\in [\rho^{\min}_i,\rho^{\max}_i]$ and
$\phi_i\in [\phi^{\min}_i,\phi^{\max}_i]$, the first stage of a
sequential MCL algorithm
\cite{hajduk2021monolithic,kuzmin2020monolithic} limits
$\Delta \numflux{f}^{\rho}_{(i,j)}$  using \eqref{eq:antiDiffCons}. The 
second stage limits $\Delta \numflux{f}^{\rho\phi}_{(i,j)}$ using a
discrete version of the product rule $(\rho\phi)'=\rho'\phi+\phi'\rho$.
The bar states 
\begin{equation}
    \overline{\phi}_{(i,j)} = \frac{\overline{(\rho \phi)}_{(i,j)} + \overline{(\rho \phi)}_{(j,i)}}{\overline{\rho}_{(i,j)} + \overline{\rho}_{(j,i)}}
    \eqLGL
    \frac{\overline{(\rho \phi)}_{(i,j)}}{\overline{\rho}_{(i,j)}}
\end{equation}
of the $\phi$ variable are symmetric in the LGL-DGSEM version. This
property is a further advantage compared to
Bernstein-basis DG methods \cite{hajduk2021monolithic}.

The inequality constraints to be enforced in the second stage are given by
\begin{align} \label{eq:bounsSeqG}
    \overline{\rho}^{Lim}_{(i,j)} \phi^{\min}_i
    \le 
    \overline{\rho}^{Lim}_{(i,j)} \overline{\phi}_{(i,j)} + \frac{\Delta \numflux{g}^{\phi,Lim}_{(i,\ii)}}
    {\lambda^{\max}_{(i,\ii)}}
    \le
    \overline{\rho}^{Lim}_{(i,j)} \phi^{\max}_i,
\end{align}
where $\Delta \numflux{g}^{\phi,Lim}_{(i,\ii)}$ is a limited approximation to
\begin{equation} 
    \Delta \numflux{g}^{\phi}_{(i,j)}
    =\Delta \numflux{f}^{\rho \phi}_{(i,j)} - \lambda^{\max}_{(i,\ii)} 
    \left(\overline{\rho}^{Lim}_{(i,j)} \overline{\phi}_{(i,j)} - \overline{(\rho \phi)}_{(i,j)} \right).
\end{equation}
It is easy to verify that conditions \eqref{eq:bounsSeqG} are equivalent to
\begin{align} \label{eq:bounsSeq}
    \overline{\rho}^{Lim}_{(i,j)} \phi^{\min}_i
    \le 
    \overline{(\rho \phi)}_{(i,j)} + \frac{\Delta \numflux{f}^{\rho\phi,Lim}_{(i,\ii)}}{\lambda^{\max}_{(i,\ii)}}
    \le
    \overline{\rho}^{Lim}_{(i,j)} \phi^{\max}_i,
\end{align}
where
\begin{equation}\label{def:frophilim}
  \Delta \numflux{f}^{\rho \phi,Lim}_{(i,j)} =
  \Delta \numflux{g}^{\phi,Lim}_{(i,j)}+
  \lambda^{\max}_{(i,\ii)} 
    \left(\overline{\rho}^{Lim}_{(i,j)} \overline{\phi}_{(i,j)} - \overline{(\rho \phi)}_{(i,j)} \right).
\end{equation}
We use the bounding fluxes
\begin{align*}
  \Delta \numflux{g}^{\phi,+}_{(i,\ii)}&= \lambda^{\max}_{(i,\ii)} \overline{\rho}^{Lim}_{(i,j)} \min \lbrace \phi^{\max}_i - \overline{\phi}_{(i,j)}, \overline{\phi}_{(i,j)} - \phi^{\min}_j  \rbrace,\\ 
   \Delta \numflux{g}^{\phi,-}_{(i,\ii)}&= \lambda^{\max}_{(i,\ii)} \overline{\rho}^{Lim}_{(i,j)} \max \lbrace \phi^{\min}_i - \overline{\phi}_{(i,j)}, \overline{\phi}_{(i,j)} - \phi^{\max}_j  \rbrace
\end{align*}
to define
\begin{equation} \label{eq:antiDiffSeq}
   \Delta \numflux{g}^{\phi,Lim}_{(i,\ii)} =
   \begin{cases}
   \min \lbrace \Delta \numflux{g}^{\phi}_{(i,\ii)}, 
   \Delta \numflux{g}^{\phi,+}_{(i,\ii)}
   \rbrace & \mathrm{if}~ \Delta \numflux{g}^{\phi}_{(i,\ii)} \ge 0, \\
   \max \lbrace \Delta \numflux{g}^{\phi}_{(i,\ii)}, 
  \Delta \numflux{g}^{\phi,-}_{(i,\ii)}
   \rbrace
   & \mathrm{otherwise.}
   \end{cases}
\end{equation}
This definition, which is similar to \eqref{eq:antiDiffCons}, 
guarantees the  validity of \eqref{eq:bounsSeq} and of the
corresponding constraints for the flux-corrected bar state
$\overline{(\rho \phi)}_{(j,i)}^{Lim}$. The limited anti-diffusive
flux $\Delta \numflux{f}^{\rho \phi,Lim}_{(i,j)}$
is calculated using formula \eqref{def:frophilim}.

\begin{remark}
  By definition \eqref{eq:antiDiffSeq}, there
  exists an effective limiting factor
  $\alpha^{\phi}\in[0,1]$ such that
  $\Delta \numflux{g}^{\phi,Lim}_{(i,\ii)}=\alpha^{\phi}
  \Delta \numflux{g}^{\phi}_{(i,\ii)}$. However, the
  value of $\alpha^{\rho\phi}$ corresponding to the identity
  $\Delta \numflux{f}^{\rho\phi,Lim}_{(i,\ii)}=\alpha^{\rho\phi}
  \Delta \numflux{f}^{\rho\phi}_{(i,\ii)}$ does not necessarily
   lie in the range $[0,1]$.
\end{remark}

\paragraph{Pressure limiter} \label{sec:pressureLimiter}
When solving the compressible Euler equations of gas dynamics (Appendix~\ref{app:euler}), we require the pressure and internal energy to be non-negative at all times. Positivity preservation is guaranteed if the limited bar states satisfy
\begin{equation}\label{eq:egeq0}
  \overline\rho_{(i, \ii)}^{Lim} \overline{(\rho E)}_{(i, \ii)}^{Lim} - \frac{\norm{\overline{(\rho \vec{v})}_{(i, \ii)}^{Lim}}^2}{2} \geq 0.
\end{equation}

To enforce \eqref{eq:egeq0}, we apply a synchronized limiting factor
$\alpha_{(i,j)}^p\in[0,1]$ to all components of $\Delta \numfluxb{f}^{Lim}_{(i,j)}$.
  The limited bar states become
\begin{equation}
    \overline{\stateL{u}}^{Lim}_{(i,j)} = \overline{\stateL{u}}_{(i,j)} + \frac{\alpha^p_{(i,j)} \Delta \numfluxb{f}^{Lim}_{(i,j)}}{\lambda^{\max}_{(i,j)}},
        ~~~~
        \overline{\stateL{u}}^{Lim}_{(j,i)} = \overline{\stateL{u}}_{(i,j)} - \frac{\alpha^p_{(i,j)} \Delta \numfluxb{f}^{Lim}_{(i,j)}}{\lambda^{\max}_{(i,j)}}
\end{equation}
and the prelimited anti-diffusive fluxes $\Delta \numfluxb{f}^{Lim}_{(i,j)}$ are replaced with $\alpha^p_{(i,j)} \Delta \numfluxb{f}^{Lim}_{(i,j)}$.

Dropping the superscript $p$ for better readability and introducing the
 scaled bar states $\overline{\stateL{w}}_{(i,j)}:=\lambda^{\max}_{(i,j)}\overline{\stateL{u}}_{(i,j)}$, we translate \eqref{eq:egeq0} into the quadratic inequalities
\begin{align} \label{eq:quadIneqPressure}
  A_{(i,j)}\alpha_{(i,j)}^2\pm B_{(i,j)}\alpha  \le Q_{(i,j)},
\end{align}
where
\begin{align*} \label{eq:quadIneqPressure}
 A_{(i,j)} &=\frac{\|\Delta \vec{\numflux{f}}^{\rho \vec{v}, Lim}_{(i,j)}\|^2}{2} - \Delta \numflux{f}^{\rho, Lim}_{(i,j)} \Delta \numfluxb{f}^{\rho E, Lim}_{(i,j)},\\
 B_{(i,j)} &= \vec{\overline{w}}^{\rho \vec{v}}_{(i,j)} \cdot \Delta \vec{\numflux{f}}^{\rho \vec{v}, Lim}_{(i,j)} - \overline w_{(i,j)}^{\rho} \Delta \numflux{f}^{\rho E, Lim}_{(i,j)} - \overline{w}_{(i,j)}^{\rho E} \Delta \numflux{f}^{\rho, Lim}_{(i,j)},\\
 Q_{(i,j)} &= \overline{w}_{(i,j)}^{\rho}\overline{w}_{(i,j)}^{\rho E} + \frac{\| \vec{\overline{w}}_{(i,j)}^{\rho \vec{v}} \|^2}{2}.
\end{align*}
Following Kuzmin \cite{kuzmin2020monolithic}, we notice that
$\alpha^2\le\alpha$ for
$\alpha\in[0,1]$. Therefore,
\eqref{eq:quadIneqPressure} holds under the linear sufficient condition
$P_{(i,j)}\alpha \le Q_{(i,j)}$, where
\begin{equation} \label{eq:pressureLimiterExact}
    P_{(i,j)} =  \max \lbrace 0,   A_{(i,j)}\rbrace+
    \rvert B_{(i,j)}\rvert. 
\end{equation}
We conclude that the pressure fix can be performed using
\begin{equation} \label{eq:pressureLimiter}
    \alpha_{(i,j)} =
    \begin{cases}
    \frac{Q_{(i,j)}}{P_{(i,j)}} & \mathrm{if}~P_{(i,j)} > Q_{(i,j)}, \\
    1 & \mathrm{otherwise}.
    \end{cases}
\end{equation}

This definition exploits the property that the bar states of low-order
LGL-DGSEM are symmetric. The general formula for $\alpha_{(i,j)}$ is more
involved \cite{kuzmin2020monolithic}.

To ensure continuous dependence of the limited fluxes
$\alpha^p_{(i,j)} \Delta \numfluxb{f}_{(i,j)}$ on the data,
one may
replace $P_{(i,j)}$ with the upper bound \cite{kuzmin2020monolithic}
\begin{align}
   P^{\max}_{(i,j)} &=  \max \lbrace 0,   A_{(i,j)}\rbrace+
    \rvert \vec{\overline{{w}}}^{\rho \vec{v}}_{(i,j)} \rvert \cdot \rvert \Delta \vec{\numflux{f}}^{\rho \vec{v}, Lim}_{(i,j)} \rvert+ \rvert \overline w_{(i,j)}^{\rho} \Delta \numflux{f}^{\rho E, Lim}_{(i,j)} \rvert \nonumber\\
    &
    + \rvert \overline w_{(i,j)}^{\rho E} \Delta \numflux{f}^{\rho, Lim}_{(i,j)}\rvert\ge P_{(i,j)}. \label{eq:pressureLimiterSimp}
\end{align}
We explore this possibility in the present paper. In the descriptions of
our numerical experiments, we call the pressure limiter  \eqref{eq:pressureLimiter} that uses \eqref{eq:pressureLimiterExact} ``sharp''. The one that uses 
\eqref{eq:pressureLimiterSimp} instead of \eqref{eq:pressureLimiterExact} is referred to as
``cautious''. As we show in the \nameref{sec:results} section, the cautious pressure fix can add much more numerical dissipation to the scheme than its sharp counterpart.

\paragraph{Semi-discrete entropy limiter}

In this work, we also use the semi-discrete entropy limiter developed in \cite{kuzmin2022limiter} for MCL schemes. Semi-discrete entropy stability of a DG or FV  method is guaranteed if the numerical fluxes satisfy  Tadmor's shuffle condition \cite{tadmor1986minimum,tadmor1983entropy,Tadmor2003}
\begin{equation} \label{eq:Tadmor}
\jump{\entVar}_{(i,j)}^T \numfluxb{f}_{(i,j)} \le \jump{\Psi}_{(i,j)},
\end{equation}
where $\jump{\Psi}_{(i,j)} := \Psi_j - \Psi_i$ denotes the jump operator, $\Psi$ is the so-called entropy-flux potential, and $\entVar$ is the vector of entropy variables.
See Appendix~\ref{app:euler} for the definition of these quantities for the Euler equations.

Let $\Delta \numfluxb{f}^{Lim}_{(i,\ii)}$ be a limited anti-diffusive flux that
is constrained to preserve local and/or global bounds for all scalar
quantities of interest. Define
\begin{equation} \label{eq:alpha_s_sync}
    \numfluxb{f}_{(i,\ii)} =
    \numfluxb{f}^{\FV}_{(i,\ii)}
    + {\alpha}^{s}_{(i,\ii)}  \frac{\Delta \numfluxb{f}^{Lim}_{(i,\ii)}}{i-j}
\end{equation}
using a correction factor $\alpha^s_{(i,j)}\in[0,1]$ such that Tadmor's condition \eqref{eq:Tadmor} is fulfilled. The substitution of 
\eqref{eq:alpha_s_sync} into \eqref{eq:Tadmor} yields a linear
inequality constraint for $\alpha^s_{(i,j)}\in[0,1]$. We enforce
this constraint using (cf. \cite{kuzmin2022limiter})
\begin{equation}
    \alpha^s_{(i,j)} = 
    \begin{cases}
    \frac{\jump{\entVar}_{(i,j)} \numfluxb{f}^{\FV}_{(i,j)} - \jump{\Psi}_{(i,j)} + \epsilon}
    {\jump{\entVar}_{(i,j)} \Delta \numfluxb{f}^{Lim}_{(i,j)}/(i-j) + \epsilon}
    & \mathrm{if}~ \Pi_{(i,j)} > \jump{\Psi}_{(i,j)}, 
    \\
    1 & \mathrm{otherwise,}
    \end{cases}
\end{equation}
where  $\epsilon$ is a small positive number and
$$
\Pi_{(i,j)}=\jump{\entVar}_{(i,j)} \left(\numfluxb{f}_{(i,j)}^{\FV} +
    \frac{\Delta \numfluxb{f}^{Lim}_{(i,j)}}{i-j}\right)
$$
is the rate of entropy production before the application of $\alpha^s_{(i,j)}$.

\begin{remark}
All correction tools described in this section lead to closed-form expressions for limited fluxes or correction factors. This distinguishes our approach from  FCT/IDP alternatives that require solving nonlinear equations \cite{RUEDARAMIREZ2022,Pazner2020,guermond2018}.
\end{remark}

\subsubsection{Definition of Bounds} \label{sec:bounds}

We distinguish between global and local bounds. Global bounds enforce physical admissibility conditions, such as positivity of the density $\rho$ and pressure $p$ in the case of the compressible Euler equations of gas dynamics.
Preservation of these  bounds is a prerequisite for running challenging simulations without crashing. If a lower bound $\rho_i^{\min}\ge 0$ is used in the limiter for $\rho$ and a (sharp or cautious) pressure limiter is applied in the final stage of the sequential limiting procedure, then the MCL scheme is positivity preserving in this sense.

The imposition of local bounds, on the other hand, makes it possible to avoid spurious oscillations within the global bounds and to improve the shock-capturing capabilities of the method. The corresponding numerical admissibility conditions are frequently formulated as local maximum or minimum principles. The inequality constraints of our MCL method are feasible if they are satisfied by the low-order bar states. Therefore, these states must be built into the definition of the upper and lower bounds. For instance, the value of a conservative or primitive quantity $\phi$ at node $i$ may be treated as numerically admissible if it is bounded by
\begin{equation} \label{eq:barStateBounds0}
    \phi^{\min}_i = \min \left\{ \phi_i, \min_{j \in \NN^*(i)} \overline{\phi}_{(i,j)} \right\},
    ~~~~
    \phi^{\max}_i = \max \left\{ \phi_i, \max_{j \in \NN^*(i)} \overline{\phi}_{(i,j)} \right\},
\end{equation}
where $\NN^*(i)=\{i-1,i+1\}$ and $\phi_i$ is the solution at node $i$ of the previous time step. 
In the multidimensional case, the integer set
$\NN^*(i)$ contains the indices of all nodes $j\ne i$ such that the flux
$\numfluxb{f}_{(i,j)}$ appears in the evolution equation for the
nodal state $\stateL{u}_i$.

We can rewrite \eqref{eq:barStateBounds0} as
\begin{equation} \label{eq:barStateBounds}
    \phi^{\min}_i = \min_{j \in \NN(i)} \overline{\phi}_{(i,j)},
    ~~~~
    \phi^{\max}_i = \max_{j \in \NN(i)} \overline{\phi}_{(i,j)},
\end{equation}
where $\NN(i)=\{i-1, i, i+1\}$ is the integer set containing $i$ and the indices of all neighboring nodes of $i$, as $\overline{\phi}_{(i,i)} = \phi_i$ by definition \eqref{def:barstates}.

If the MCL bounds \eqref{eq:barStateBounds}
are too tight, the flux-corrected scheme may fail to achieve the optimal
order of accuracy in smooth regions. Wider bounds can be constructed
by including the values of $\phi^{\FV, n+1}$ at node $i$ and its neighbors belonging
to cells that physically contain the nodal point \cite{lohmann2017,hajduk2020}.
The inclusion of extrapolated states makes it possible to guarantee
linearity preservation on general meshes \cite{kuzmin2020monolithic}.
Alternatively,
local bounds can be relaxed using
smoothness indicators \cite{kuzmin2020g} to
avoid a potential loss of accuracy due to
unnecessary limiting.

 The use of the low-order bar states to define the bounds is also common in FCT/IDP methods \cite{kuzmin2010a,kuzmin2012,Guermond2016,guermond2018}. 
Similarly to the MCL version, further
nodal states $\phi_j^{\FV, n+1}$ can be incorporated into the definition of $\phi^{\min}_i$
and $\phi^{\max}_i$. The use of smoothness indicators is also an option
\cite{guermond2018,Pazner2020}.
In FCT/IDP methods, it is also possible to define the bounds with the low-order solution at the next time step instead of using the low-order bar states, e.g. \cite{Pazner2020,kuzmin2020limiting,RUEDARAMIREZ2022},
\begin{equation} \label{eq:lowOrderBoundsFCT}
    \phi^{\min}_i = \min_{j \in \NN(i)} {\phi}^{\FV,n+1}_{j},
    ~~~~
    \phi^{\max}_i = \max_{j \in \NN(i)} {\phi}^{\FV,n+1}_{j}.
\end{equation}
This leads to more dissipative schemes, but the time-step restriction might be relaxed in some situations (see next section for details).
A peculiarity of FCT/IDP schemes is the fact that the bounds for the limited anti-diffusive fluxes are inversely proportional to the time step. As a consequence, the quality of flux-limited approximations often exhibits strong dependence on the CFL number.

In this paper, we use the tight local bounds \eqref{eq:barStateBounds} for both MCL and FCT methods without any relaxation.

\subsubsection{Explicit Time-Step Restrictions} \label{sec:timestep}

The monolithic convex limiting method enforces nodal bounds by limiting the interface fluxes, such that all high-order bar states satisfy the bounds.
This strategy enforces the bounds of the discrete solution if \eqref{eq:convexSumBarStates} is a convex sum.
As a result, we obtain the following CFL-like time-step restriction in 1D:
\begin{equation}
\label{eq:CFLcondition}
\Delta t \le
\min_i
\frac{J\,\omega_i}{\lambda^{\max}_{(i,i-1)}+\lambda^{\max}_{(i,i+1)}}.
\end{equation}

In general, the time-step restriction scales as
\begin{equation} \label{eq:CFLconditionMuldiD}
    \Delta t \lesssim \min_i \frac{m_i}{2 \sum_{d=1}^{\numDims} \lambda^{\max,d}_i},
\end{equation}
where $\numDims$ is the number of spatial dimensions of the problem, $m_i$ is the diagonal mass matrix entry of node $i$, and $\lambda^{\max,d}_i$ is the maximum wave speed at node $i$ in the coordinate direction $d$.
For instance, with $\numDims=1$ we have $\lambda^{\max,1}_i = \max \{ \lambda_{(i,i-1)}^{\max},\lambda_{(i,i+1)}^{\max} \}$.

Since the CFL condition \eqref{eq:CFLconditionMuldiD} is derived from the bar-states representation, MCL and FCT/IDP methods that use bar-state bounds \eqref{eq:barStateBounds} need to fulfill it to be able to keep the solution within bounds.
Unfortunately, \eqref{eq:CFLconditionMuldiD} is more restrictive than the typical CFL stability condition of the low-order method, especially for multiple space dimensions \cite{toro2013riemann}:
\begin{equation}
     \Delta t \lesssim \min_i
     \left\lbrace
     \min_{1\le d \le \numDims} \left\lbrace
     \frac{\Delta x^d_i}{\lambda^{\max,d}_i}
     \right\rbrace
     \right\rbrace,
\end{equation}
where $\Delta x_i^d$ is the size of the subcell in direction $d$.

To the authors' knowledge, the only way to circumvent the strict CFL condition \eqref{eq:CFLconditionMuldiD} and still keep the solution within prescribed bounds is to use FCT/IDP methods and avoid the bar-state bounds altogether.
For instance, one can use bounds computed from the (robust) low-order solution at the next time step \eqref{eq:lowOrderBoundsFCT}. 
Although this strategy can reduce the computing time in some situations, it might lead to non-physical solutions since the low order method is only provably positivity preserving (see, e.g., \cite{perthame1996positivity}) when the strict CFL condition \eqref{eq:CFLconditionMuldiD} is fulfilled.    \\

Since FCT/IDP methods apply limiting to enforce a fully discrete fulfillment of the bounds (i.e. after the time update is done), their anti-diffusion correction (and hence the spatial discretization) depends on the time-step size. On the other hand, MCL applies the limiting at the semi-discrete (spatial) level without considering the temporal discretization. As a result, the amount of dissipation depends on the CFL number for FCT/IDP methods, but not for the MCL approach. In fact, we have included a numerical investigation of the CFL dependence of the schemes and demonstrate that there is no obvious time convergence observable for FCT/IDP. 

\section{Numerical Results} \label{sec:results}

To test the convergence and robustness properties of the MCL/DGSEM schemes, we run simulations with the compressible Euler equations of gas dynamics (see Appendix~\ref{app:euler}).

In all cases, we use the Rusanov (LLF) numerical flux for the compatible robust low-order subcell scheme as well as for the surface fluxes in the high-order DG method. The CFL condition is given by \eqref{eq:CFLconditionMuldiD}.

\subsection{Convergence Test}

We simulate the advection of a density wave with initial condition
\begin{equation}
\label{eq:entropywave2D}
\rho(x) = 2+A\,\sin(2\,\pi\,(x+y)),\quad A=0.98,\quad \vec{v} = (0.1,0.2,0),\quad p=20,
\end{equation}
in the computational domain $\Omega = [-1,1]^2$, and quantify the $L_2$ error of the solution as the mesh is refined.

We use the entropy-conserving and kinetic energy preserving flux of Ranocha \cite{ranocha2018generalised} for the volume numerical flux of the split-form DGSEM method, choose a heat capacity ratio of $\gamma=1.4$, and use CFL=0.9.

We first study the effect of imposing global bounds on the solution with the MCL approach.
To do so, we impose strict positivity of density and pressure for all limited bar states,
\begin{equation} \label{eq:boundsPositivity}
    \overline{\rho}^{Lim}_{(i,j)} \ge 0,
    ~~~~
    p(\overline{\stateL{u}}^{Lim}_{(i,j)}) \ge 0, ~~~ \forall i, j\in \NN(i).
\end{equation}
The positivity of density is enforced with a one-sided MCL limiter for conservative variables, and the positivity of pressure with the sharp pressure positivity limiter.

Tables \ref{tab:eoc_n3_positivity} and \ref{tab:eoc_n4_positivity} show the $L_2$ error of all solution quantities and the experimental order of convergence (EOC) for the MCL/DGSEM scheme that imposes global bounds (positivity) for density and pressure.
The EOC is computed for four different Cartesian meshes with $N_e$ elements per spatial direction. We observe an $\mathrm{EOC} \approx N+1$ for approximations with polynomial degree $N$.

We now study the effect of imposing local bounds with the MCL approach.
To do so, we impose local minima and maxima on the density, velocity and specific total energy of all bar states using the sequential MCL limiter,
\begin{align} \label{eq:boundsSequential}
    \min_{j \in \NN (i)} \overline{\rho}_{(i,j)}
    \le \overline{\rho}_{(i,j)}^{Lim} \le
    \max_{j \in \NN (i)} \overline{\rho}_{(i,j)},&
    \min_{j \in \NN (i)} \overline{E}_{(i,j)}
    \le \overline{E}_{(i,j)}^{Lim} \le
    \max_{j \in \NN (i)} \overline{E}_{(i,j)}
    \nonumber\\
    \min_{j \in \NN (i)} \overline{v}_{1(i,j)}
    \le \overline{v}_{1(i,j)}^{Lim} \le
    \max_{j \in \NN (i)} \overline{v}_{1(i,j)},&
    \min_{j \in \NN (i)} \overline{v}_{2(i,j)}
    \le \overline{v}_{2,(i,j)}^{Lim} \le
    \max_{j \in \NN (i)} \overline{v}_{2(i,j)}
\end{align}
Tables \ref{tab:eoc_n3_sequential} and \ref{tab:eoc_n4_sequential} show the $L_2$ error of all solution quantities and the EOC for the MCL/DGSEM scheme that imposes local bounds (sequential limiting).
In this case, the experimental order of convergence is at most second order, independent of the polynomial degree, which indicates that the local bounds are too strict to achieve an $EOC>2$. Hence, what is done typically, one needs to relax the strict bounds to restore high-order accuracy, for instance by combining the local bounds of the MCL limiter with a smoothness sensor/indicator to avoid limiting smooth extrema of the approximate solution. 

\begin{table}[]
    \centering
    \caption{$L_2$ errors and EOC for the convergence test with positivity limiting (global bounds), $N=3$.}
    \resizebox{\columnwidth}{!}{
    \begin{tabular}{ccccccccc}
\hline
$N_e$ &
$\norm{\epsilon_{\rho}}^2$ & EOC &
$\norm{\epsilon_{\rho v_1}}^2$ & EOC &
$\norm{\epsilon_{\rho v_2}}^2$ & EOC &
$\norm{\epsilon_{\rho E}}^2$ & EOC \\
\hline
$4$ & $1.29 \times 10^{-1}$ & $-$ & $1.29 \times 10^{-2}$ & $-$ & $2.58 \times 10^{-2}$ & $-$ & $3.22 \times 10^{-3}$ & $-$ \\
$8$ & $1.15 \times 10^{-2}$ & $3.48$ & $1.15 \times 10^{-3}$ & $3.48$ & $2.30 \times 10^{-3}$ & $3.48$ & $2.88 \times 10^{-4}$ & $3.48$ \\
$16$ & $9.35 \times 10^{-4}$ & $3.62$ & $9.35 \times 10^{-5}$ & $3.62$ & $1.87 \times 10^{-4}$ & $3.62$ & $2.34 \times 10^{-5}$ & $3.62$ \\
$32$ & $2.25 \times 10^{-5}$ & $5.37$ & $2.25 \times 10^{-6}$ & $5.37$ & $4.51 \times 10^{-6}$ & $5.37$ & $5.64 \times 10^{-7}$ & $5.37$ \\
\hline
mean & & $4.16$ & & $4.16$ & & $4.16$ & & $4.16$ \\
\hline
    \end{tabular}
    }
    \label{tab:eoc_n3_positivity}

    \caption{$L_2$ errors and EOC for the convergence test with positivity limiting (global bounds), $N=4$.}
    \resizebox{\columnwidth}{!}{
    \begin{tabular}{ccccccccc}
\hline
$N_e$ &
$\norm{\epsilon_{\rho}}^2$ & EOC &
$\norm{\epsilon_{\rho v_1}}^2$ & EOC &
$\norm{\epsilon_{\rho v_2}}^2$ & EOC &
$\norm{\epsilon_{\rho E}}^2$ & EOC \\
\hline
$4$ & $1.75 \times 10^{-2}$ & $-$ & $1.75 \times 10^{-3}$ & $-$ & $3.50 \times 10^{-3}$ & $-$ & $4.37 \times 10^{-4}$ & $-$ \\
$8$ & $1.61 \times 10^{-3}$ & $3.44$ & $1.61 \times 10^{-4}$ & $3.44$ & $3.22 \times 10^{-4}$ & $3.44$ & $4.03 \times 10^{-5}$ & $3.44$ \\
$16$ & $2.55 \times 10^{-5}$ & $5.98$ & $2.55 \times 10^{-6}$ & $5.98$ & $5.11 \times 10^{-6}$ & $5.98$ & $6.38 \times 10^{-7}$ & $5.98$ \\
$32$ & $8.97 \times 10^{-7}$ & $4.83$ & $8.97 \times 10^{-8}$ & $4.83$ & $1.79 \times 10^{-7}$ & $4.83$ & $2.24 \times 10^{-8}$ & $4.83$ \\
\hline
mean & & $4.75$ & & $4.75$ & & $4.75$ & & $4.75$ \\
\hline
    \end{tabular}
    }
    \label{tab:eoc_n4_positivity}
\end{table}

\begin{table}[]
    \centering
    \caption{$L_2$ errors and EOC for the convergence test with sequential limiting (local bounds), $N=3$.}
    \resizebox{\columnwidth}{!}{
    \begin{tabular}{ccccccccc}
\hline
$N_e$ &
$\norm{\epsilon_{\rho}}^2$ & EOC &
$\norm{\epsilon_{\rho v_1}}^2$ & EOC &
$\norm{\epsilon_{\rho v_2}}^2$ & EOC &
$\norm{\epsilon_{\rho E}}^2$ & EOC \\
\hline
$4$ & $2.51 \times 10^{-1}$ & $-$ & $2.51 \times 10^{-2}$ & $-$ & $5.02 \times 10^{-2}$ & $-$ & $6.28 \times 10^{-3}$ & $-$ \\
$8$ & $7.74 \times 10^{-2}$ & $1.70$ & $7.74 \times 10^{-3}$ & $1.70$ & $1.55 \times 10^{-2}$ & $1.70$ & $1.92 \times 10^{-3}$ & $1.71$ \\
$16$ & $1.91 \times 10^{-2}$ & $2.02$ & $1.91 \times 10^{-3}$ & $2.02$ & $3.82 \times 10^{-3}$ & $2.02$ & $4.77 \times 10^{-4}$ & $2.01$ \\
$32$ & $5.70 \times 10^{-3}$ & $1.74$ & $5.70 \times 10^{-4}$ & $1.74$ & $1.14 \times 10^{-3}$ & $1.74$ & $1.42 \times 10^{-4}$ & $1.75$ \\
\hline
mean & & $1.82$ & & $1.82$ & & $1.82$ & & $1.82$\\
\hline
    \end{tabular}
    }
    \label{tab:eoc_n3_sequential}

    \caption{$L_2$ errors and EOC for the convergence test with sequential limiting (local bounds), $N=4$.}
    \resizebox{\columnwidth}{!}{
    \begin{tabular}{ccccccccc}
\hline
$N_e$ &
$\norm{\epsilon_{\rho}}^2$ & EOC &
$\norm{\epsilon_{\rho v_1}}^2$ & EOC &
$\norm{\epsilon_{\rho v_2}}^2$ & EOC &
$\norm{\epsilon_{\rho E}}^2$ & EOC \\
\hline
$4$ & $1.27 \times 10^{-1}$ & $-$ & $1.27 \times 10^{-2}$ & $-$ & $2.53 \times 10^{-2}$ & $-$ & $3.18 \times 10^{-3}$ & $-$ \\
$8$ & $4.19 \times 10^{-2}$ & $1.59$ & $4.19 \times 10^{-3}$ & $1.60$ & $8.38 \times 10^{-3}$ & $1.59$ & $1.04 \times 10^{-3}$ & $1.61$ \\
$16$ & $1.19 \times 10^{-2}$ & $1.81$ & $1.19 \times 10^{-3}$ & $1.81$ & $2.39 \times 10^{-3}$ & $1.81$ & $2.98 \times 10^{-4}$ & $1.80$ \\
$32$ & $4.47 \times 10^{-3}$ & $1.42$ & $4.47 \times 10^{-4}$ & $1.42$ & $8.94 \times 10^{-4}$ & $1.42$ & $1.12 \times 10^{-4}$ & $1.42$ \\
\hline
mean & & $1.61$ & & $1.61$ & & $1.61$ & & $1.61$\\
\hline
    \end{tabular}
    }
    \label{tab:eoc_n4_sequential}
\end{table}

\subsection{Kelvin-Helmholtz Instability}

We consider the inviscid two-dimensional Kelvin-Helmholtz instability (KHI) setup, e.g., presented in \cite{RR2021,RUEDARAMIREZ2022}.
Due to its high density contrast and compressibility effects, the test case is challenging for nodal high-order methods when the (under-resolved) vortical structures of the KHI develop and evolve. In fact, the standard LGL-DGSEM method requires limiting to ensure robustness (positivity in this case).

The initial condition is given by
\begin{align} \label{eq:KHI_IniCond}
\rho (x,y) &= \frac{1}{2}
+ \frac{3}{4} B,
&
p (x,y) &= 1,
\nonumber \\
v_{1} (x,y) &= \frac{1}{2} \left( B-1 \right),
&
v_{2} (x,y) &= \frac{1}{10} \sin(2 \pi x),
\end{align}
with
$
B=
\tanh \left( 15 y + 7.5 \right) - \tanh(15y-7.5).
$

We tessellate the simulation domain, $\Omega=[-1,1]^2$, using $64 \times 64$ quadrilateral elements, use periodic boundary conditions, represent the solution with polynomials of degree $N=7$, and run the simulation until the final time $t=10$.
Moreover, we discretize the Euler equations using the split-form DGSEM and the entropy-conserving and kinetic energy preserving flux of Ranocha \cite{ranocha2018generalised} for the volume numerical fluxes, and select CFL=0.9

We first study the effect of imposing global bounds on the solution with the MCL and FCT/IDP approaches.
With the MCL method, it is possible to impose strict positivity of density and pressure for all limited bar states \eqref{eq:boundsPositivity} using a one-sided MCL limiter for conservative quantities and the sharp pressure positivity limiter. 
However, with the FCT/IDP approach presented in \cite{RUEDARAMIREZ2022} a positive threshold greater than zero is necessary, as otherwise some nodes might get an invalid vacuum state. 
Hence, for the FCT/IDP variant, we consider the heuristic positivity-preserving method of Rueda-Ramírez \cite{RR2021} in a subcell-wise manner, i.e., we impose lower bounds for density and pressure that depend on the FV solution,
\begin{equation} \label{eq:rhop_restriction}
\rho_{i} \ge \beta \rho^{\FV}_{i},
~~~~
p_{i} \ge \beta p^{\FV}_{i},
\end{equation}
with $\beta = 0.1$. 
We note, that this is a somewhat stricter requirement than strict positivity.

\begin{figure}[h!]
\centering
\includegraphics[trim=350 700 350 0 ,clip,width=0.5\linewidth]{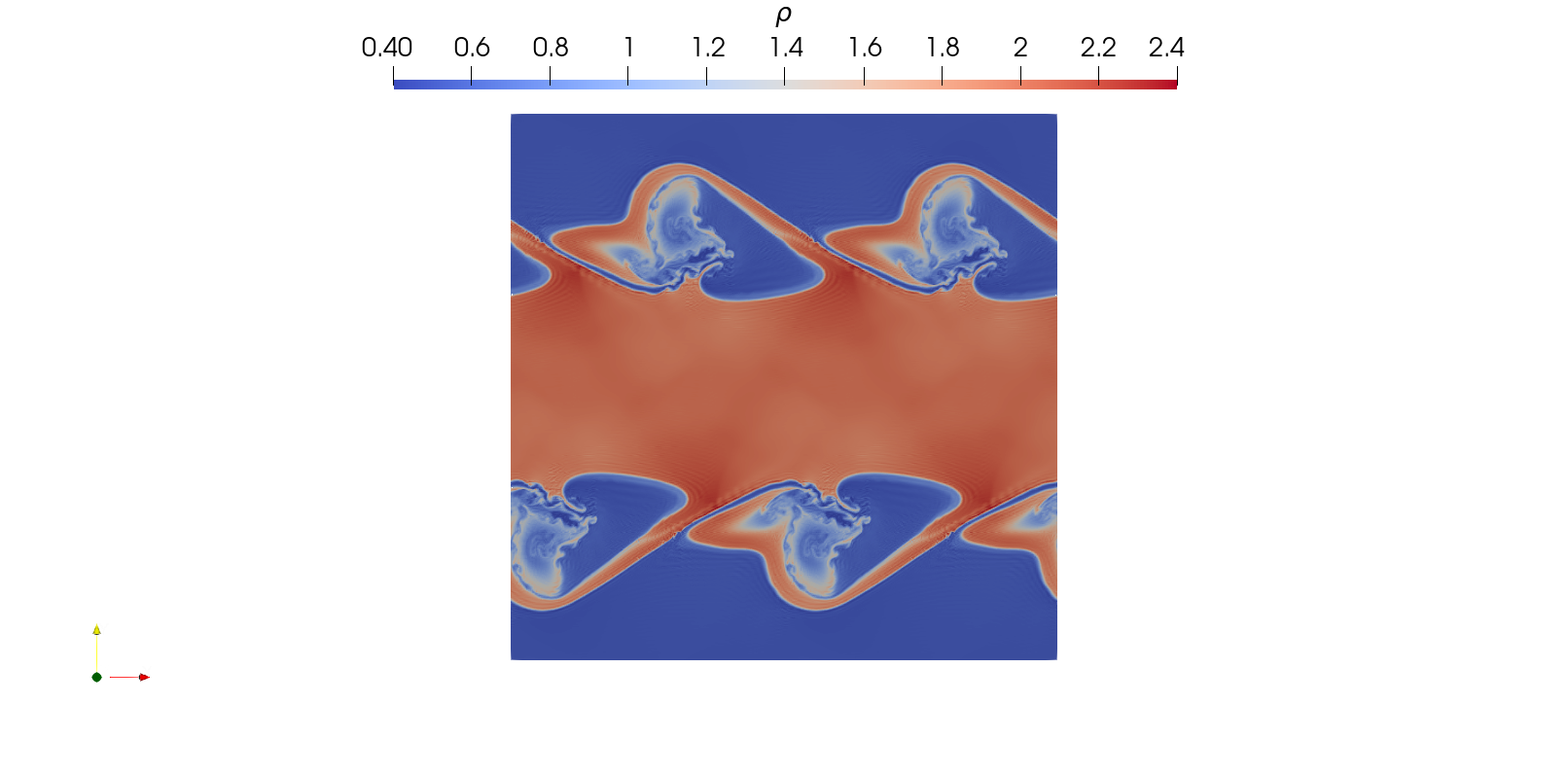}
\\
\begin{subfigure}[b]{0.33\linewidth}
    \centering
	\includegraphics[trim=520 120 520 125 ,clip,width=\linewidth]{figures/kelvinhelmholtz/mcl_positivity_dens_pres_exact_t3.7_cfl_0.9.png}
	\hspace{50pt}
	\caption{MCL, $t=3.7$}
\end{subfigure}
\begin{subfigure}[b]{0.33\linewidth}
    \centering
	\includegraphics[trim=520 120 520 125,clip,width=\linewidth]{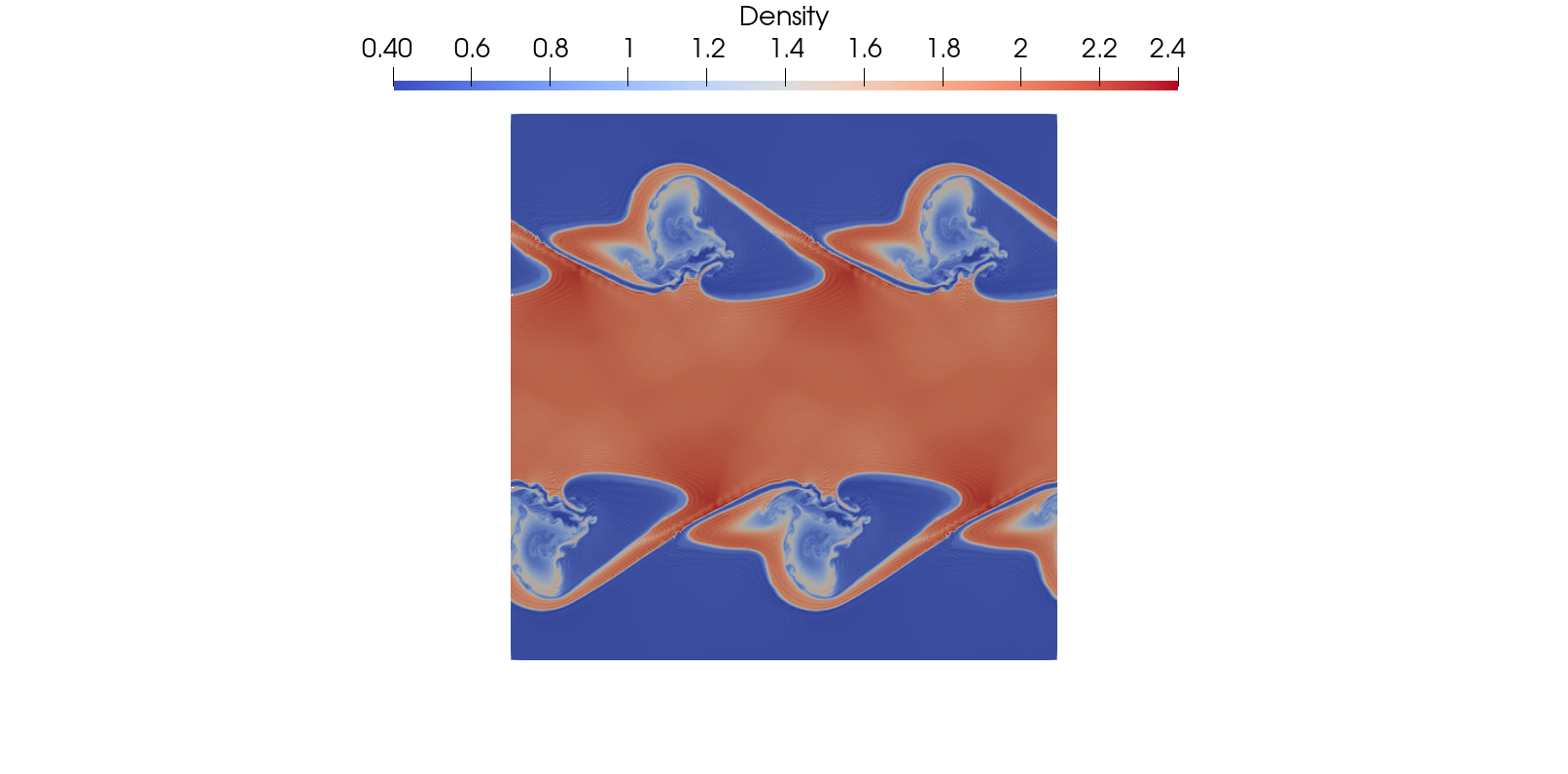}
	\hspace{50pt}
	\caption{FCT/IDP, $t=3.7$}
\end{subfigure}
\begin{subfigure}[b]{0.33\linewidth}
    \centering
	\includegraphics[trim=520 120 520 125 ,clip,width=\linewidth]{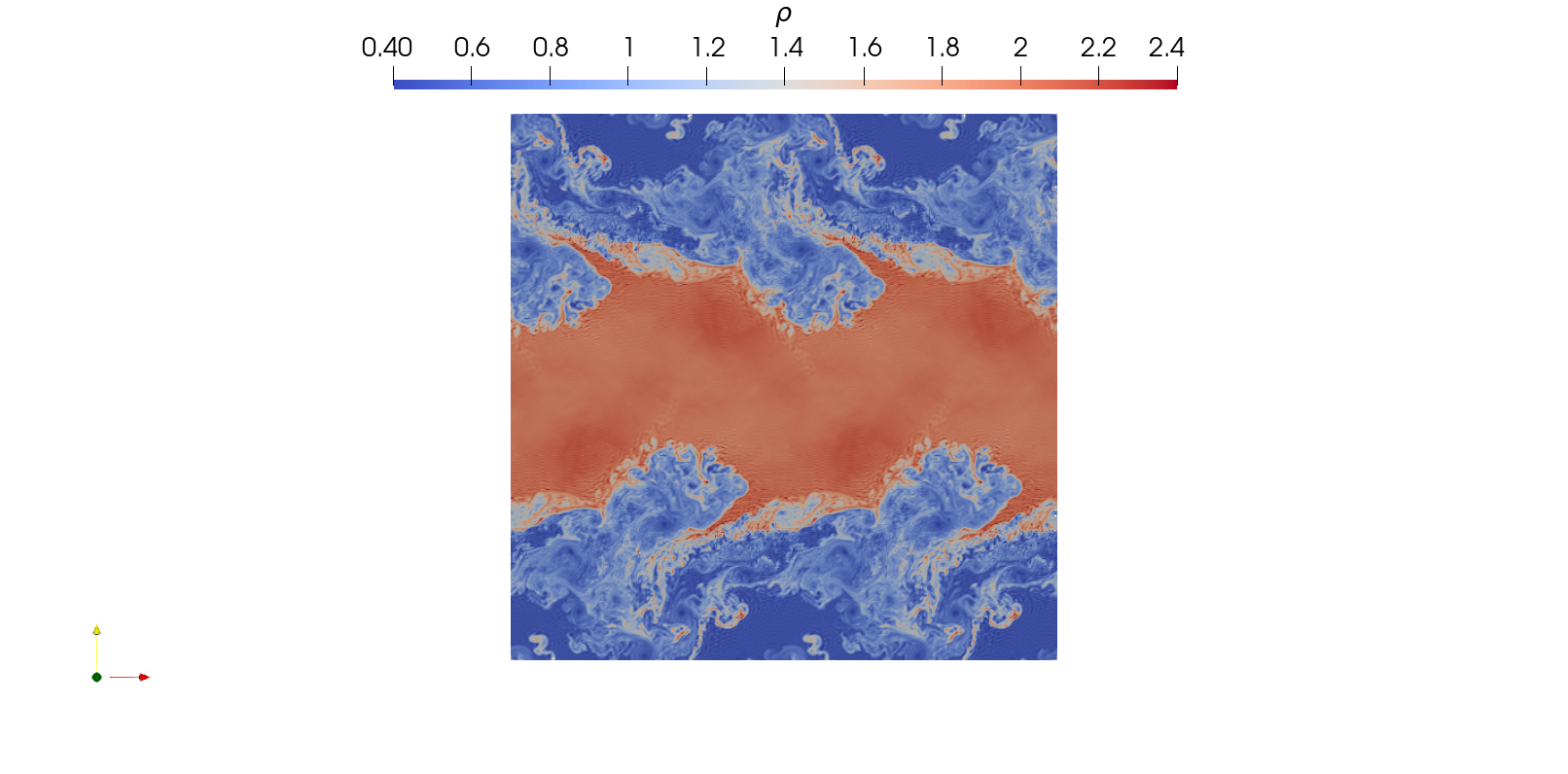}
	\hspace{50pt}
	\caption{MCL, $t=6.7$}
\end{subfigure}
\begin{subfigure}[b]{0.33\linewidth}
    \centering
	\includegraphics[trim=520 120 520 125 ,clip,width=\linewidth]{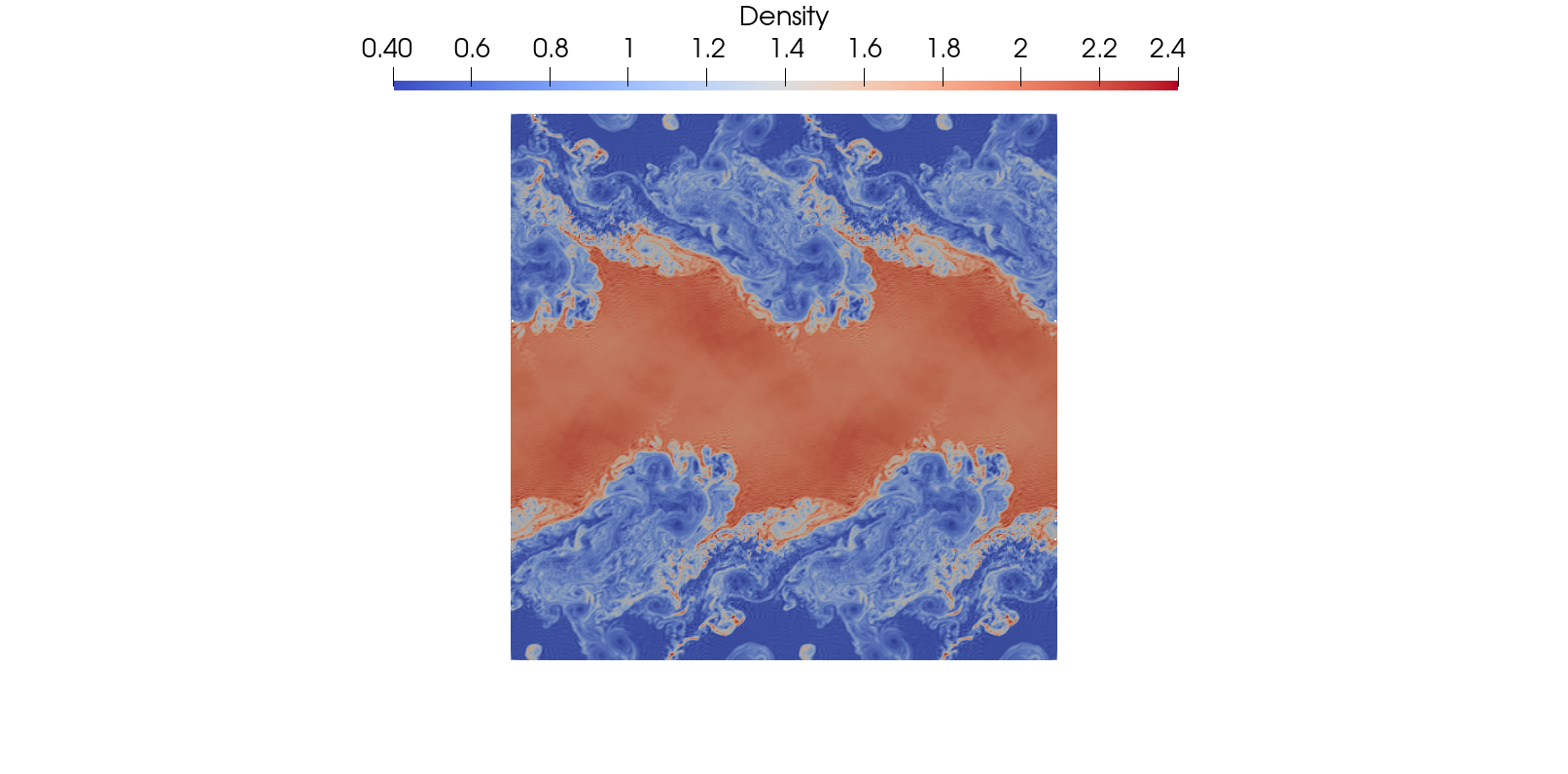}
	\hspace{50pt}
	\caption{FCT/IDP, $t=6.7$}
\end{subfigure}
\begin{subfigure}[b]{0.33\linewidth}
    \centering
	\includegraphics[trim=520 120 520 125 ,clip,width=\linewidth]{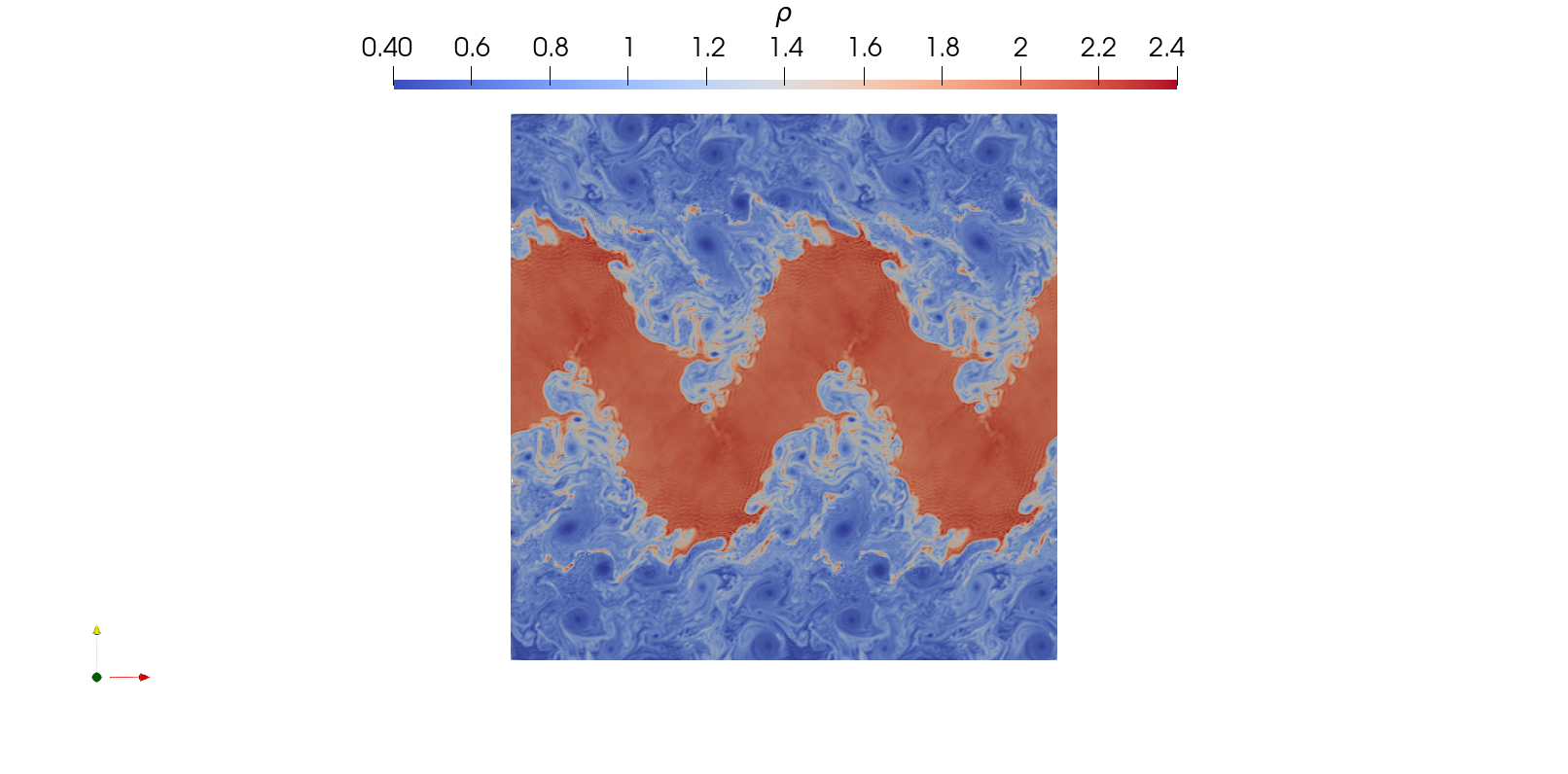}
	\hspace{50pt}
	\caption{MCL, $t=10$}
\end{subfigure}
\begin{subfigure}[b]{0.33\linewidth}
    \centering
	\includegraphics[trim=520 120 520 125 ,clip,width=\linewidth]{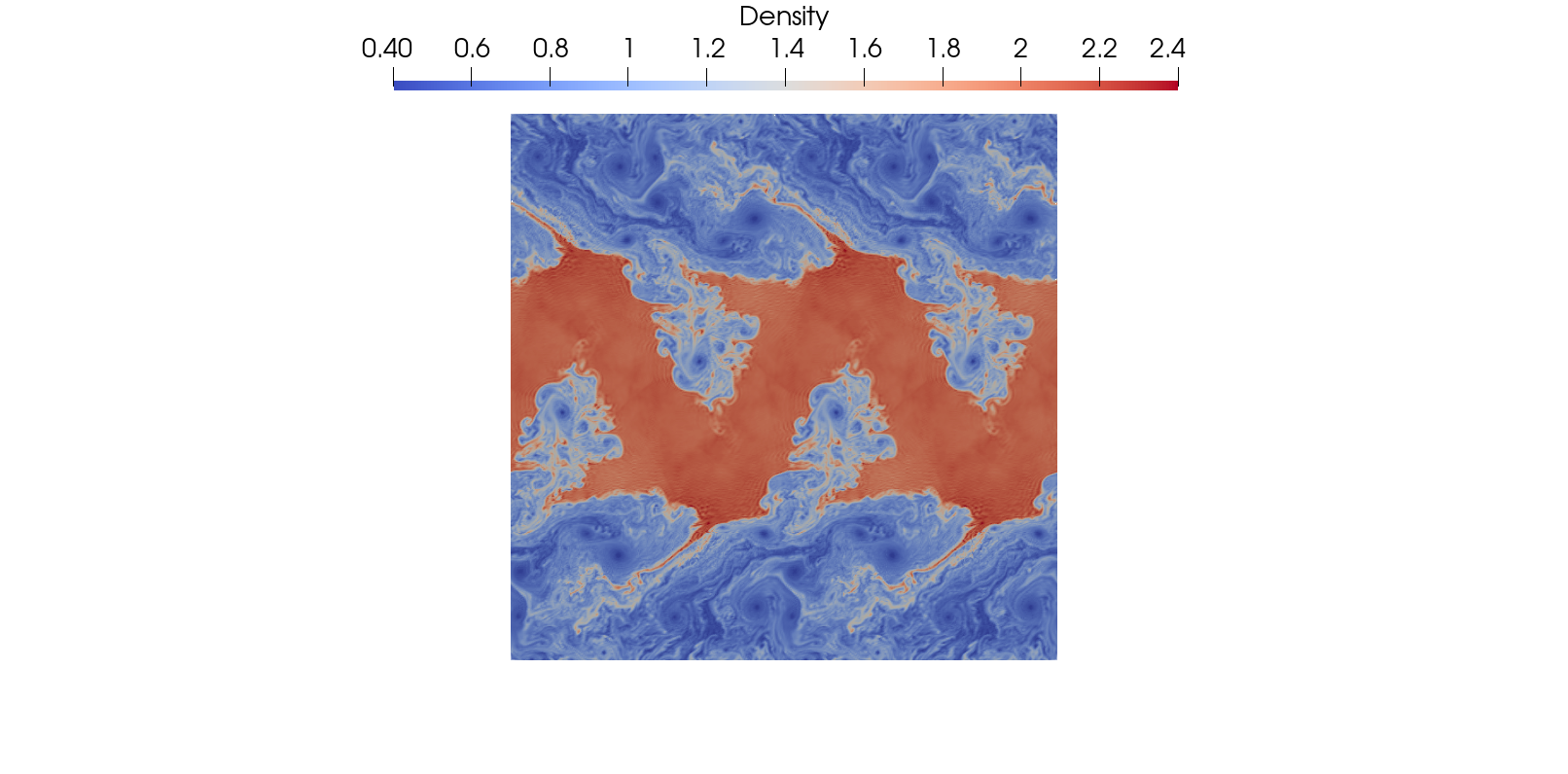}
	\hspace{50pt}
	\caption{FCT/IDP, $t=10$}
\end{subfigure}

\caption{Density contours for the Kelvin-Helmholtz simulations with MCL and FCT/IDP (positivity limiting for density and pressure).
DGSEM results with polynomial degree $N=7$ and $64\times 64$ elements.}
\label{fig:kelvinhelmholtz_MCL_positivity}
\end{figure}

Figure~\ref{fig:kelvinhelmholtz_MCL_positivity} shows the density contours at different stages of the KHI simulation using the MCL and FCT/IDP limiters with global bounds - both approaches run stably until the final time.
Even though this simulation setup is extremely sensitive to the discretization scheme \cite{RUEDARAMIREZ2022} regarding the shape and form of the vortex roll-ups, the two approaches to impose global bounds produce remarkably similar looking results.

We now study the effect of imposing local bounds on the solution with the MCL and FCT/IDP approaches.
With the MCL method, we use the standard sequential limiting to impose local minima and maxima on the density, velocity and total energy \eqref{eq:boundsSequential}.
With the FCT/IDP method, we impose local minima and maxima on the density, and local minima on the specific entropy,
\begin{equation} \label{eq:IDPbounds}
    \min_{j \in \NN (i)} \bar{\rho}_{(i,j)}
    \le \rho_{i} \le
    \max_{j \in \NN (i)} \bar{\rho}_{(i,j)},
    ~~~~
    \min_{j \in \NN (i)} \eta(\bar{\state{u}}_{(i,j)})
    \le \eta({\state{u}}_{i}),
\end{equation}
where the condition on the modified specific entropy,
$\eta = e \rho^{1- \gamma}$,
guarantees the fulfillment of a discrete entropy inequality \cite{guermond2019}. 
Moreover, $\eta$ is an efficient choice since it is computationally cheaper to evaluate than the specific entropy $s=\ln(p/\rho^{\gamma})$ and it improves the convergence of the Newton method that is used in FCT/IDP methods to solve the non-linear equation to obtain the limiting factor \cite{guermond2019,maier2021efficient}.
Note that condition \eqref{eq:boundsSequential} is typical for MCL and condition \eqref{eq:IDPbounds} is standard for FCT/IDP.

\begin{figure}[h!]
\centering
\includegraphics[trim=350 700 350 0 ,clip,width=0.5\linewidth]{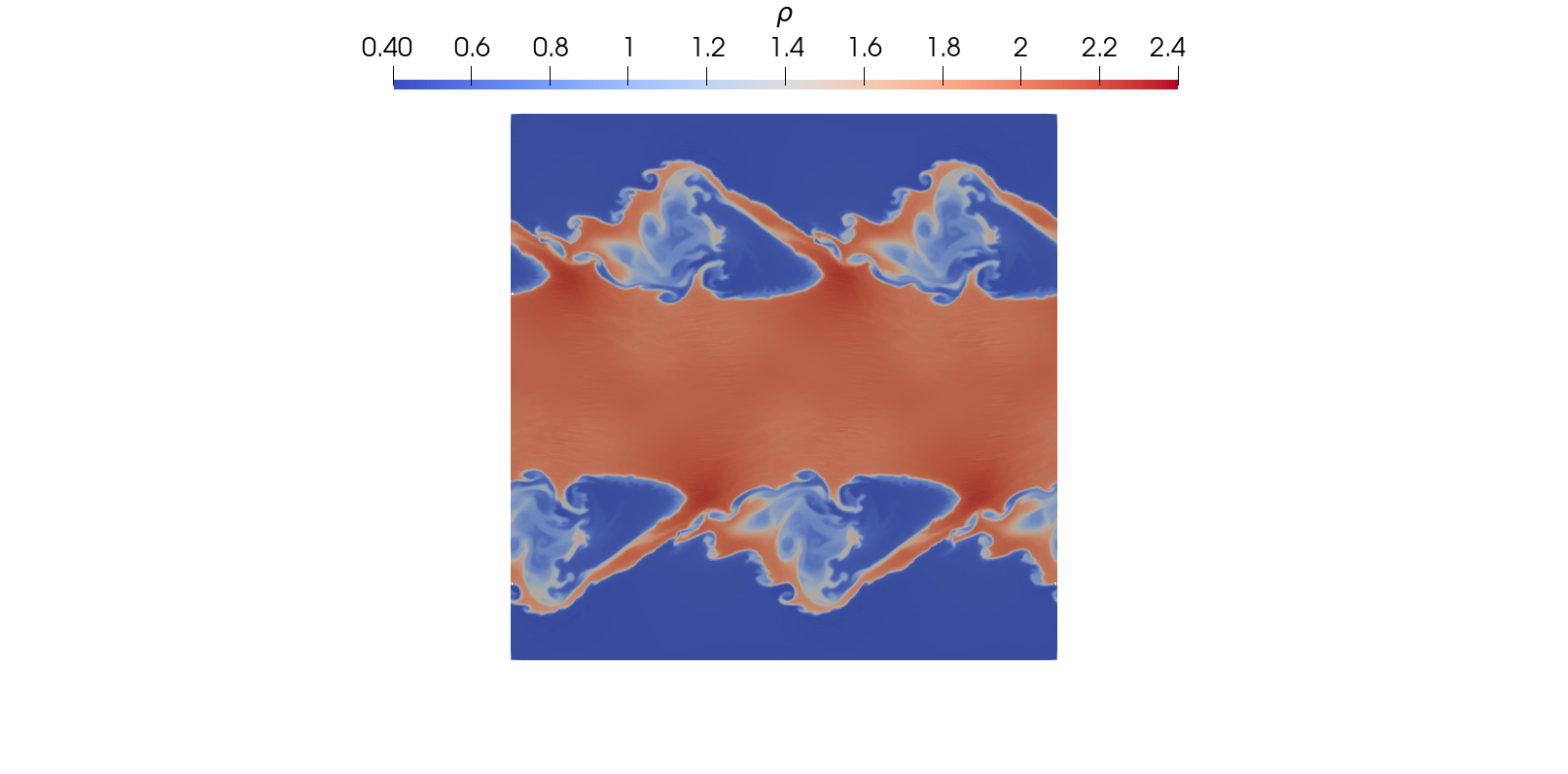}
\\
\begin{subfigure}[b]{0.33\linewidth}
    \centering
	\includegraphics[trim=520 120 520 125 ,clip,width=\linewidth]{figures/kelvinhelmholtz/mcl_sequential_t3.7_cfl_0.9.png}
	\hspace{50pt}
	\caption{MCL, $t=3.7$}
\end{subfigure}
\begin{subfigure}[b]{0.33\linewidth}
    \centering
	\includegraphics[trim=520 120 520 125 ,clip,width=\linewidth]{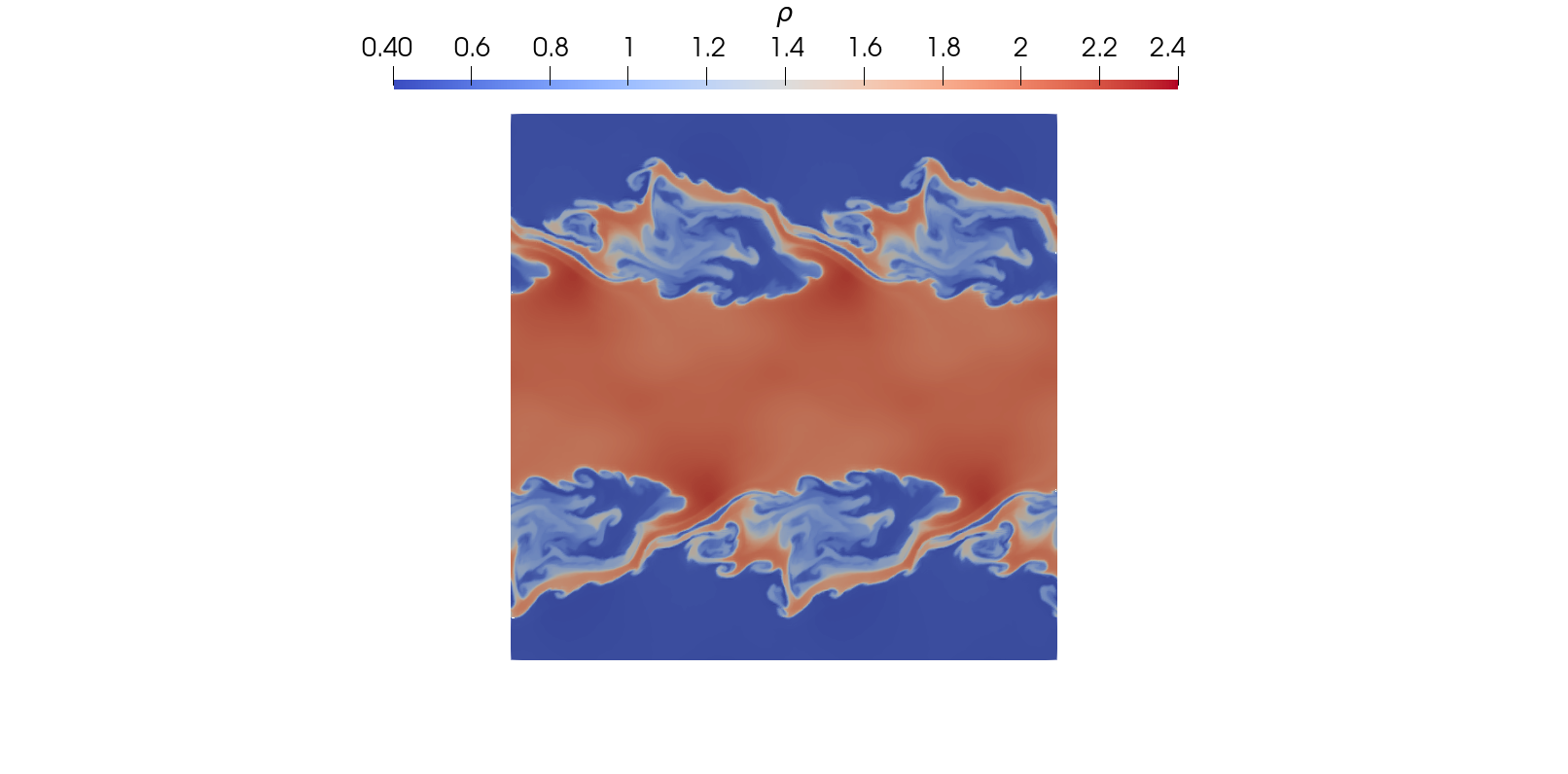}
	\hspace{50pt}
	\caption{FCT/IDP, $t=3.7$}
\end{subfigure}
\begin{subfigure}[b]{0.33\linewidth}
    \centering
	\includegraphics[trim=520 120 520 125 ,clip,width=\linewidth]{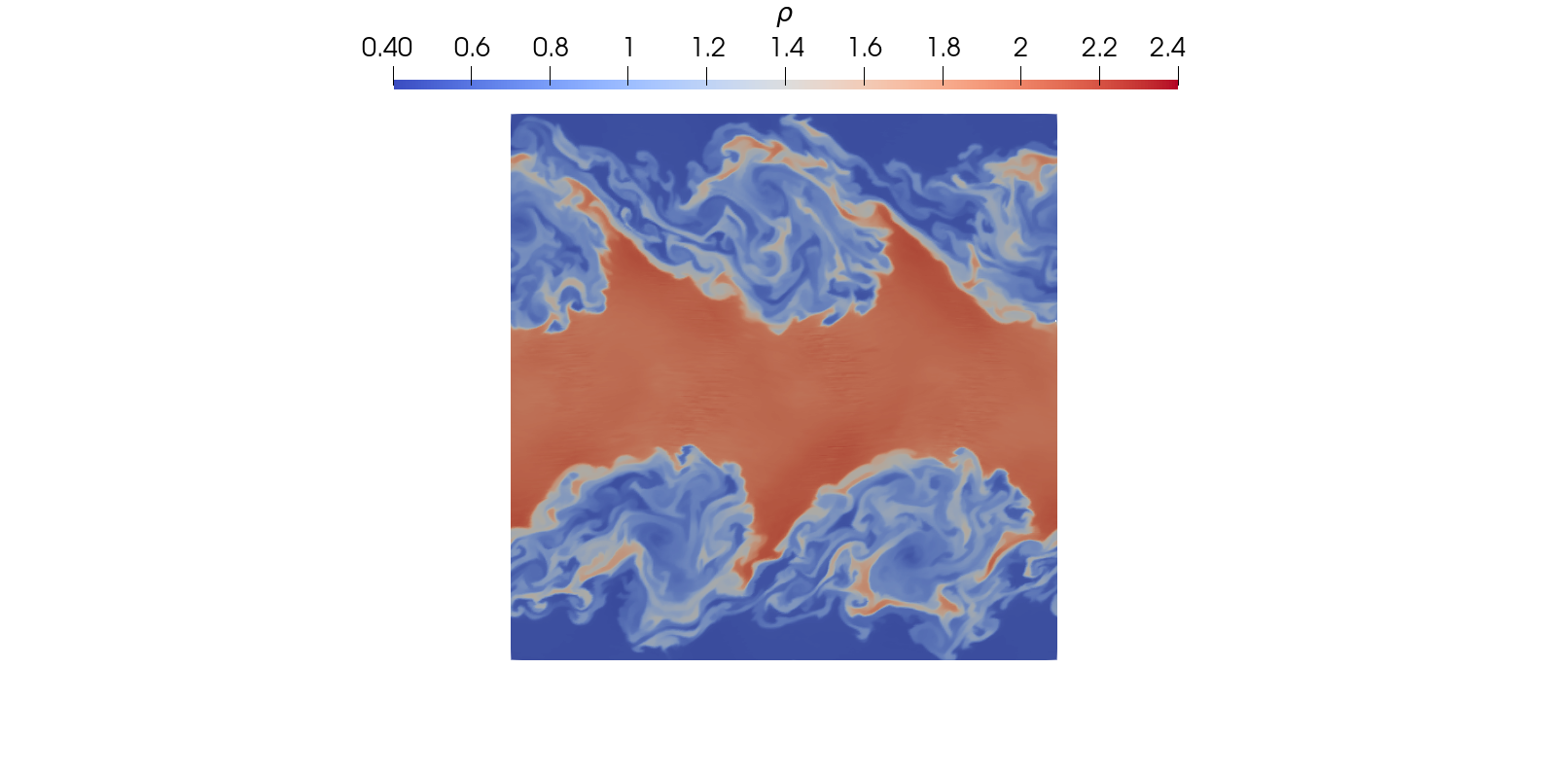}
	\hspace{50pt}
	\caption{MCL, $t=6.7$}
\end{subfigure}
\begin{subfigure}[b]{0.33\linewidth}
    \centering
	\includegraphics[trim=520 120 520 125 ,clip,width=\linewidth]{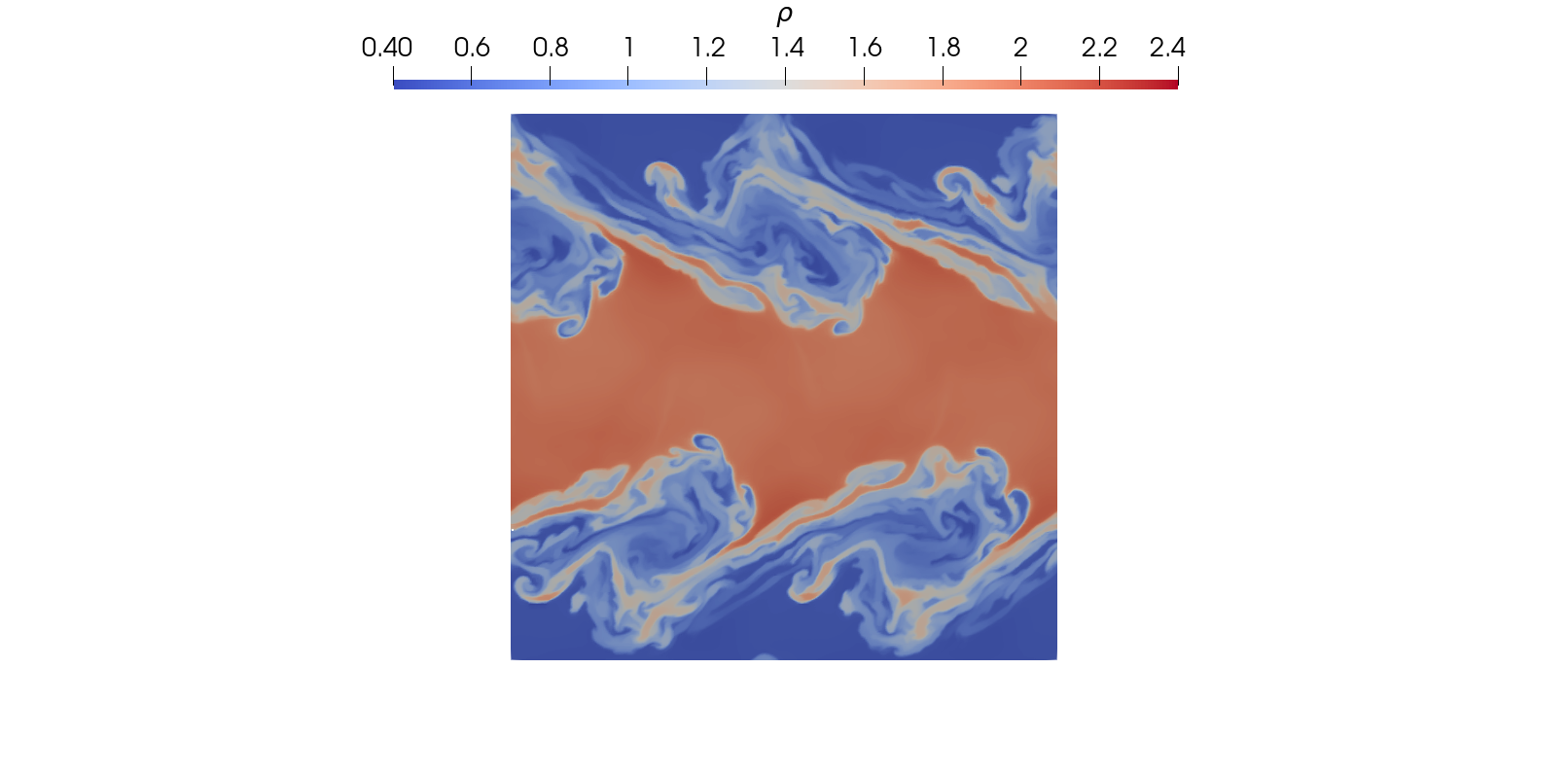}
	\hspace{50pt}
	\caption{FCT/IDP, $t=6.7$}
\end{subfigure}
\begin{subfigure}[b]{0.33\linewidth}
    \centering
	\includegraphics[trim=520 120 520 125 ,clip,width=\linewidth]{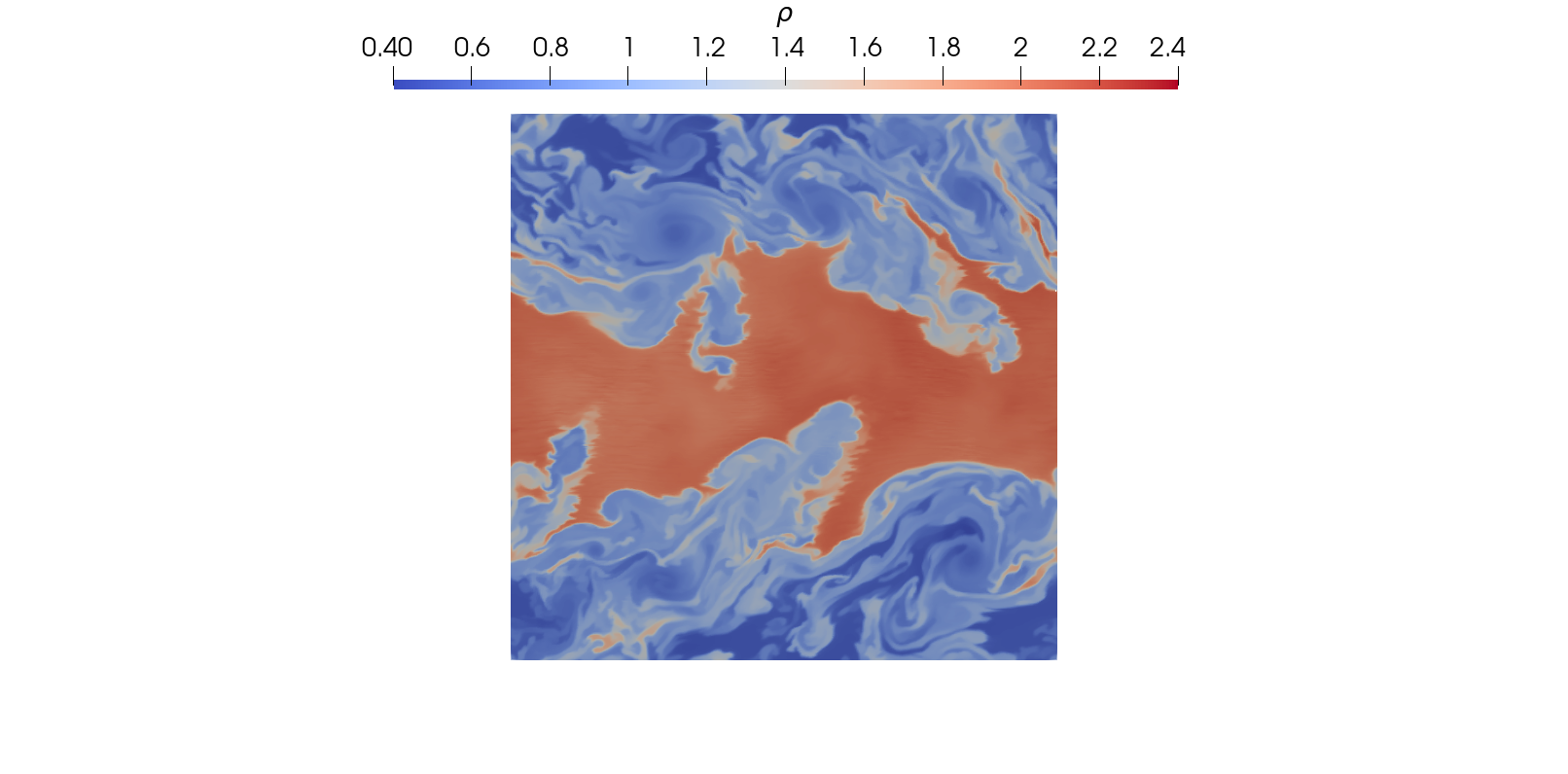}
	\hspace{50pt}
	\caption{MCL, $t=10$}
\end{subfigure}
\begin{subfigure}[b]{0.33\linewidth}
    \centering
	\includegraphics[trim=520 120 520 125 ,clip,width=\linewidth]{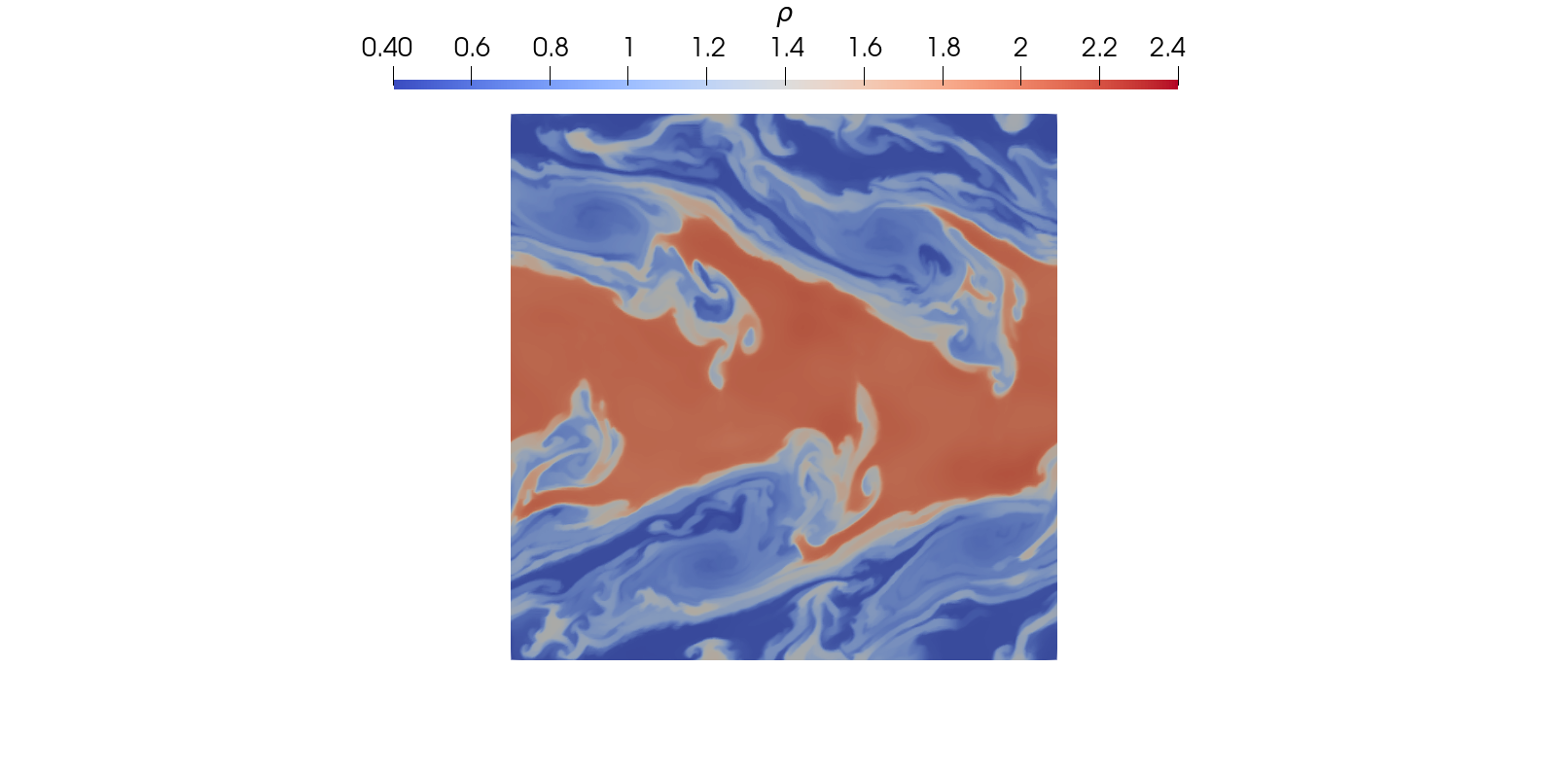}
	\hspace{50pt}
	\caption{FCT/IDP, $t=10$}
\end{subfigure}

\caption{Density contours for the Kelvin-Helmholtz simulations with MCL (sequential limiting) and FCT/IDP (density and specific entropy).
DGSEM results with polynomial degree $N=7$ and $64\times 64$ elements.}
\label{fig:kelvinhelmholtz_MCL_sequential}
\end{figure}

Figure~\ref{fig:kelvinhelmholtz_MCL_sequential} shows the density contours at different stages of the KHI simulation using the MCL and FCT/IDP limiters with local bounds.
Again, the solutions are very comparable between the two approaches, even though the limiting techniques and bounds are different. When comparing Figures~\ref{fig:kelvinhelmholtz_MCL_positivity}~and~\ref{fig:kelvinhelmholtz_MCL_sequential}, it is evident that the methods that impose local bounds add more numerical dissipation than the methods with global bounds, as expected. There is a smaller range of scales apparent with local bounds, especially at larger times. 

Finally, we compute the mean limiting factor as
\begin{equation} \label{eq:alpha_bar}
    \bar \alpha (t) = \left( \frac{1}{V} \sum_{e=1}^K \sum_{i,j=0}^N J_{ij} \omega_{ij} \alpha^e_{ij}(t) \right),
\end{equation}
where $e \in [1,K]$ denotes the element index, $K$ is the number of elements of the domain, $i,j \in [0,N]$ are the node indexes, $N$ is the polynomial degree, $\alpha^e_{ij}(t)$ is the limiting factor of node $ij$ of element $e$ at time $t$, and $V$ is the total area of the domain.
Since in MCL methods the limiting is done for each interface and each equation without using a limiting factor, we first compute the effective limiting factor for each interface and each equation, and then compute a nodal $\alpha^e_{ij}$ using the average over all interfaces of node $ij$,
\begin{equation}\label{eq:nodalAlpha}
    \alpha^e_{ij} := \frac{1}{4} \sum_{k\in \NN^e(ij)} \alpha^e_{(ij,k)}.
\end{equation}

\begin{figure}[h!]
\centering
\includegraphics[trim=0 0 0 0 ,clip,width=0.65\linewidth]{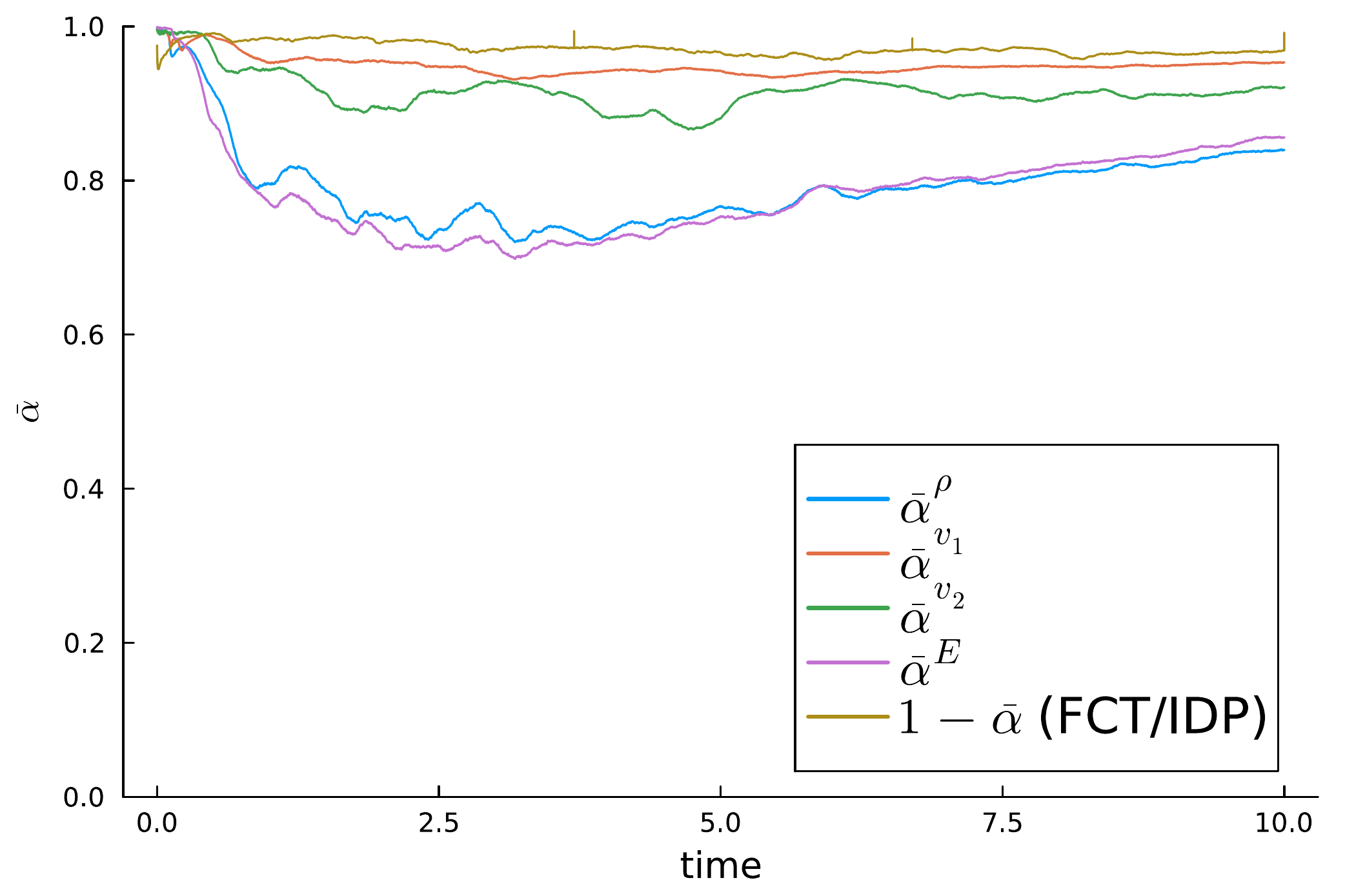}

\caption{Evolution of limiting factors for MCL (sequential limiting) and FCT/IDP (density and specific entropy).
DGSEM results with polynomial degree $N=7$ and $64\times 64$ elements.}
\label{fig:kelvinhelmholtz_MCL_alphas}
\end{figure}

We present a plot of the evolution of the mean limiting factors for the KHI simulations in Figure~\ref{fig:kelvinhelmholtz_MCL_alphas}.
To obtain limiting factors between $0$ and $1$, the factors for the momentum and energy equations of MCL are the effective scaling of the auxiliary flux, i.e., 
\begin{equation} \label{eq:alphaPrimitive}
    \alpha^{e,\phi}_{(ij,k)} :=
    \begin{cases}
    1 & \mathrm{if}~\Delta \numflux{g}^{\phi}_{(ij,k)} \approx 0, \\
    \frac{\Delta \numflux{g}^{\phi,Lim}_{(ij,k)} + \epsilon \, \texttt{sign}(\Delta \numflux{g}^{\phi}_{(ij,k)})}
    {\Delta \numflux{g}^{\phi}_{(ij,k)} + \epsilon \, \texttt{sign}(\Delta \numflux{g}^{\phi}_{(ij,k)})} & \mathrm{otherwise},
    \end{cases}
\end{equation}
where $\epsilon$ is a very small number.
We plot the quantity $(1-\alpha)$ for FCT/IDP methods since the limiting factor of our FCT/IDP methods \cite{RUEDARAMIREZ2022} is defined inversely as for MCL methods.
A mean limiting factor $\bar \alpha = 1$ means that the discretization uses the unlimited high-order scheme everywhere, whereas a mean limiting factor $\bar \alpha = 0$ means that the anti-diffusion fluxes, $\Delta \numfluxb{f}$ for FCT/IDP, or $\Delta \numflux{f}^{\rho}$ or $\Delta \numflux{g}^{\phi}$ for MCL, are set to zero everywhere.

\subsection{Inviscid Bow Shock Upstream of a Blunt Body}

We consider the supersonic flow over a 2D blunt body that produces a detached bow shock to test the performance of the MCL/DGSEM method on curvilinear grids.
This problem setup was proposed as an advanced test case for the Fifth International Workshop on High-Order CFD Methods \cite{hiocfd5}.

The left boundary of the domain is a circular segment with origin $(3.85,0)$ and radius $5.9$, the blunt body has a flat front of length $1$ connected with two quarter circles of radius $0.5$, and the right boundary is located at $x=0$.
The heat capacity ratio is set to $\gamma=1.4$ and the initial condition is the constant state
\begin{equation}
\rho(x,y) = 1.4, \qquad
p(x,y) = 1, \qquad
v_1(x,y) = 4, \qquad
v_2(x,y) = 0,
\end{equation}
which corresponds to a Mach number $\mathrm{Ma}=1.4$.

For the blunt body we use a reflecting wall boundary condition, while for the other boundaries we use characteristics-based inflow/outflow boundaries, on which the external state is selected depending on the flow conditions normal to the boundary.

We use the split-form DGSEM with the entropy-conserving and kinetic energy preserving flux of Ranocha \cite{ranocha2018generalised}, polynomial degree $N=5$, an isoparametric mapping of the geometry, MCL limiting with global \eqref{eq:boundsPositivity} and local \eqref{eq:boundsSequential} bounds, and a conforming mesh with $36$ elements distributed regularly on the inflow and wall boundaries and $24$ elements distributed regularly on the outflow boundaries.
To impose positivity of pressure, we use again the sharp pressure positivity limiter.

\begin{figure}[h!]
\centering
\begin{subfigure}[b]{0.35\linewidth}
    \centering
	\includegraphics[trim=700 0 500 0 ,clip,width=0.8\linewidth]{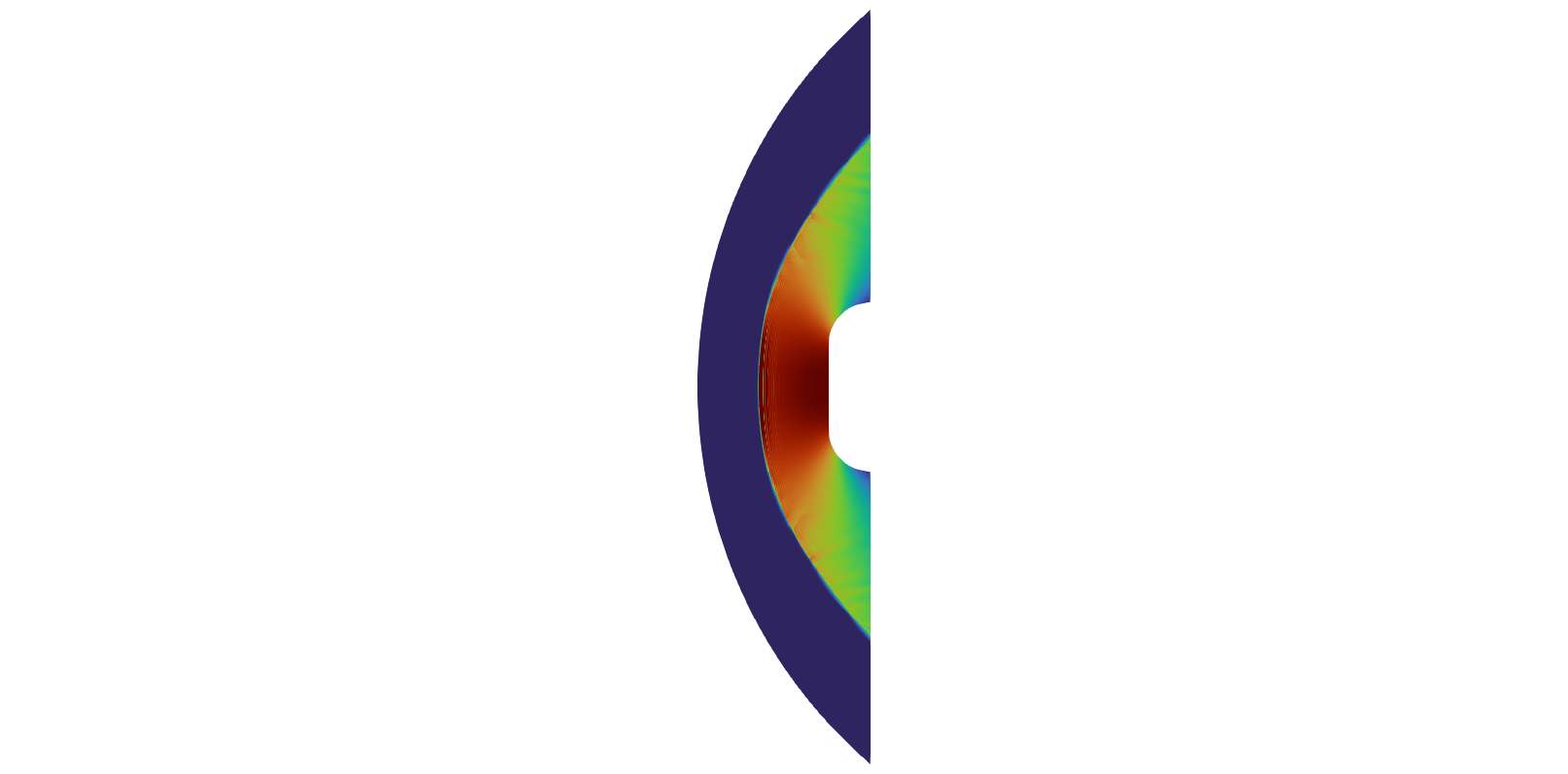}
	\hspace{50pt}
	\caption{Global bounds}
	\label{fig:bowshock_MCL_positivity}
\end{subfigure}
\begin{subfigure}[b]{0.35\linewidth}
    \centering
\includegraphics[trim=700 0 500 0 ,clip,width=0.8\linewidth]{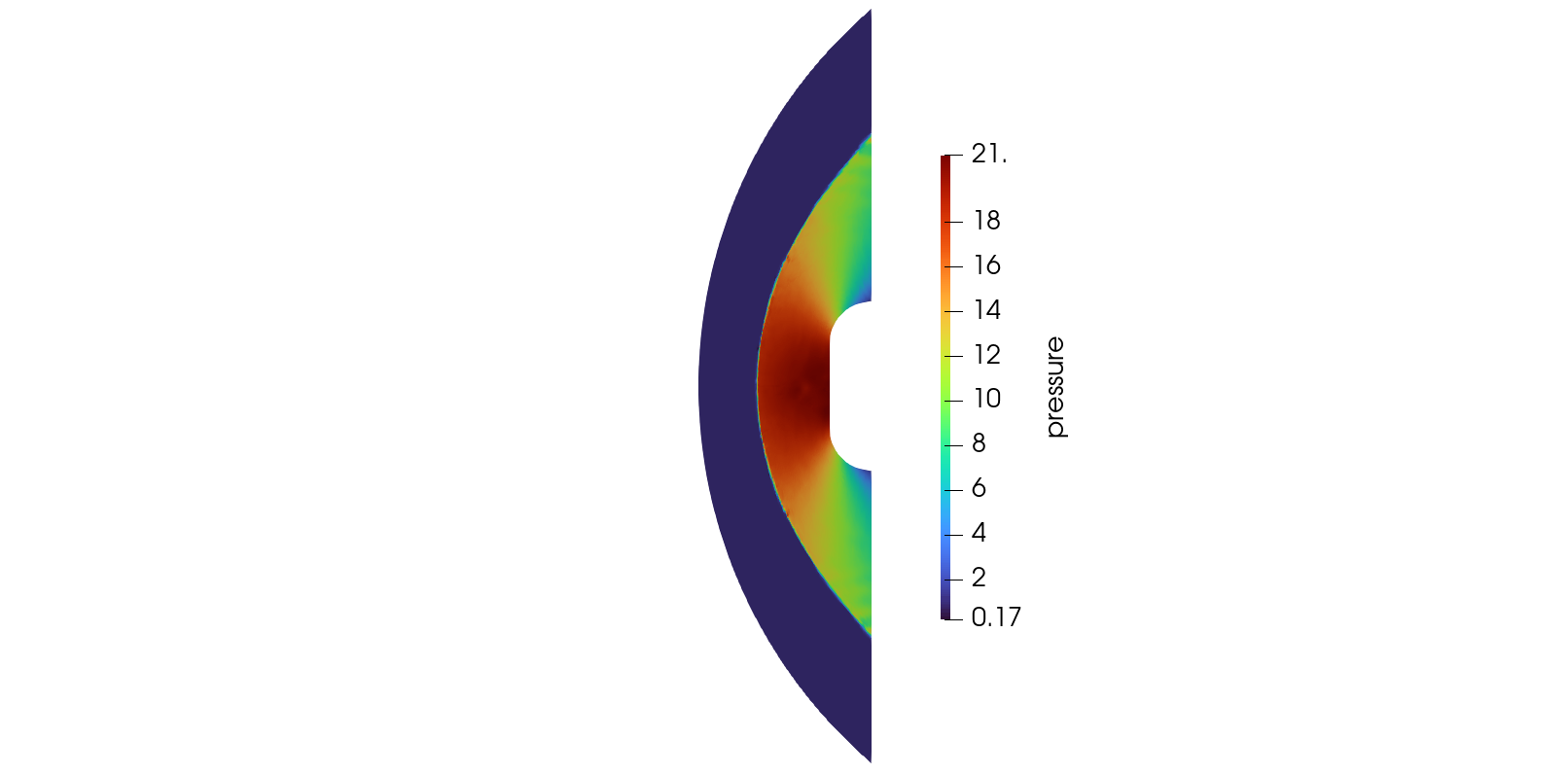}
	\hspace{50pt}
	\caption{Local bounds}
	\label{fig:bowshock_MCL_sequential}
\end{subfigure}

\caption{Pressure contours for the bow shock simulations with MCL at time $t=10$.
DGSEM results with polynomial degree $N=5$ and $24 \times 36$ elements. Left plot shows the simulation with global bounds only (positivity) and right plot shows the result with local bounds (sequential limiting).}
\label{fig:bowshock_MCL}
\end{figure}

Figure~\ref{fig:bowshock_MCL} shows the pressure contours for the bow shock simulations at time $t=10$ using MCL limiters with global and local bounds.
Figure~\ref{fig:bowshock_MCL_positivity} shows that global positivity bounds are enough to keep the simulation running until the end time, but spurious oscillations appear near the shock.
The use of local bounds removes the spurious oscillations near the shock, as can be observed in Figure~\ref{fig:bowshock_MCL_sequential}.

\subsection{Sedov Blast Explosion}

The Sedov blast problem is a very challenging simulation setup with a strong circular shock that describes the evolution of a symmetrical blast wave expanding from an initial concentration of density and pressure into a gas at rest.

For the initial condition, we use the standard setup from the FLASH astrophysical code
\cite{fryxell2000flash}.
The gas is initially at rest, $v_1 (t=0) = v_2 (t=0) =0$, the density is constant $\rho (t=0)=1$, the atmospheric pressure is $p_0 = 10^{-5}$, and we insert a quantity of dimensionless energy $e=1$ into a small region of radius $r_0=0.21875$ at the center of the grid,
\begin{align} \label{eq:Sedov_IniCond}
p (t=0) = 
\begin{cases}
p_0 & \mathrm{if}~r\ge r_0, \\
\frac{(\gamma - 1) e}{\pi r_0^2} & \mathrm{otherwise,}
\end{cases}
\end{align}
with $r=\sqrt{x^2+y^2}$.

We tessellate the simulation domain, $\Omega=[-2,2]^2$, with $64 \times 64$ quadrilateral elements, use periodic boundary conditions (however the final time is small enough, such that this does not matter), run the simulations with the split-form DGSEM and the entropy-conserving and kinetic energy preserving flux of Chandrashekar~\cite{Chandrashekar2013} for the volume numerical fluxes, represent the solution with polynomials of degree $N=3$, run the simulation until $t=3$, and use CFL=0.9.

We test different variants of the MCL limiter with local bounds to identify which variants of the MCL limiter can handle the strong shocks of this test, and visualize the nodal limiting factors computed with \eqref{eq:nodalAlpha}.

The first variant of the MCL limiter that we test (from now on referred to as \textbf{MCL limiter A}) uses the standard MCL sequential limiter \eqref{eq:boundsSequential} first, then the sharp positivity limiter with global bounds \eqref{eq:boundsPositivity}, and the semi-discrete entropy limiter \eqref{eq:alpha_s_sync} at last.

\begin{figure}
    \begin{minipage}{0.38\textwidth}
    \includegraphics[trim=460 729 460 0 ,clip,width=\linewidth]{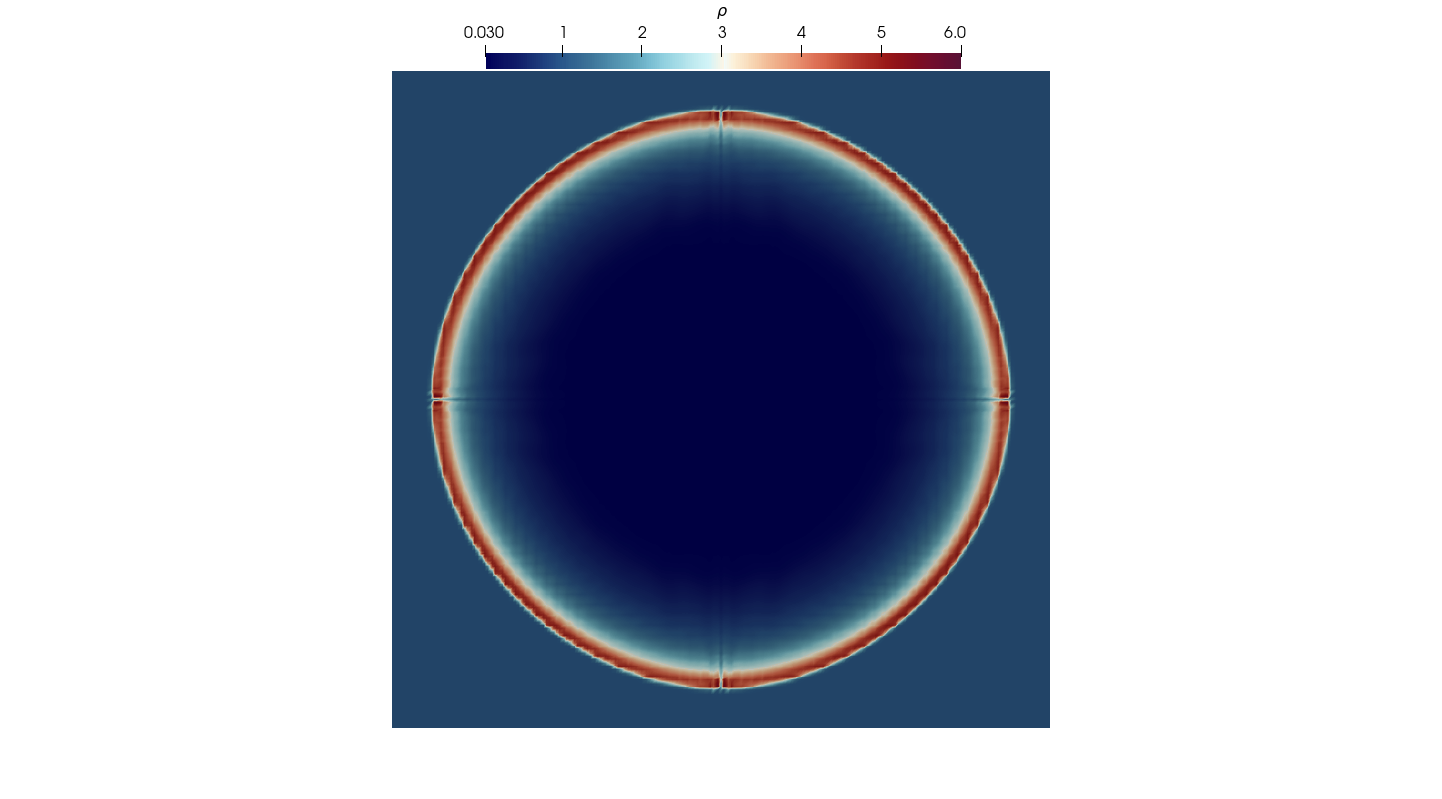}
    
    \includegraphics[trim=390 69 390 69 ,clip,width=0.96\linewidth]{figures/sedov/dens_seq_rho.png}
    \end{minipage}
     \begin{minipage}{0.57\textwidth}
     \centering
     \includegraphics[trim=460 729 460 0 ,clip,width=0.67\linewidth]{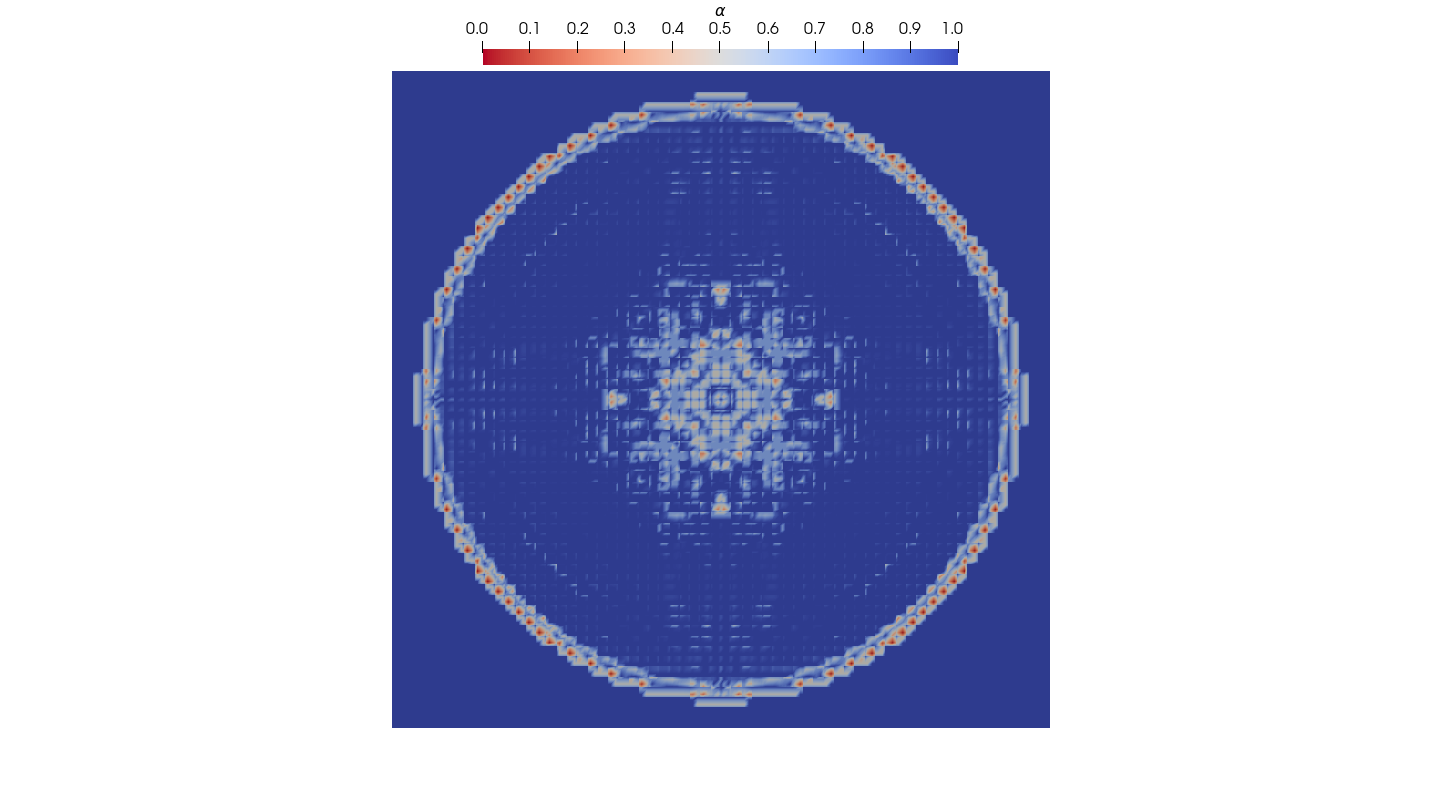}
     
     \begin{overpic}[trim=390 69 390 69 ,clip,width=0.32\textwidth]{figures/sedov/dens_seq_alpha_rho.png}
     \put (0,84) {\color{white}$\Delta \numflux{f}^{\rho}$}
     \end{overpic}
     \begin{overpic}[trim=390 69 390 69 ,clip,width=0.32\linewidth]{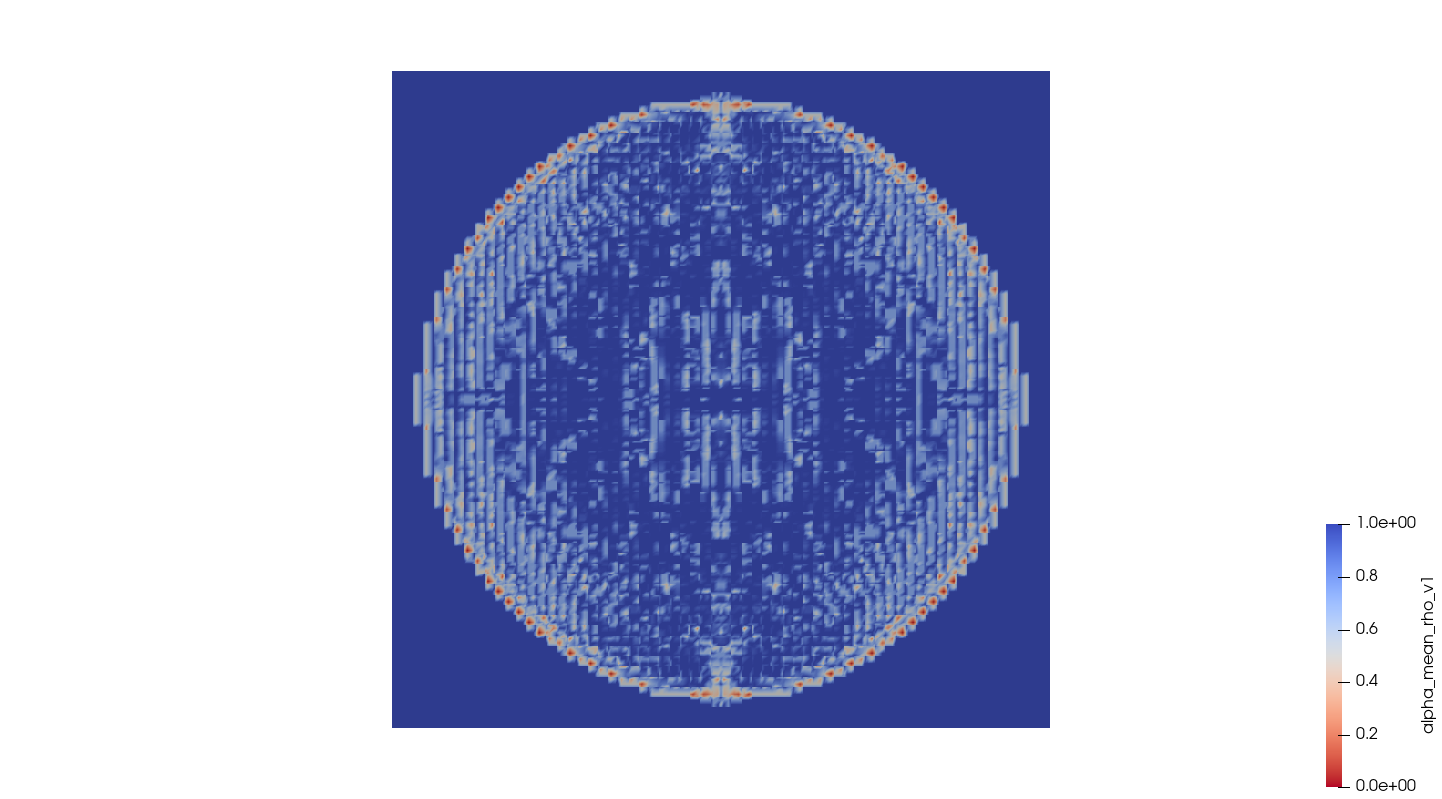}
     \put (0,84) {\color{white}$\Delta \numflux{g}^{v_1}$}
     \end{overpic}
     \begin{overpic}[trim=390 69 390 69 ,clip,width=0.32\linewidth]{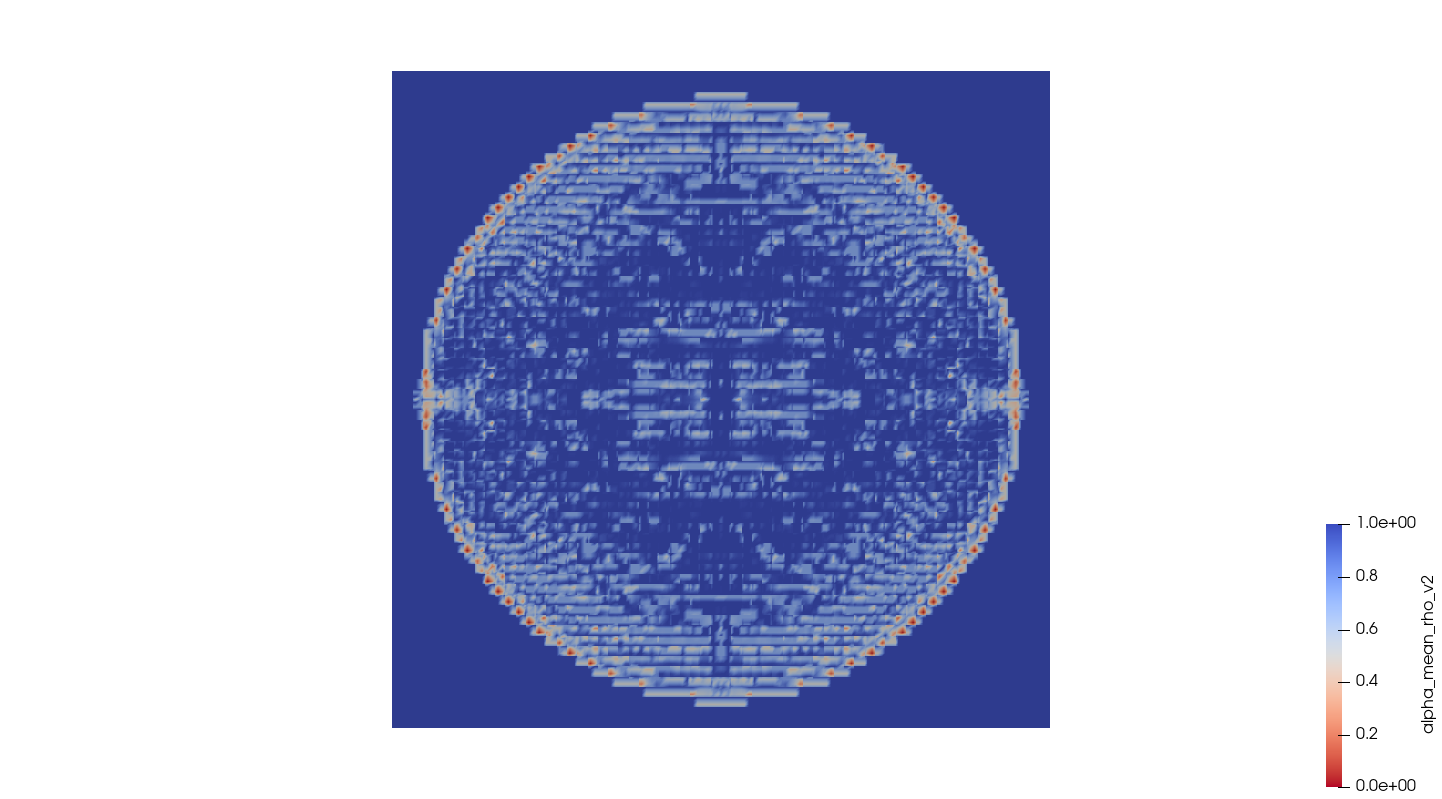}
     \put (0,84) {\color{white}$\Delta \numflux{g}^{v_2}$}
     \end{overpic}
     
     \begin{overpic}[trim=390 69 390 69 ,clip,width=0.32\linewidth]{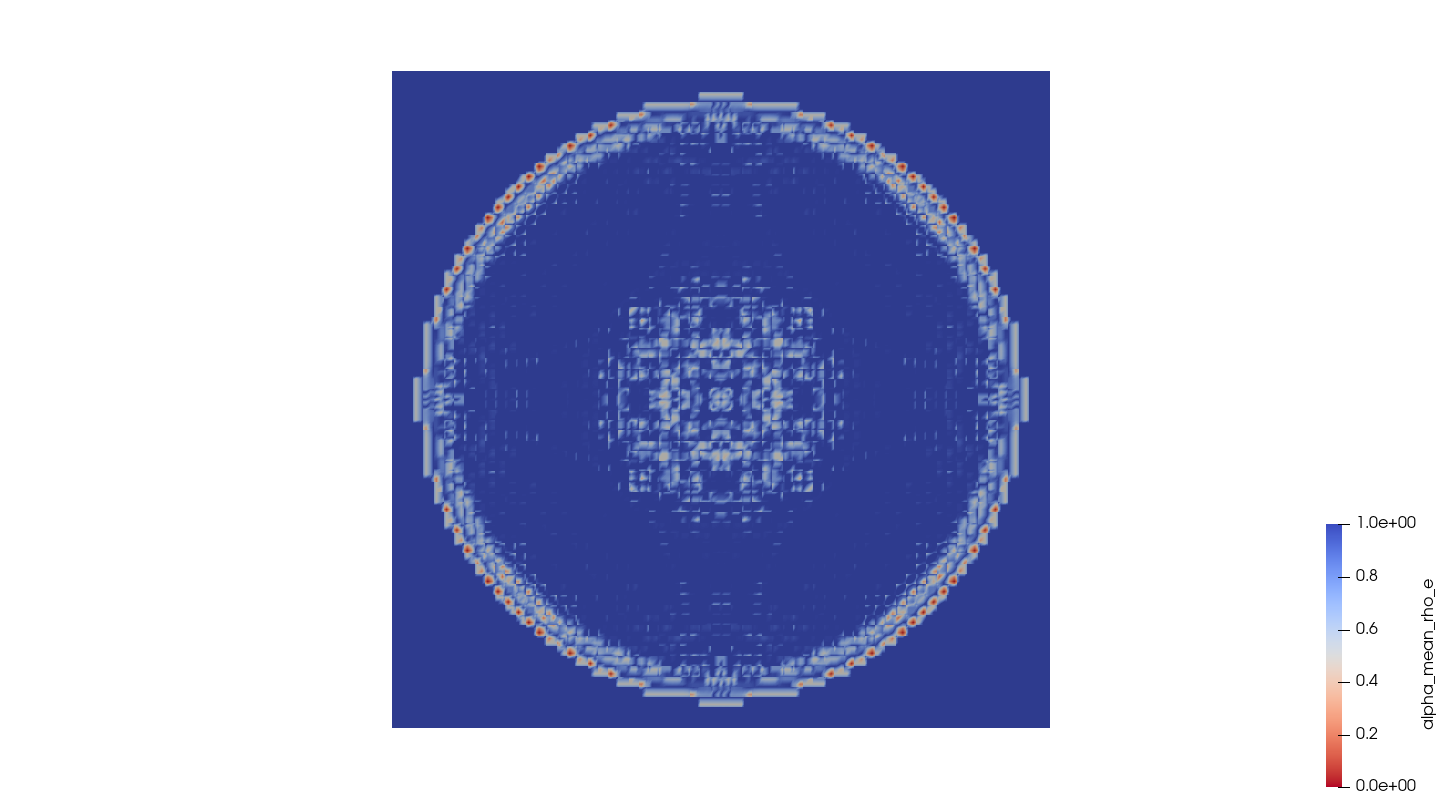}
     \put (0,85) {\color{white}$\Delta \numflux{g}^{E}$}
     \end{overpic}
     \begin{overpic}[trim=390 69 390 69 ,clip,width=0.32\linewidth]{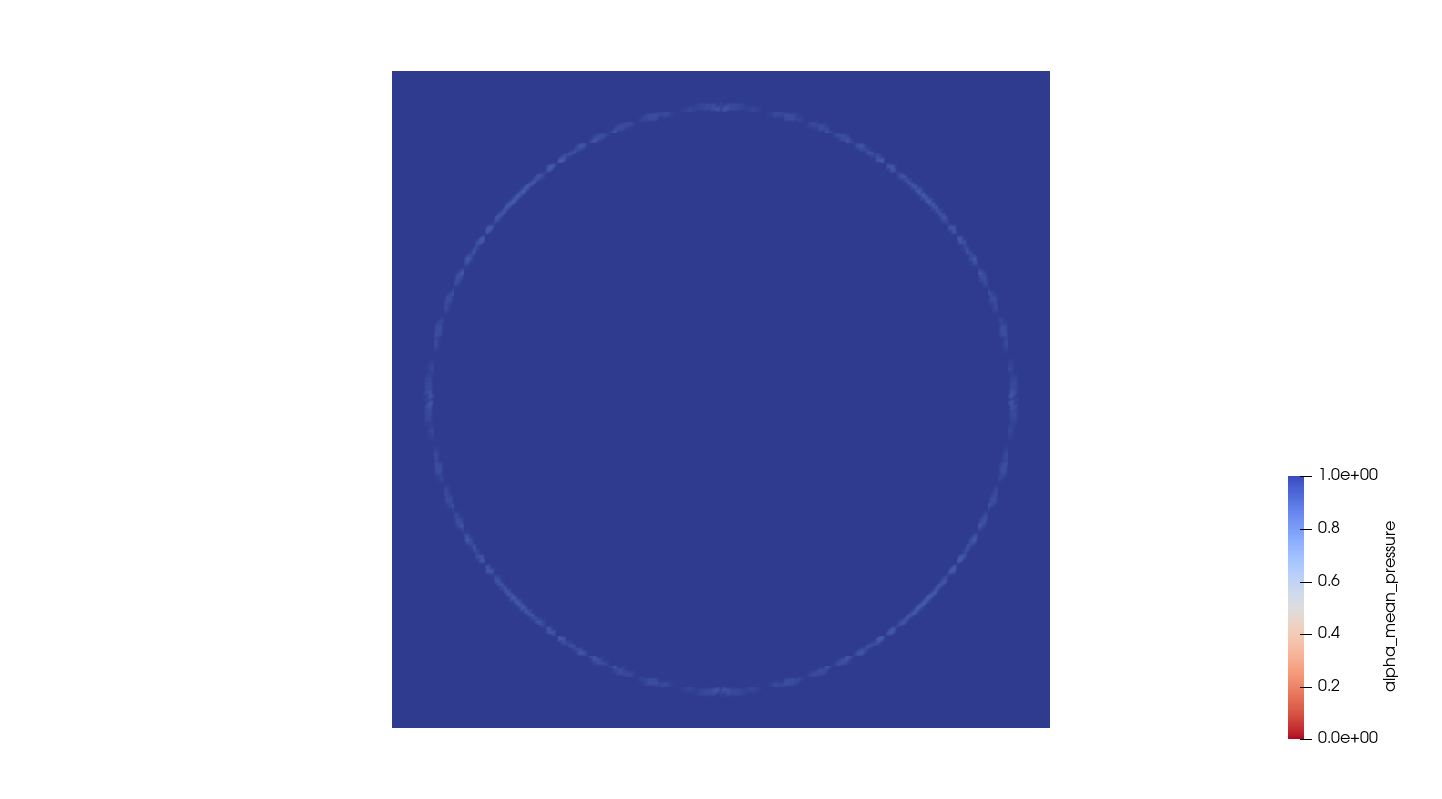}
     \put (0,84) {\color{white}$\Delta \numfluxb{f}^{(p)}$}
     \end{overpic}
     \begin{overpic}[trim=390 69 390 69 ,clip,width=0.32\linewidth]{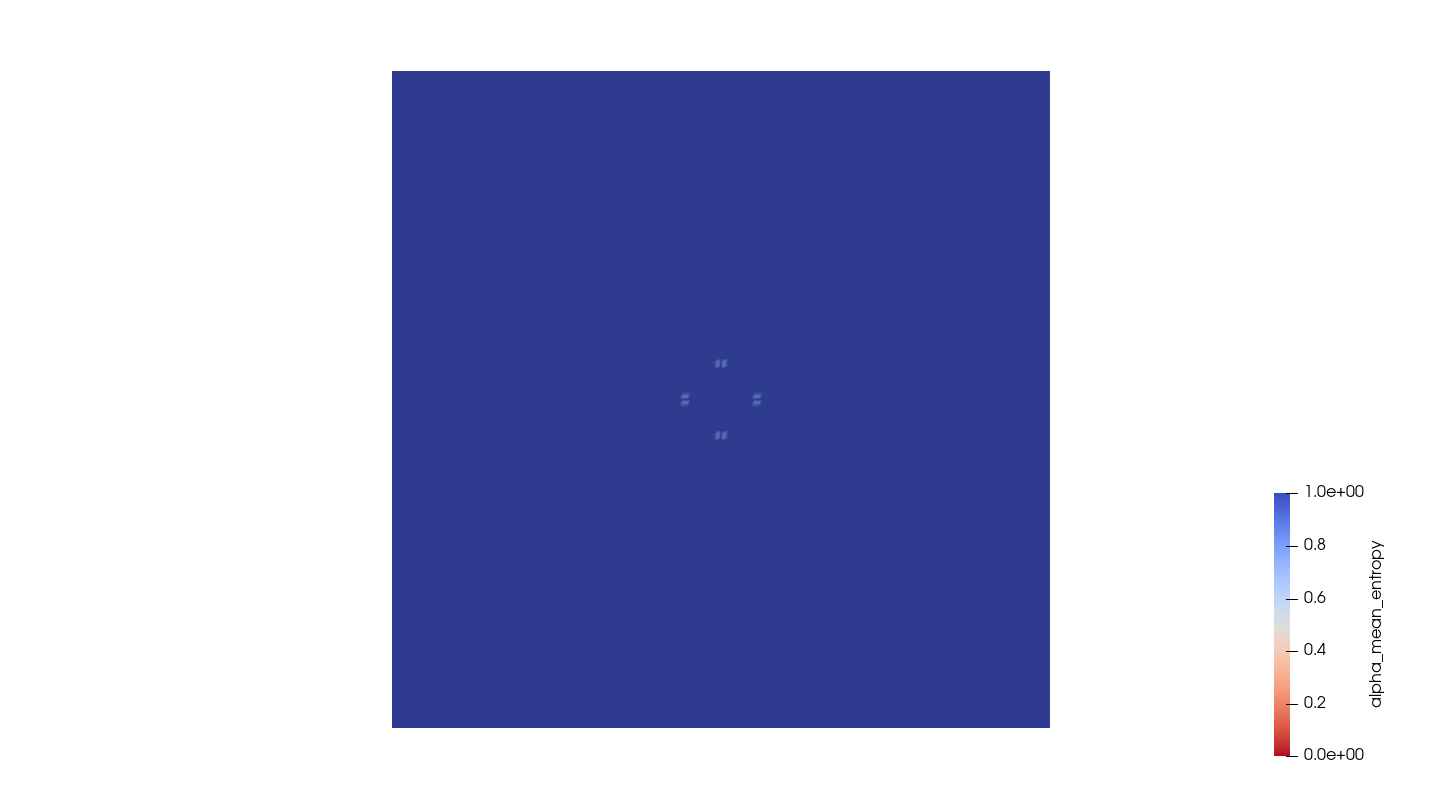}
     \put (0,84) {\color{white}$\Delta \numfluxb{f}^{(ds/dt)}$}
     \end{overpic}
     \end{minipage}
    
    \caption{Contours of the density and limiting factors for the DGSEM simulation with MCL limiter \textbf{A} (sequential limiting, pressure and entropy limiters) of the Sedov blast test at $t=3$.
    The text on the limiting factor plots indicates to which flux the factors are applied.}
    \label{fig:sedov_seq}
\end{figure}

Figure~\ref{fig:sedov_seq} shows the density contours and limiting factors for the Sedov blast simulation using the MCL limiter \textbf{A} at the final time $t=3$.
The text on the limiting factor plots indicates to which flux the factors are applied.
Since both the pressure positivity limiter and the semi-discrete entropy limiter act on all components of the anti-diffusive flux, $\Delta \numfluxb{f}$, we indicate in brackets if the limiting factor corresponds to the pressure ($p$) or the entropy ($ds/dt$) limiters.
Some artifacts at the shock front of the blast wave can be observed with this standard version of the MCL limiter.

The second variant of the MCL limiter that we test (from now on referred to as \textbf{MCL limiter B}) uses the MCL limiter for conservative quantities with local bar-state bounds for the density,
\begin{equation} \label{eq:boundsDensity}
    \min_{j \in \NN (i)} \overline{\rho}_{(i,j)}
    \le \overline{\rho}_{(i,j)}^{Lim} \le
    \max_{j \in \NN (i)} \overline{\rho}_{(i,j)},
\end{equation}
then computes the effective limiting factor for density,
\begin{equation} \label{eq:alphaDensityForAll}
    \alpha^{\rho}_{(i,j)} =
    \begin{cases}
    1 & \mathrm{if}~\Delta \numflux{f}^{\rho}_{(i,j)} \approx 0 \\
    \frac{\Delta \numflux{f}^{\rho,Lim}_{(i,j)} + \epsilon \, \texttt{sign}(\Delta \numflux{f}^{\rho}_{(ij,k)})}
    {\Delta \numflux{f}^{\rho}_{(i,j)} + \epsilon \, \texttt{sign}(\Delta \numflux{f}^{\rho}_{(ij,k)})} & \mathrm{otherwise},
    \end{cases}
\end{equation}
where $\epsilon$ is a very small number, applies $\alpha^{\rho}_{(i,j)}$ to all other conservative quantities ($\Delta \numflux{f}^{\rho v_1}$, $\Delta \numflux{f}^{\rho v_2}$, and $\Delta \numflux{f}^{\rho E}$), then uses the sharp positivity limiter with global bounds \eqref{eq:boundsPositivity}, and the semi-discrete entropy limiter \eqref{eq:alpha_s_sync} at last.

\begin{figure}
    \begin{minipage}{0.38\textwidth}
    \includegraphics[trim=460 729 460 0 ,clip,width=\linewidth]{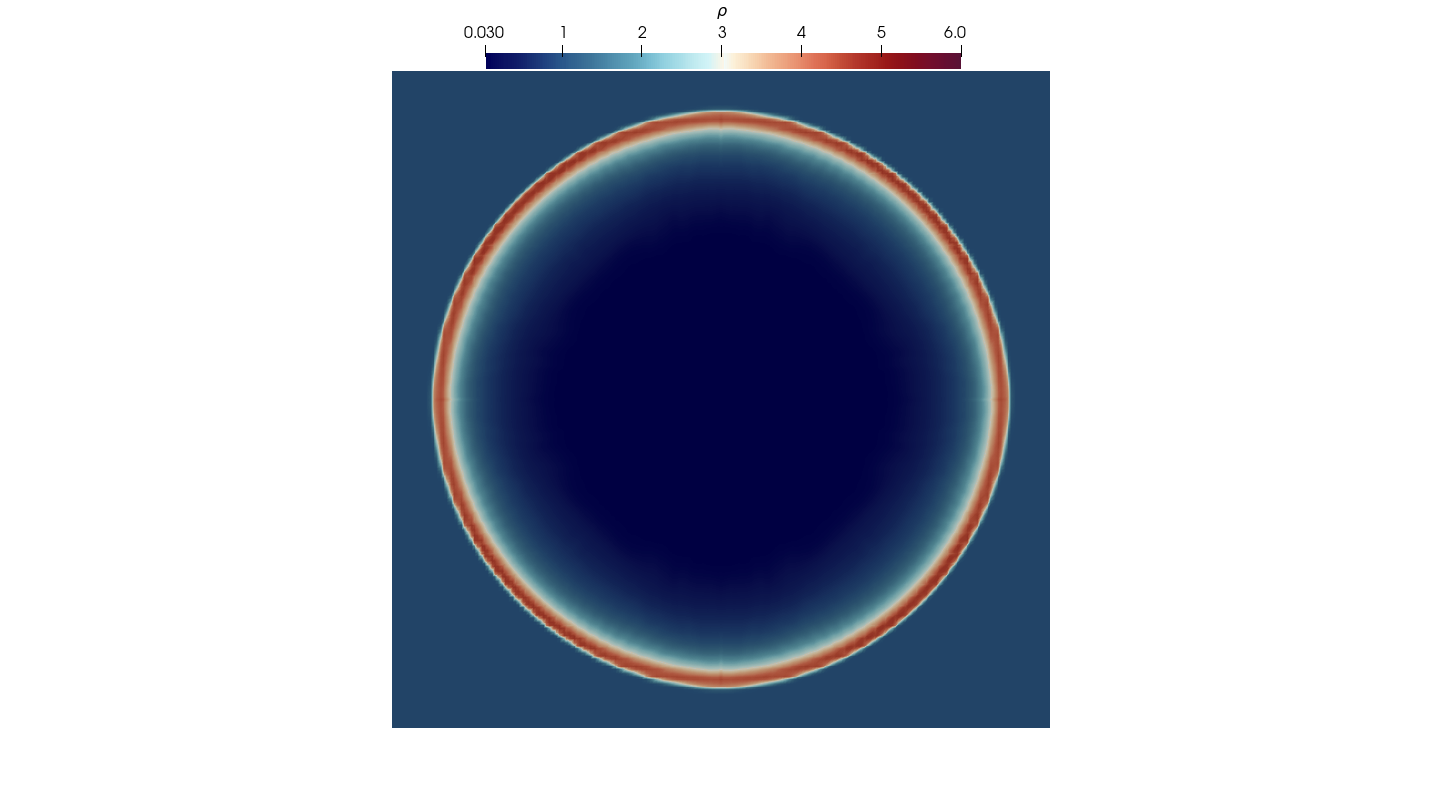}
    
    \includegraphics[trim=390 69 390 69 ,clip,width=0.96\linewidth]{figures/sedov/densforall_rho.png}
    \end{minipage}
     \begin{minipage}{0.57\textwidth}
     \centering
     \includegraphics[trim=460 729 460 0 ,clip,width=0.67\linewidth]{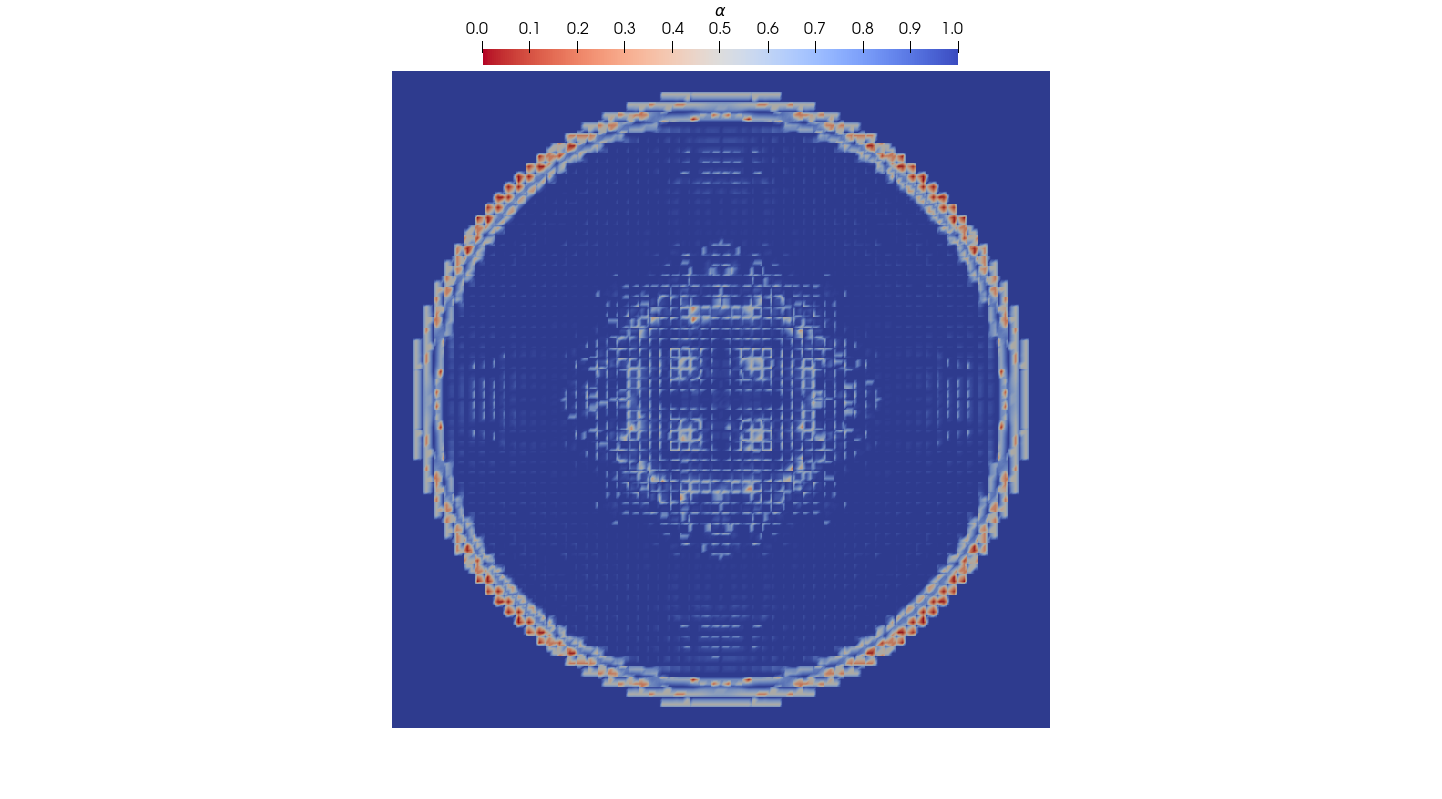}
     
     \begin{overpic}[trim=390 69 390 69 ,clip,width=0.32\textwidth]{figures/sedov/densforall_alpha_rho.png}
     \put (0,84) {\color{white}$\Delta \numflux{f}^{\rho}$}
     \end{overpic}
     \begin{overpic}[trim=390 69 390 69 ,clip,width=0.32\linewidth]{figures/sedov/densforall_alpha_rho.png}
     \put (0,84) {\color{white}$\Delta \numflux{f}^{\rho v_1}$}
     \end{overpic}
     \begin{overpic}[trim=390 69 390 69 ,clip,width=0.32\linewidth]{figures/sedov/densforall_alpha_rho.png}
     \put (0,84) {\color{white}$\Delta \numflux{f}^{\rho v_2}$}
     \end{overpic}
     
     \begin{overpic}[trim=390 69 390 69 ,clip,width=0.32\linewidth]{figures/sedov/densforall_alpha_rho.png}
     \put (0,84) {\color{white}$\Delta \numflux{f}^{\rho E}$}
     \end{overpic}
     \begin{overpic}[trim=390 69 390 69 ,clip,width=0.32\linewidth]{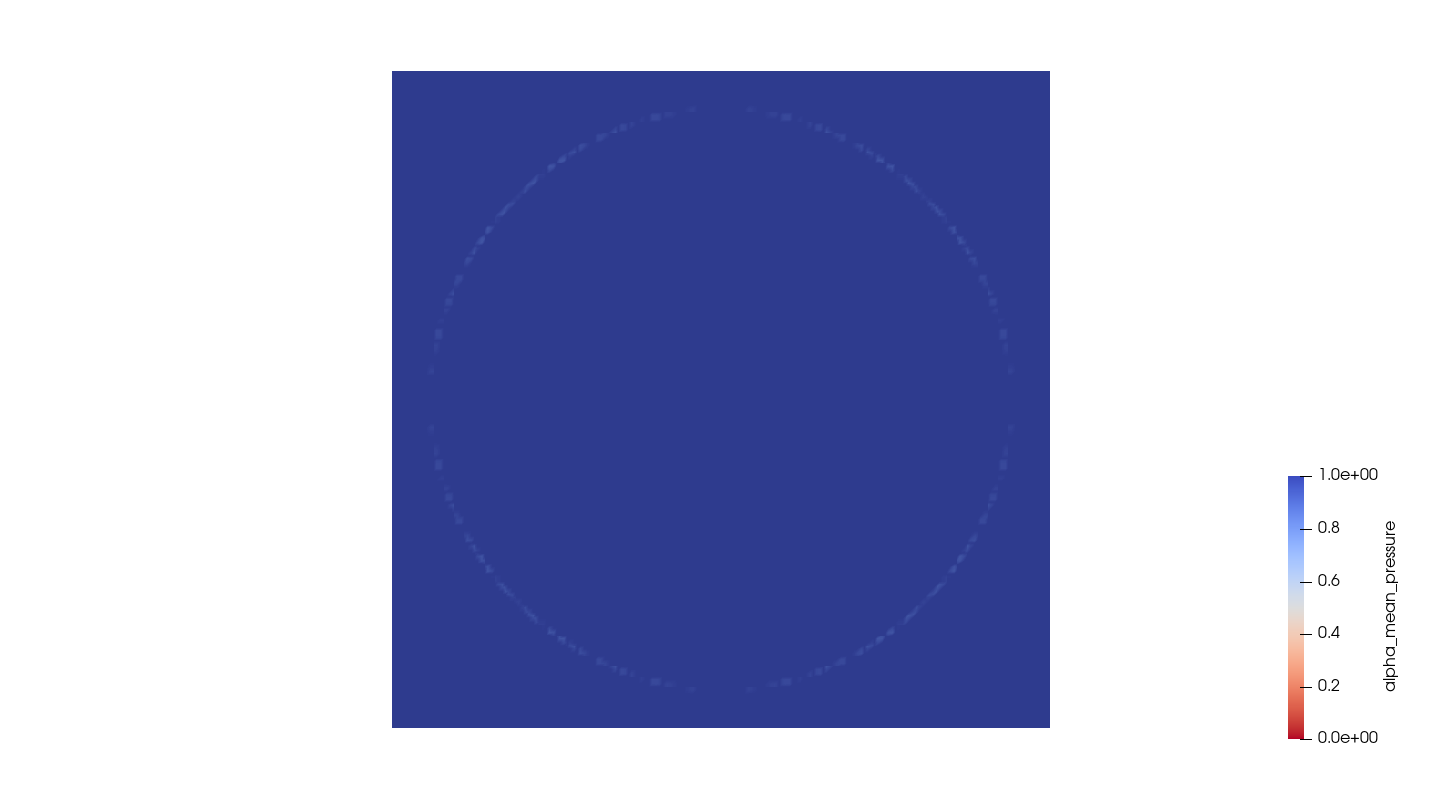}
     \put (0,84) {\color{white}$\Delta  \numfluxb{f}^{(p)}$}
     \end{overpic}
     \begin{overpic}[trim=390 69 390 69 ,clip,width=0.32\linewidth]{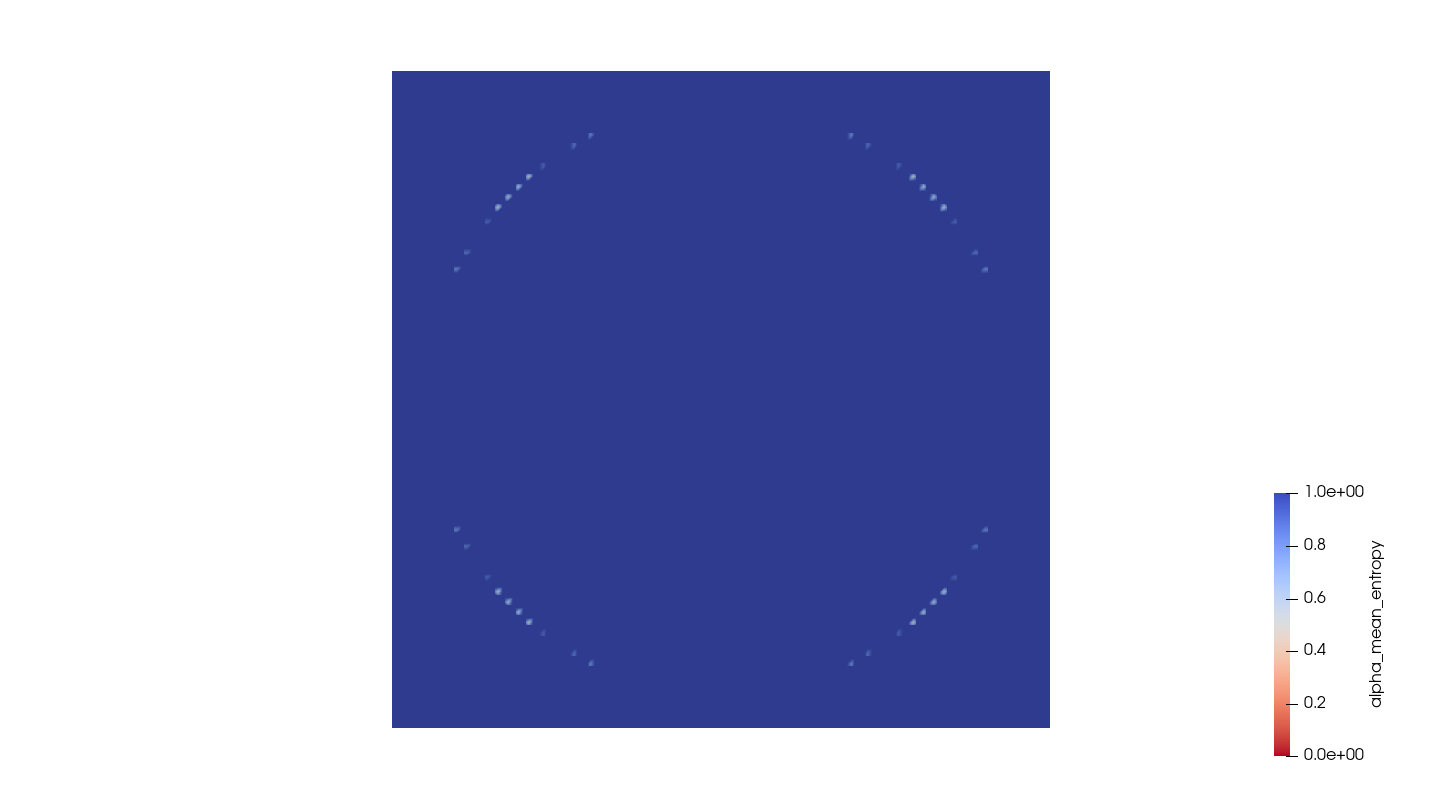}
     \put (0,84) {\color{white}$\Delta \numfluxb{f}^{(ds/dt)}$}
     \end{overpic}
     \end{minipage}
    
    \caption{Contours of the density and limiting factors for the DGSEM simulation with MCL limiter \textbf{B} (density factor for all equations, pressure and entropy limiters) of the Sedov blast test at $t=3$.
    The text on the limiting factor plots indicates to which flux the factors are applied.}
    \label{fig:sedov_densforall}
\end{figure}

Figure~\ref{fig:sedov_densforall} shows the density contours and limiting factors for the Sedov blast simulation using the MCL limiter \textbf{B} at the final time $t=3$.
In this case, the density, pressure, and the semi-discrete entropy limiter act on all components of the anti-diffusive flux.
Note that the pressure positivity limiter needs to act less, but the semi-discrete entropy limiter needs to act more than in the MCL limiter \textbf{A}.
Moreover, the MCL limiter \textbf{B} removes the Carbuncle-like artefacts at the shock fronts of the blast.

The last variant of the MCL limiter that we test (from now on referred to as \textbf{MCL limiter C}) is a combination of MCL limiters \textbf{A} and \textbf{B}.
It uses the MCL limiter for conservative quantities with local bar-state bounds for the density \eqref{eq:boundsDensity}, computes the effective limiting coefficient using \eqref{eq:alphaDensityForAll} and applies it to the other conservative quantities, then uses the sequential limiter to impose the bounds on ``primitive'' quantities \eqref{eq:boundsSequential}, then imposes global bounds on density and pressure \eqref{eq:boundsPositivity} with the sharp positivity limiter, and then applies the semi-discrete entropy limiter \eqref{eq:alpha_s_sync} at last.

\begin{figure}
    \begin{minipage}{0.38\textwidth}
    \includegraphics[trim=460 729 460 0 ,clip,width=\linewidth]{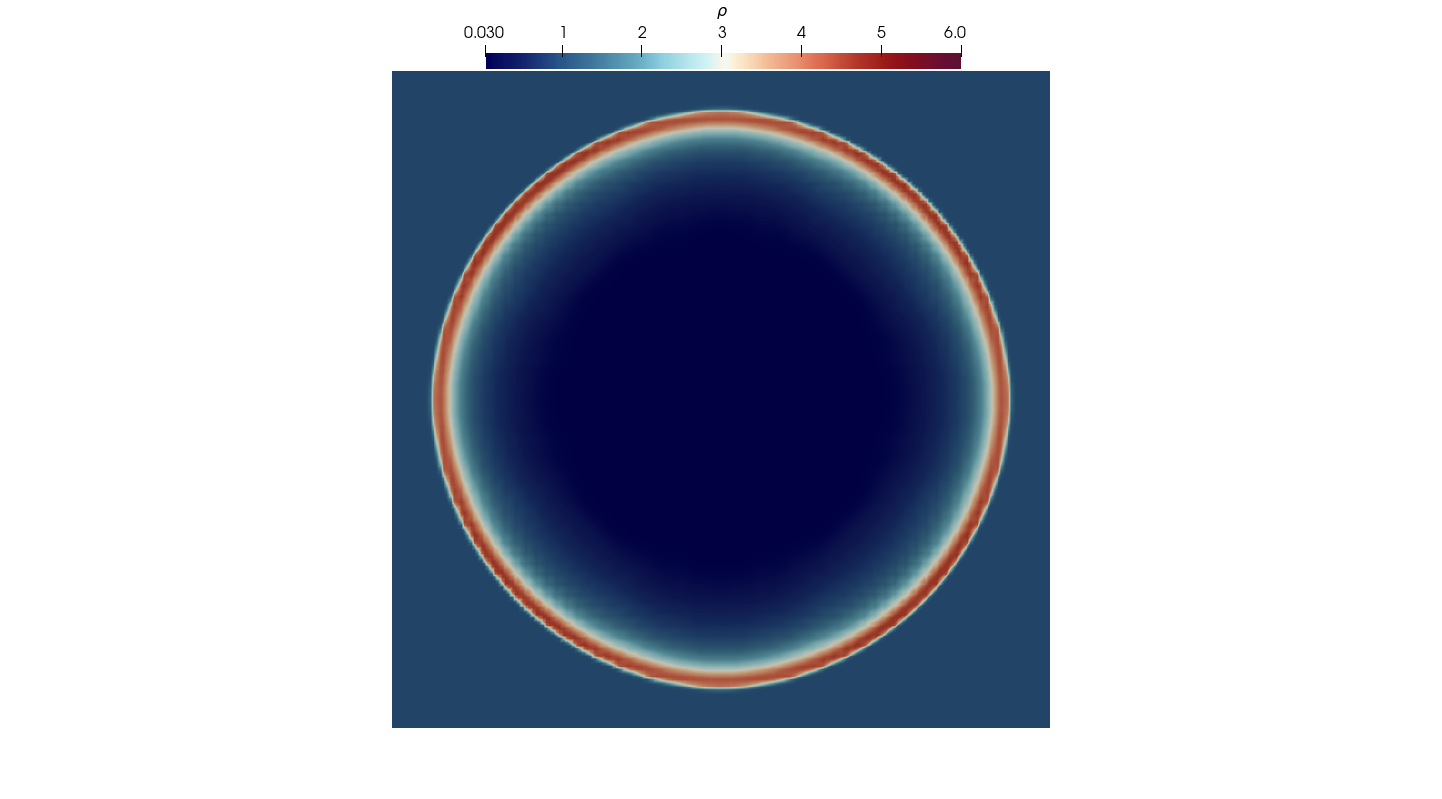}
    
    \includegraphics[trim=390 69 390 69 ,clip,width=0.96\linewidth]{figures/sedov/densforall_seq_rho.png}
    \end{minipage}
     \begin{minipage}{0.57\textwidth}
     \centering
     \includegraphics[trim=460 729 460 0 ,clip,width=0.67\linewidth]{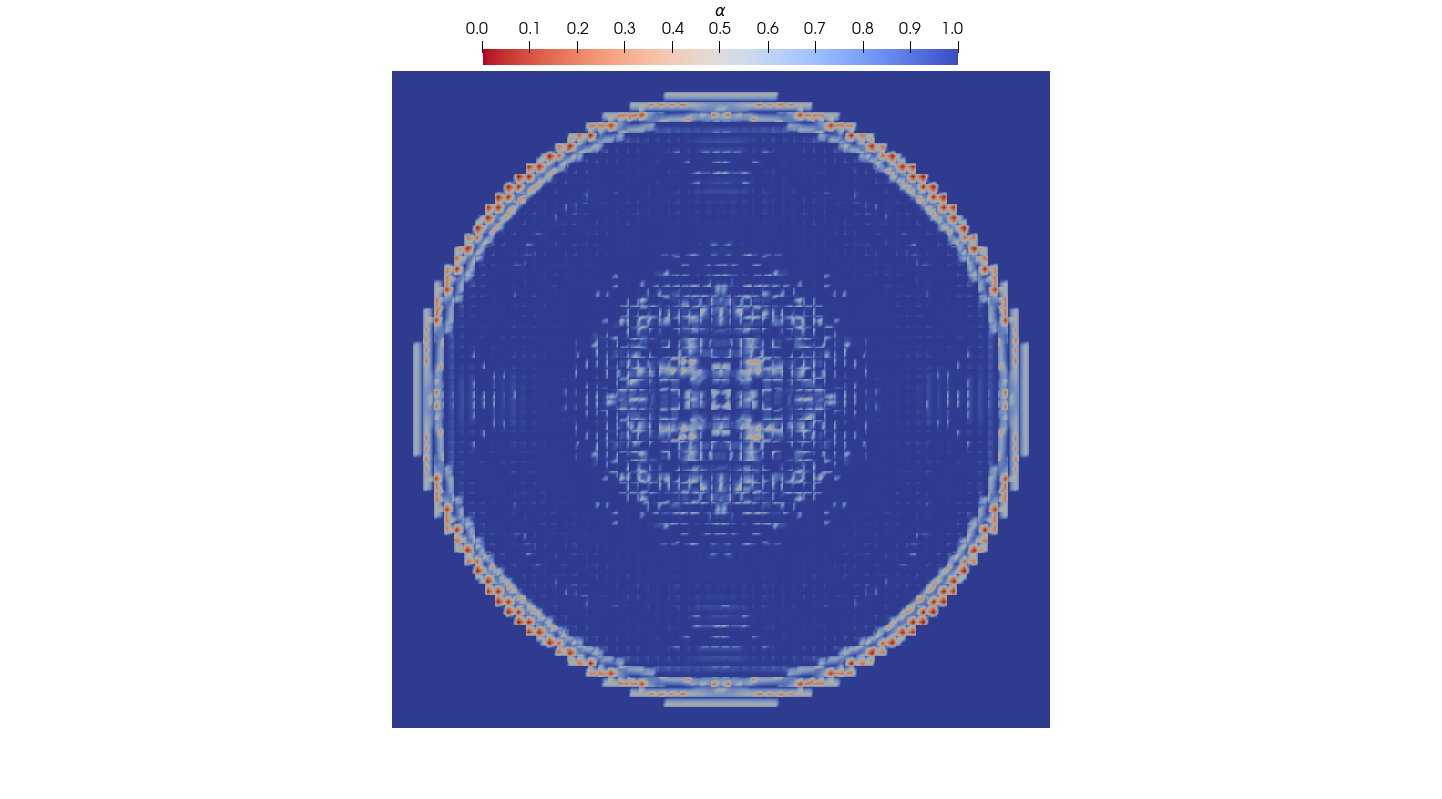}
     
     \begin{overpic}[trim=390 69 390 69 ,clip,width=0.32\textwidth]{figures/sedov/densforall_seq_alpha_rho.png}
     \put (0,84) {\color{white}$\Delta \numfluxb{f}^{(\rho)}$}
     \end{overpic}
     \begin{overpic}[trim=390 69 390 69 ,clip,width=0.32\linewidth]{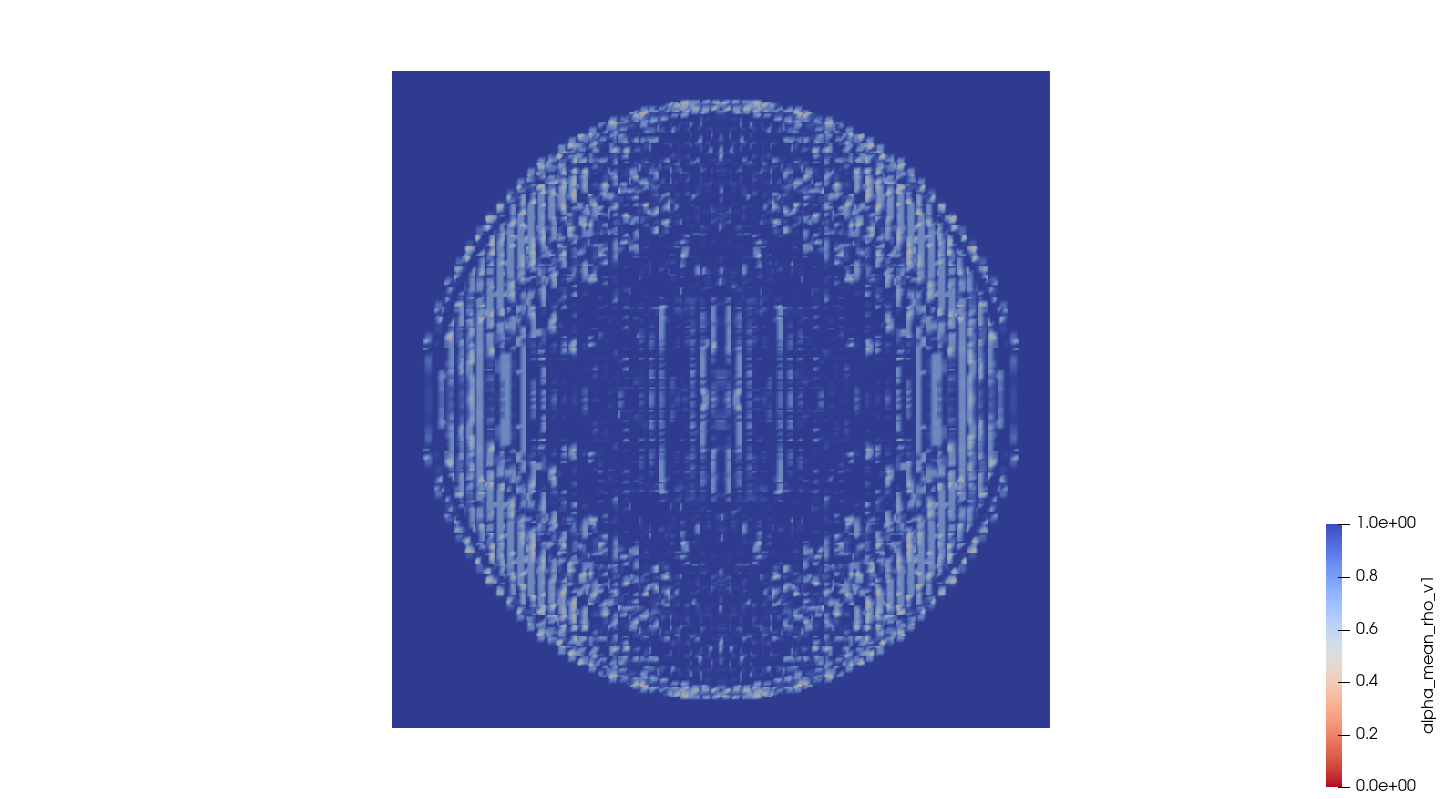}
     \put (0,84) {\color{white}$\Delta \numflux{g}^{v_1}$}
     \end{overpic}
     \begin{overpic}[trim=390 69 390 69 ,clip,width=0.32\linewidth]{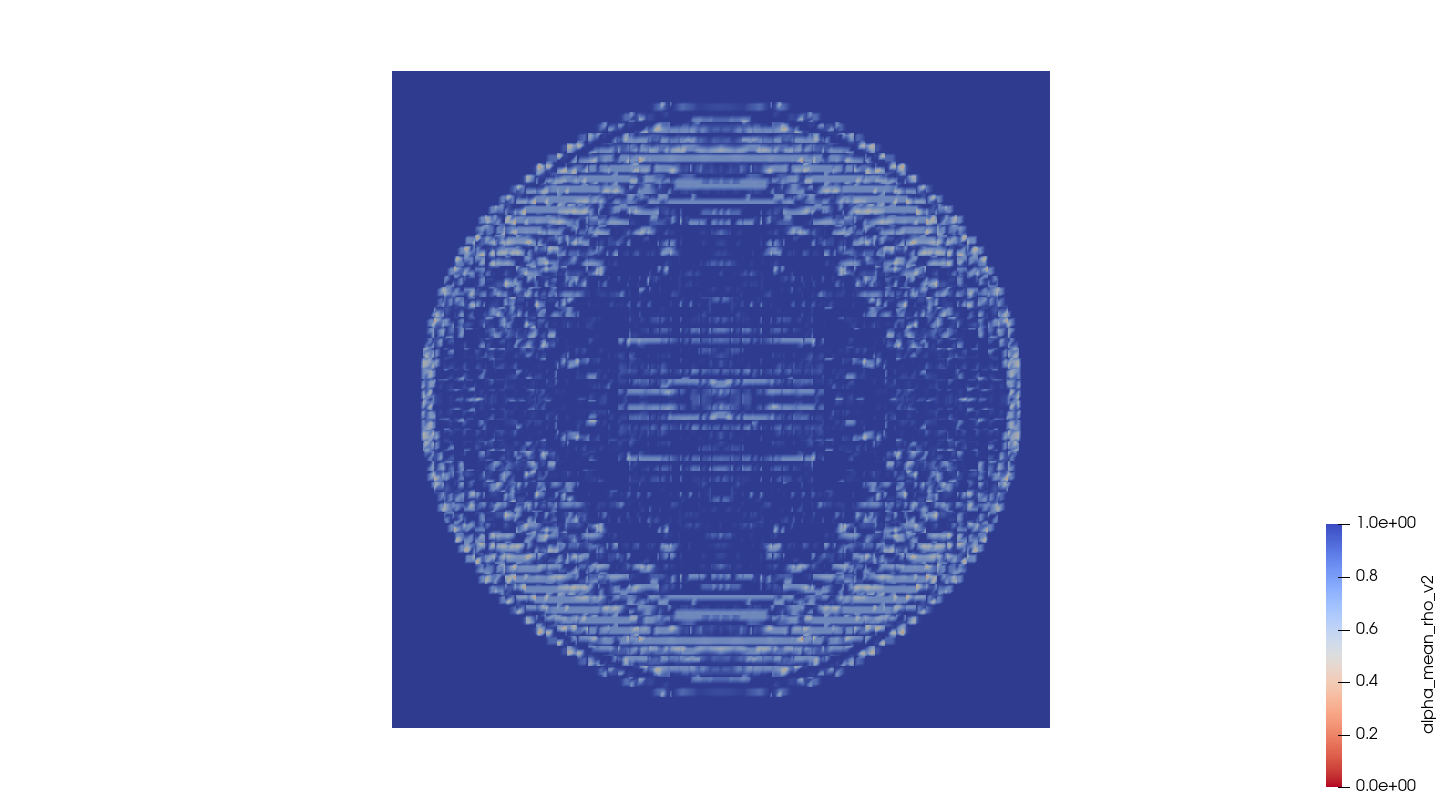}
     \put (0,84) {\color{white}$\Delta \numflux{g}^{v_2}$}
     \end{overpic}
     
     \begin{overpic}[trim=390 69 390 69 ,clip,width=0.32\linewidth]{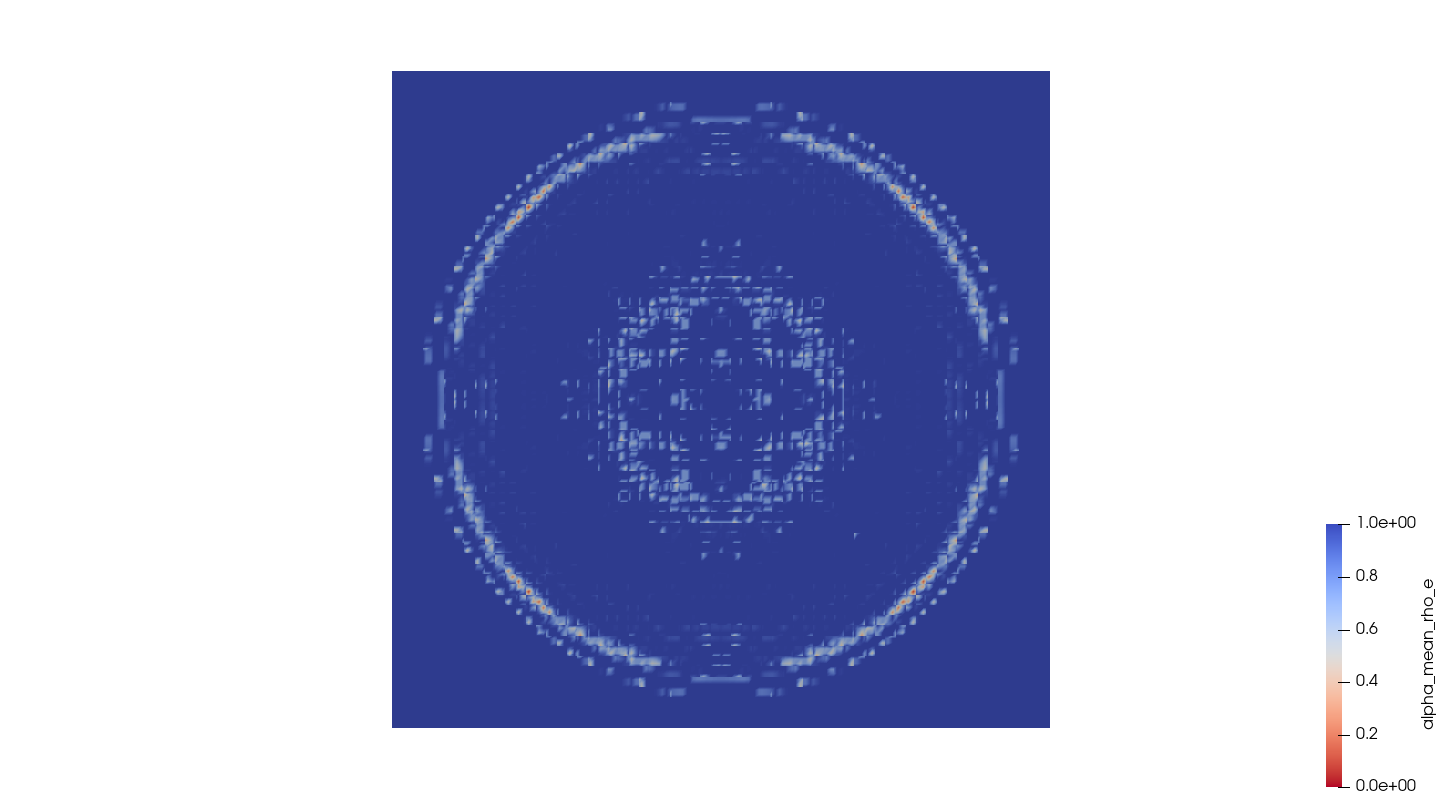}
     \put (0,84) {\color{white}$\Delta \numflux{g}^{E}$}
     \end{overpic}
     \begin{overpic}[trim=390 69 390 69 ,clip,width=0.32\linewidth]{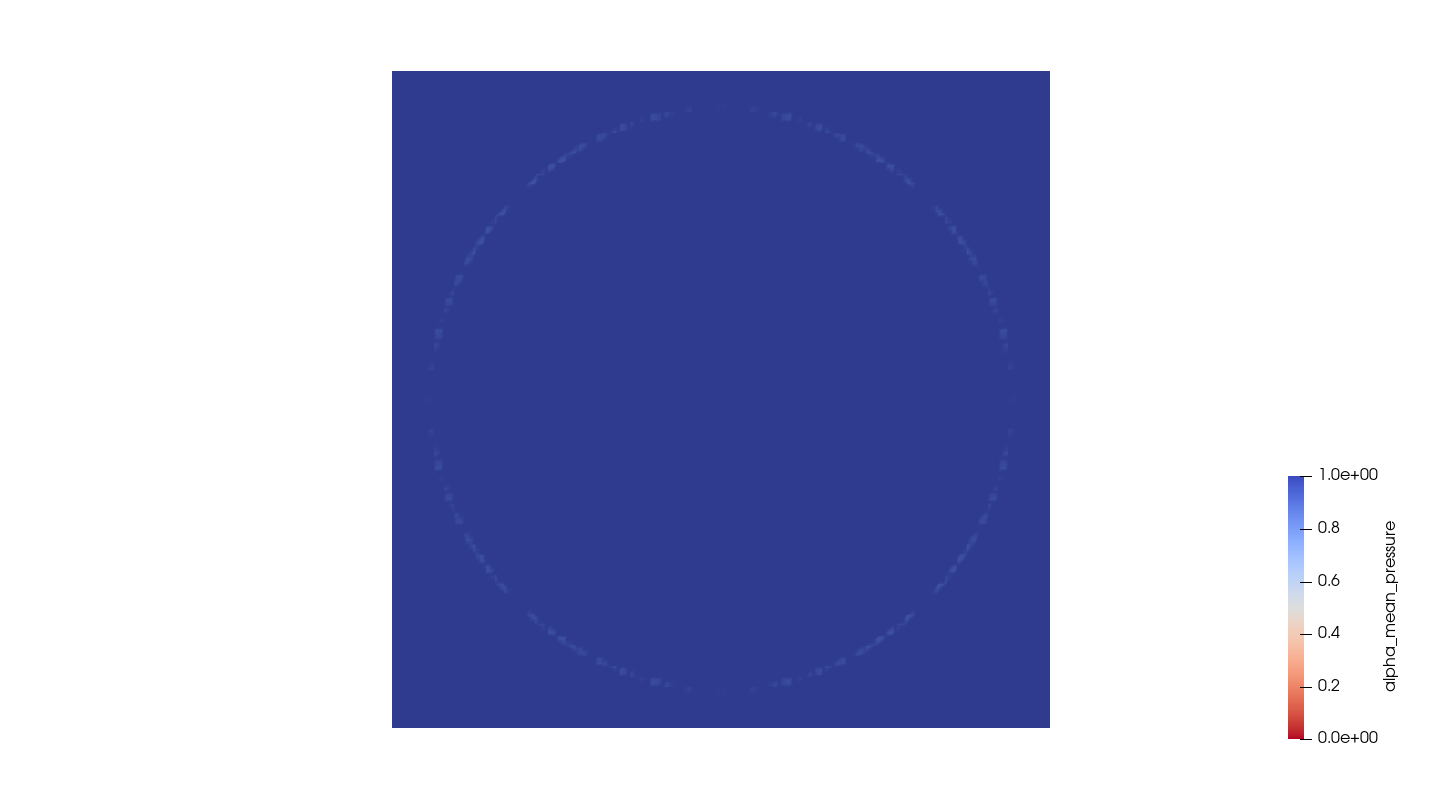}
     \put (0,84) {\color{white}$\Delta \numfluxb{f}^{(p)}$}
     \end{overpic}
     \begin{overpic}[trim=390 69 390 69 ,clip,width=0.32\linewidth]{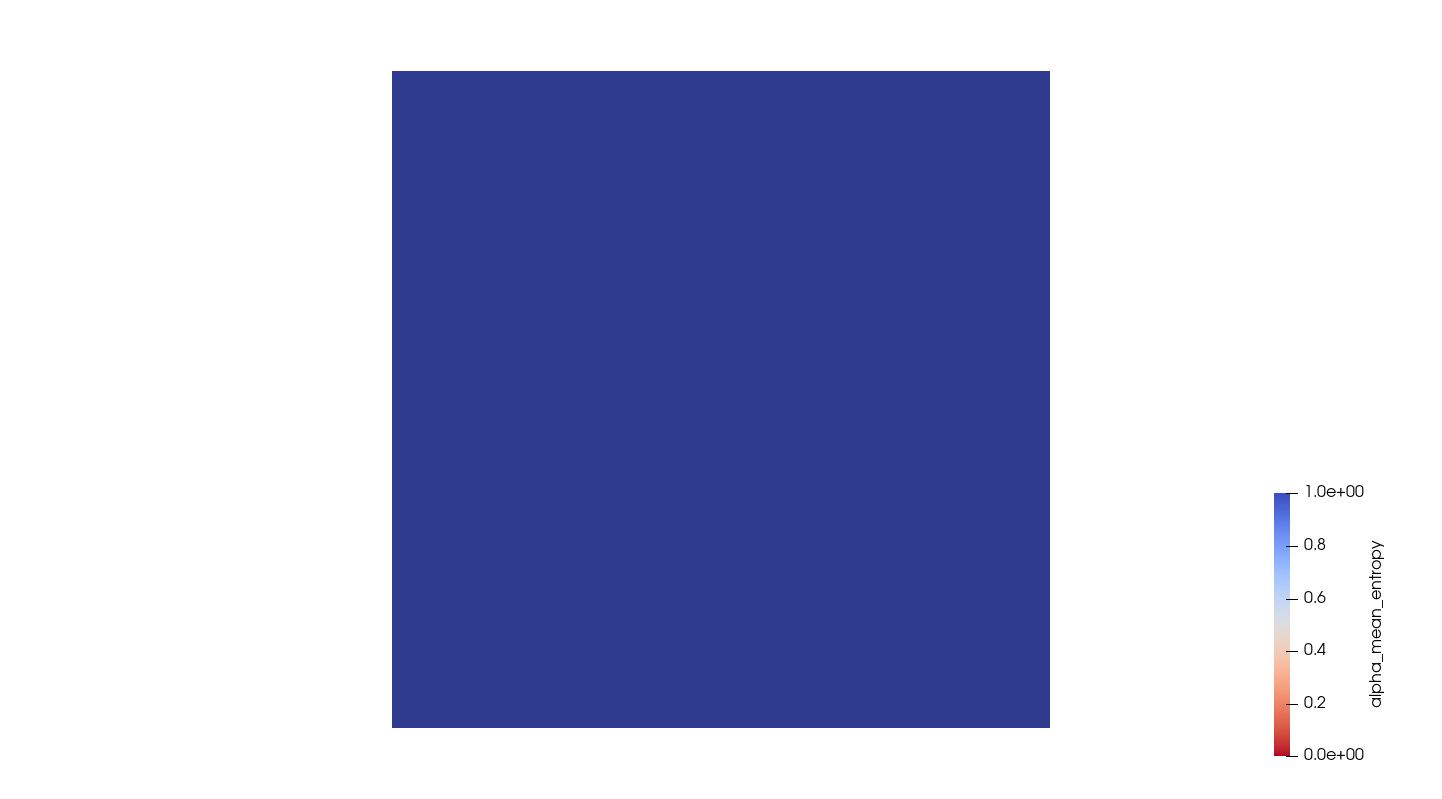}
     \put (0,84) {\color{white}$\Delta \numfluxb{f}^{(ds/dt)}$}
     \end{overpic}
     \end{minipage}
    
    \caption{Contours of the density and limiting factors for the DGSEM simulation with MCL limiter \textbf{C} (density factor for all equations, sequential limiting, pressure and entropy limiters) of the Sedov blast test at $t=3$.
    The text on the limiting factor plots indicates to which flux the factors are applied.}
    \label{fig:sedov_densforall_seq}
\end{figure}

Figure~\ref{fig:sedov_densforall_seq} illustrates the density and limiting factor contours for the Sedov blast simulation using the MCL limiter \textbf{C} at the final time $t=3$.
Again, the density, pressure, and the semi-discrete entropy limiter act on all components of the anti-diffusive flux.
Note that the sequential limiter for the velocities and the total specific energy needs to apply less limiting than in MCL limiter \textbf{A} due to the action of the density limiting coefficient on the momentum and total energy fluxes.
The pressure positivity limiter needs to apply less limiting than in MCL limiters \textbf{A} and \textbf{B}, and the semi-discrete entropy limiter does not need to act at all for the snapshot at $t=3$.
As with MCL limiter \textbf{B}, the Carbuncle-like artefacts are no longer present, but the resulting scheme is clearly more dissipative.

\subsection{High-Mach Astrophysical Jet}

To test the robustness of the MCL techniques, we simulate a setup inspired by an astrophysical jet application with Mach number $\textrm{Ma} \approx 2000$, which was originally proposed by Ha~et~al.~\cite{ha2005numerical}.
This extreme benchmark case has been used to stress-test shock-capturing techniques for high-order methods \cite{zhang2010,liuoscillation,RUEDARAMIREZ2022}.

The computational domain, $\Omega = [-0.5,0.5]^2$, is filled with a mono-atomic gas ($\gamma = 5/3$) at rest with
\begin{equation*}
\rho(x,y) = 0.5, \qquad
p(x,y) = 0.4127, \qquad
v_1(x,y) = 0, \qquad
v_2(x,y) = 0,
\end{equation*}
and on the left boundary there is a hypersonic inflow with
\begin{equation*}
\rho(-0.5,y_B) = 5, ~
p(-0.5,y_B) = 0.4127, ~
v_1(-0.5,y_B) = 800, ~
v_2(-0.5,y_B) = 0,
\end{equation*}
for $y_B \in [-0.05, 0.05]$, which corresponds to a Mach number of $\textrm{Ma}=2156.91$ with respect to the speed of sound of the jet gas, and $\textrm{Ma}=682.08$ with respect to the speed of sound of the ambient gas.

We solve this problem using $256 \times 256$ quadrilateral elements of degree $N=3$, use periodic boundary conditions for the top and bottom boundaries and characteristics-based inflow/outflow boundary conditions for the left and right boundaries, the entropy-conserving and kinetic energy preserving flux of Ranocha \cite{ranocha2018generalised} for the volume numerical flux of the DGSEM method, and different MCL and FCT/IDP limiters and CFL numbers.

\begin{figure}[h!]
\centering
\includegraphics[trim=336 846 336 0 ,clip,width=0.65\linewidth]{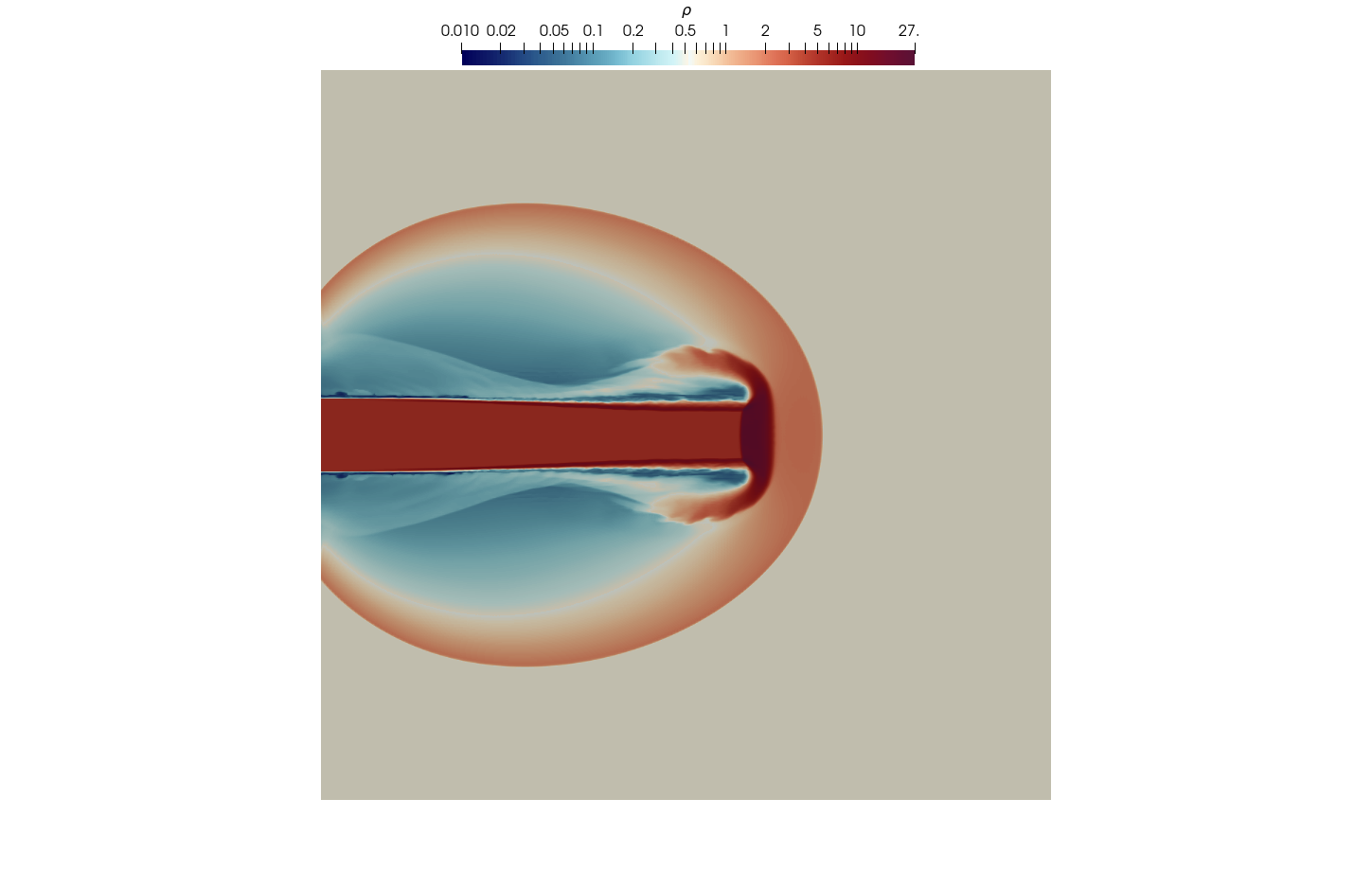}
\\
\begin{minipage}[t]{0.32\textwidth}
\begin{overpic}[trim=336 72 336 72 ,clip,width=\linewidth]{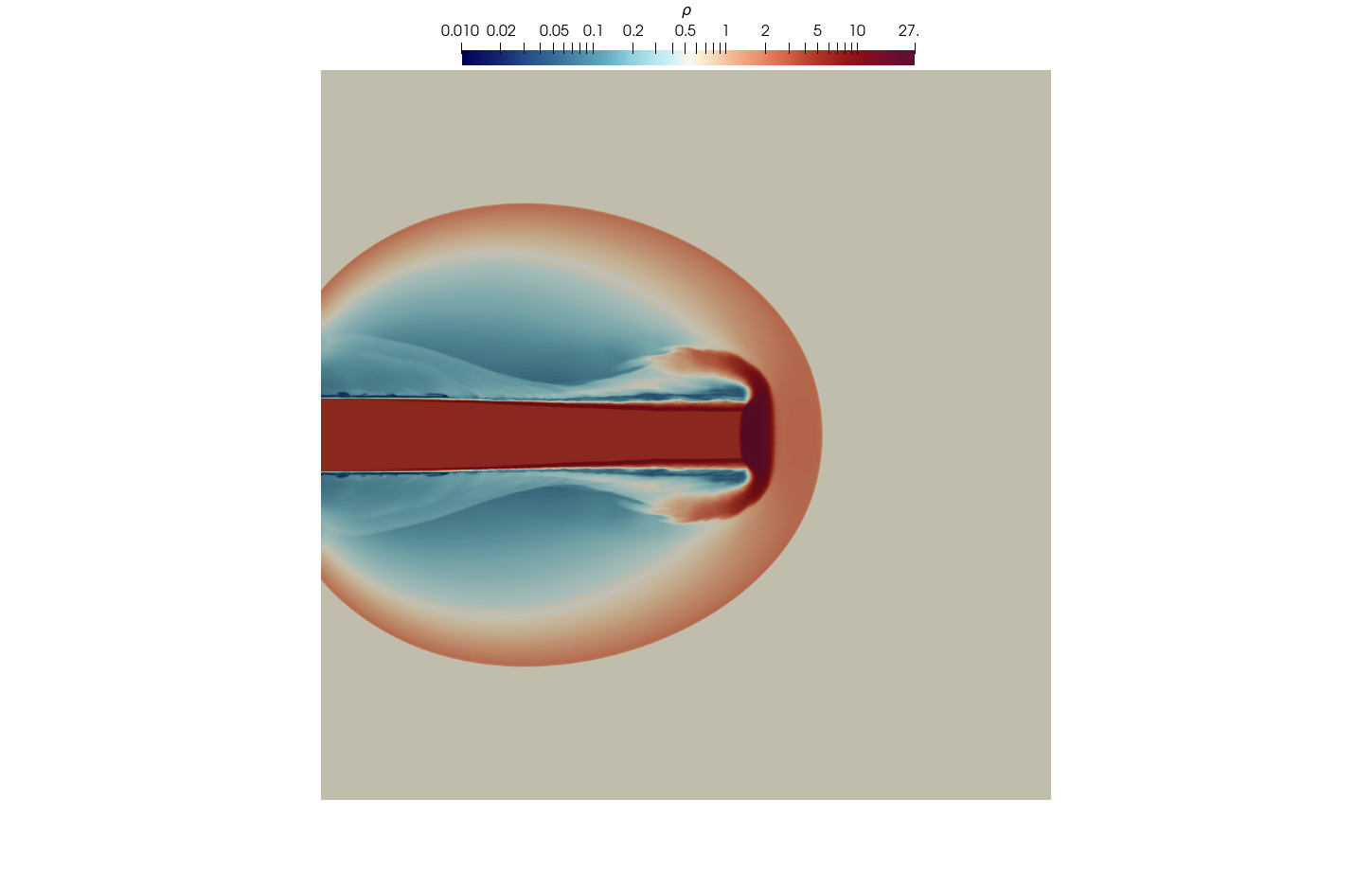}
    \put (3,90) {CFL=0.9}
\end{overpic}
\begin{overpic}[trim=336 72 336 72 ,clip,width=\linewidth]{figures/astrojet/mcl_densityforall_seq_presExact_cfl_0.45.png}
     \put (3,90) {CFL=0.45}
\end{overpic}
\begin{overpic}[trim=389 70 389 70 ,clip,width=\linewidth]{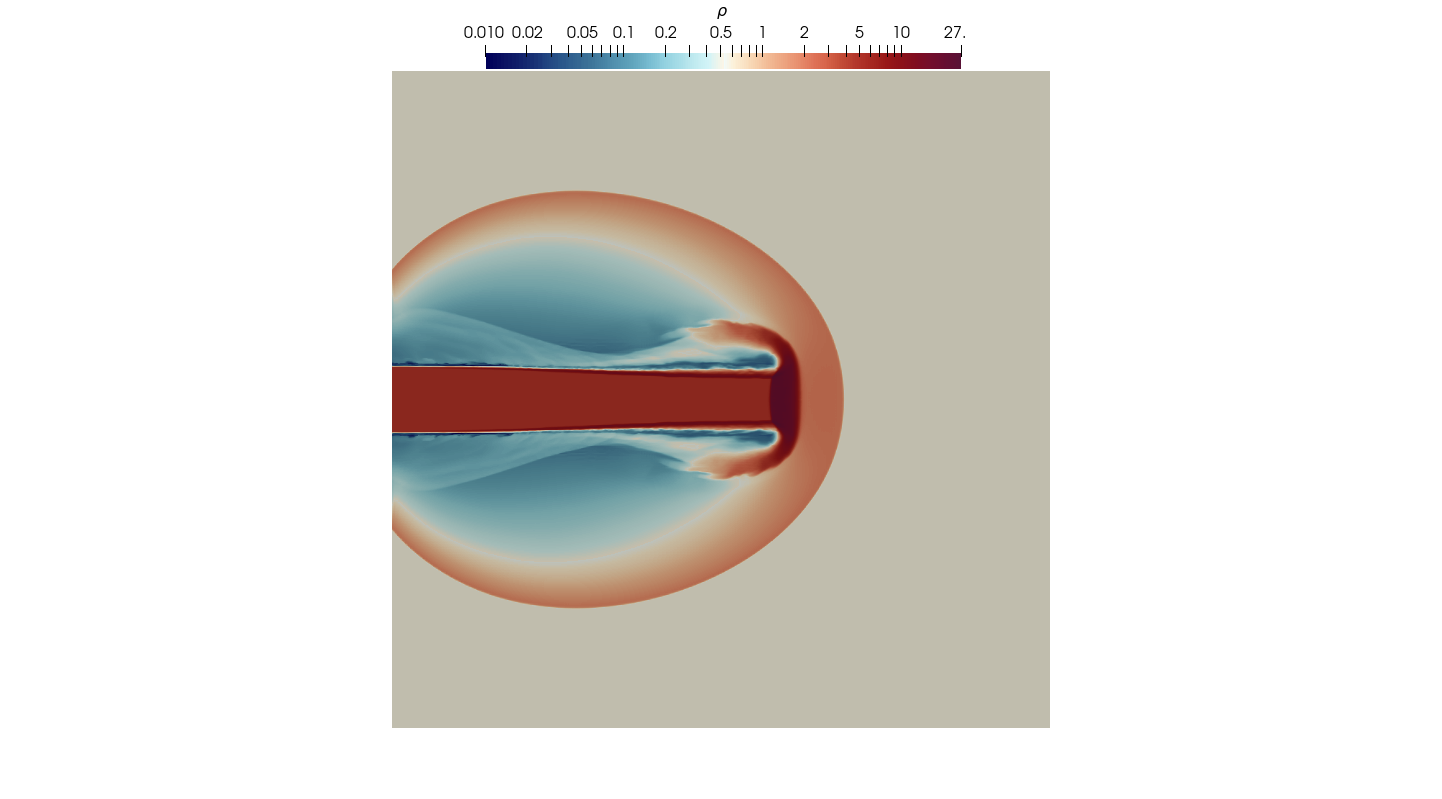}
	     \put (3,90) {CFL=0.225}
\end{overpic}
\begin{subfigure}[b]{\linewidth}
    \centering
	\begin{overpic}[trim=333 72 333 72 ,clip,width=\linewidth]{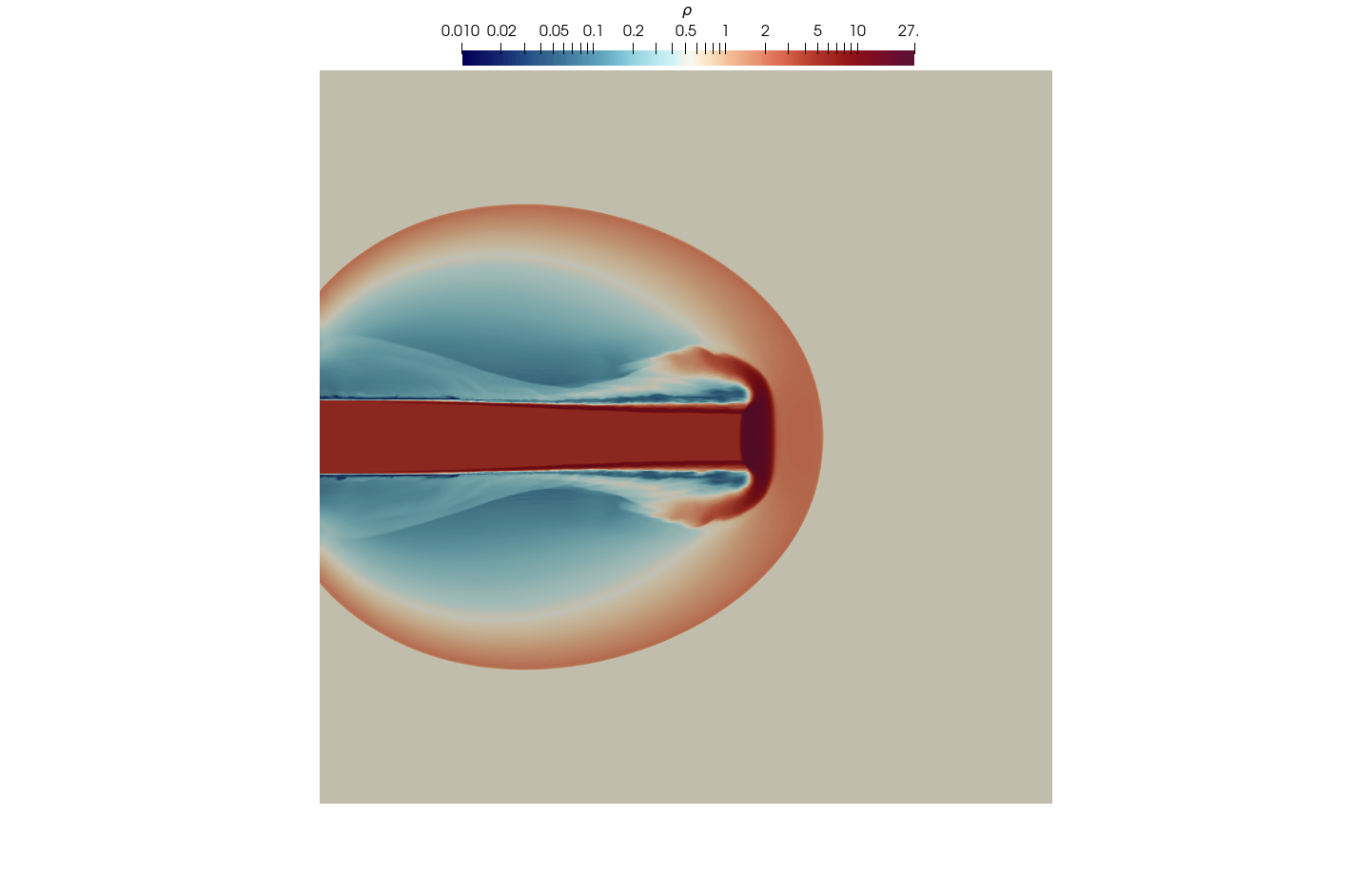}
	     \put (3,90) {CFL=0.09}
    \end{overpic}
	\caption{MCL}
\end{subfigure}
\end{minipage}
\begin{minipage}[t]{0.32\textwidth}
\begin{overpic}[trim=336 72 336 72 ,clip,width=\linewidth]{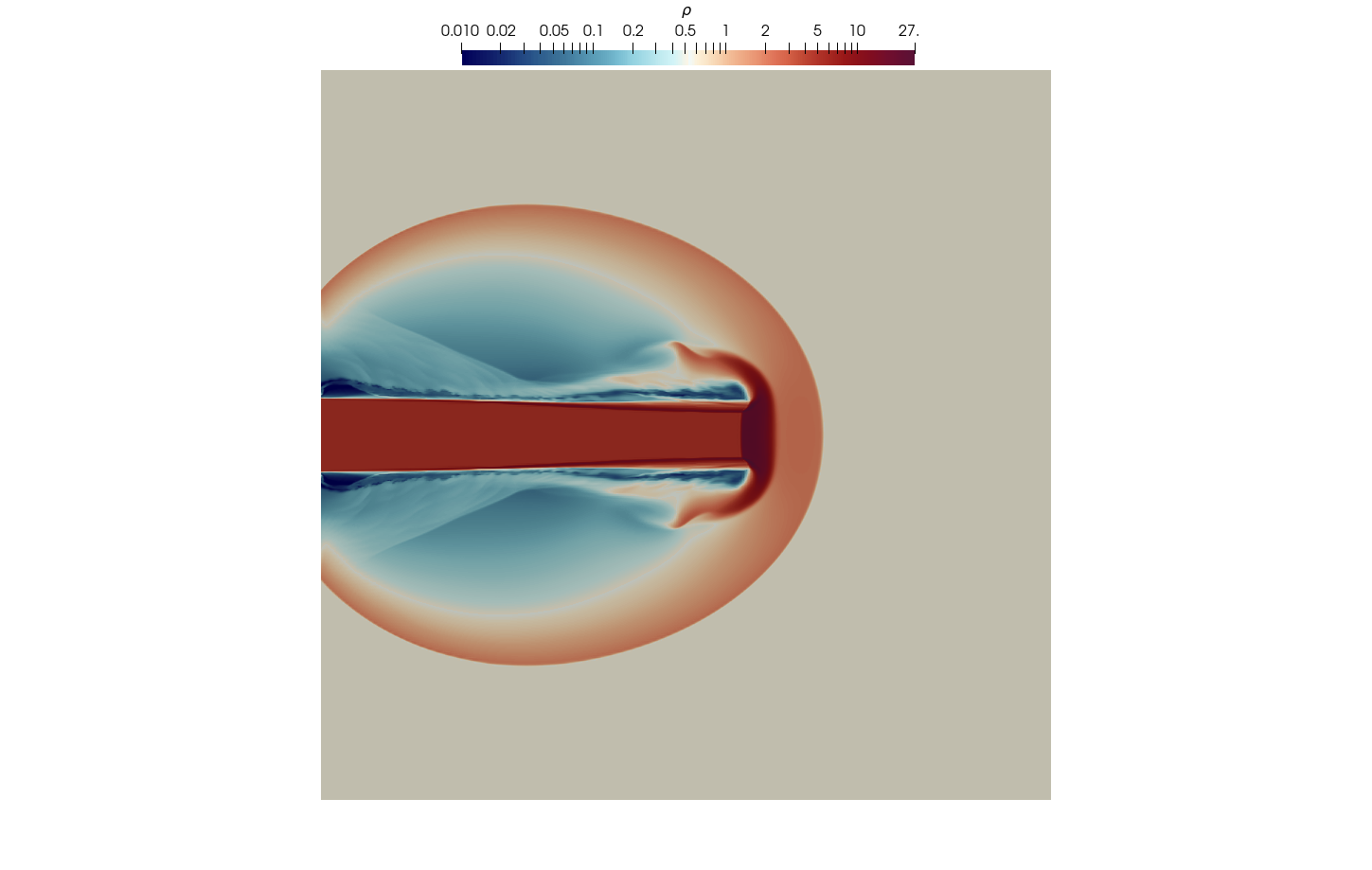}
    \put (3,90) {CFL=0.9}
\end{overpic}
\begin{overpic}[trim=336 72 336 72 ,clip,width=\linewidth]{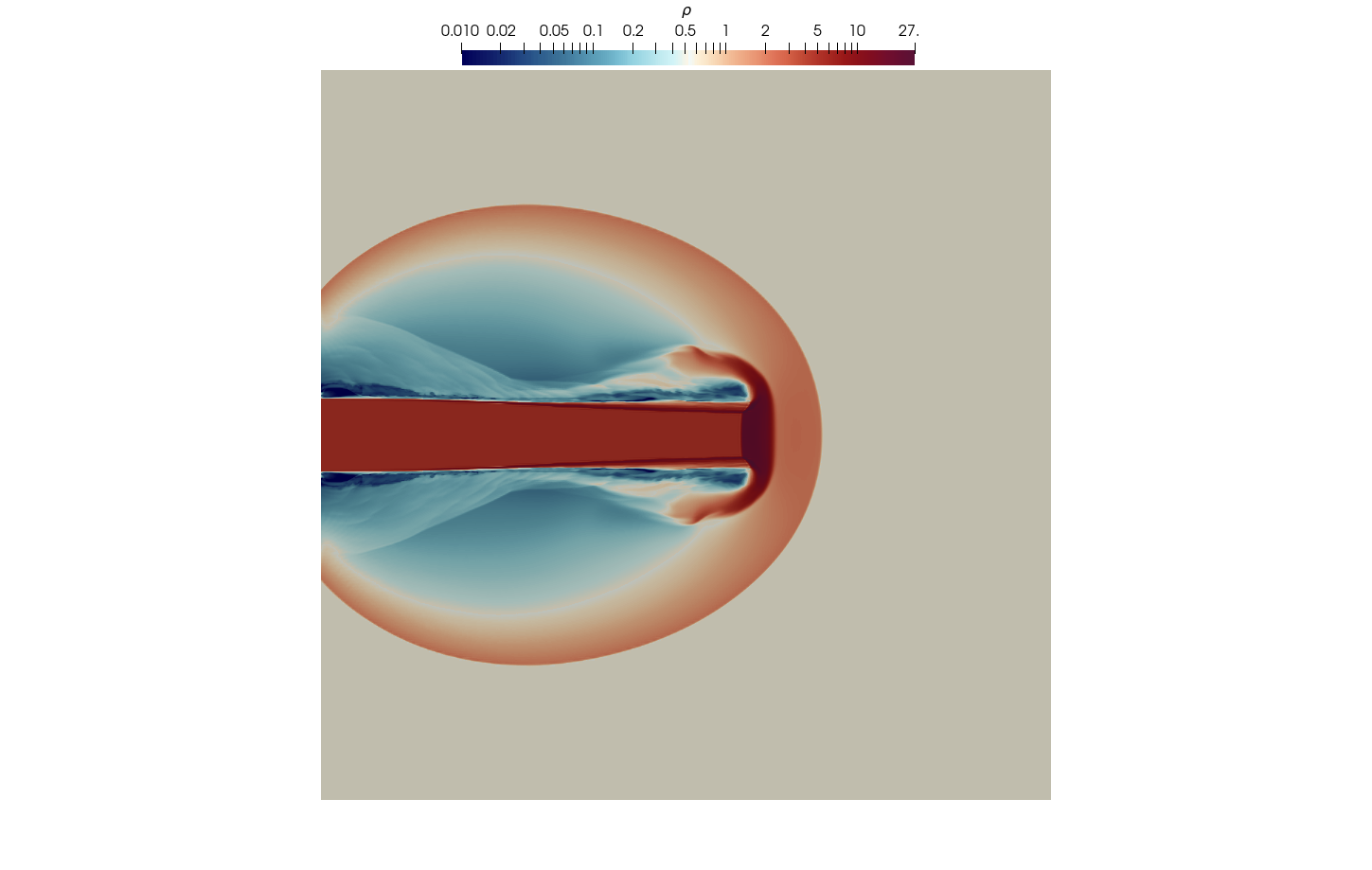}
    \put (3,90) {CFL=0.45}
\end{overpic}
\begin{overpic}[trim=336 72 336 72 ,clip,width=\linewidth]{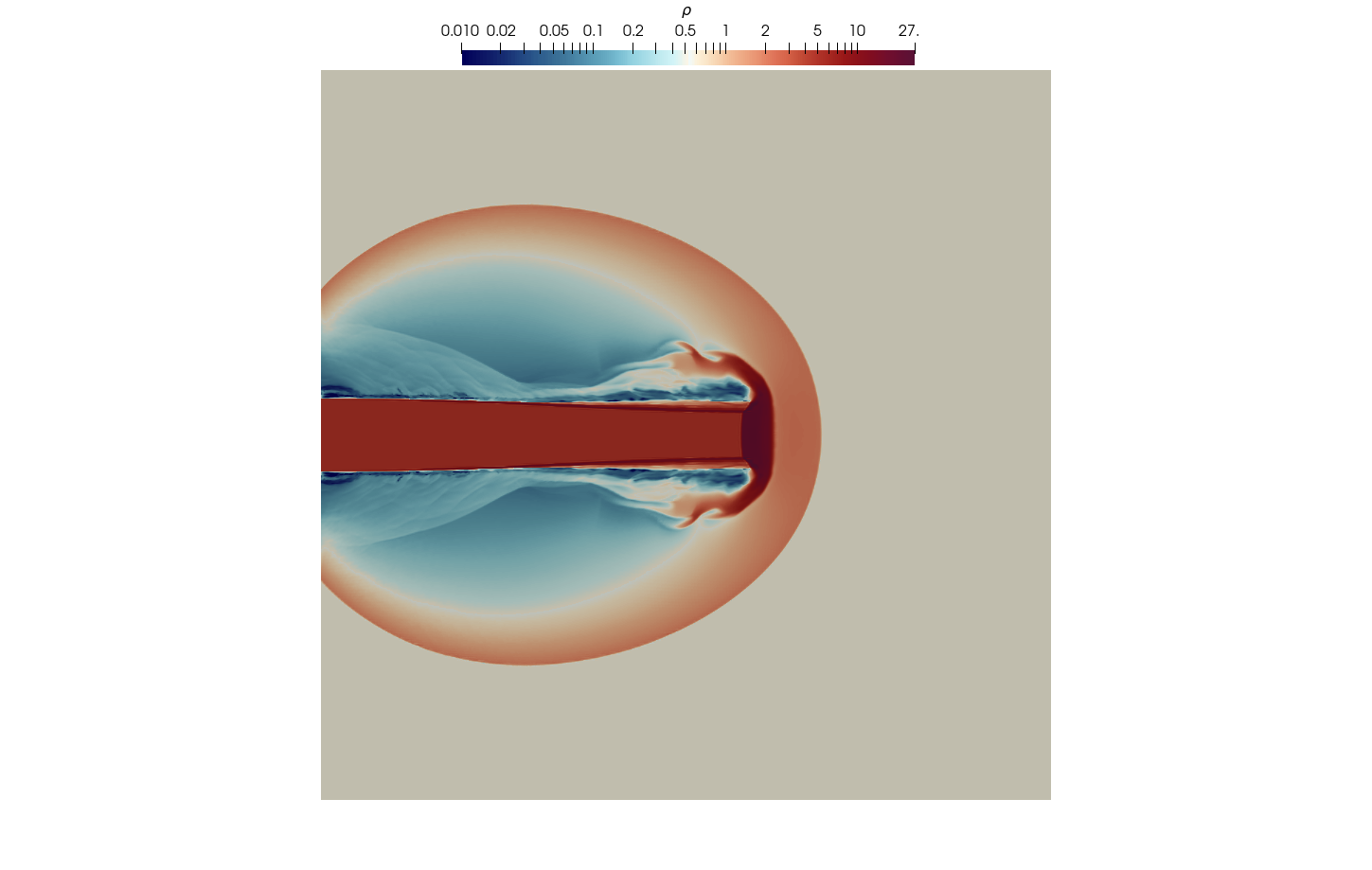}
    \put (3,90) {CFL=0.225}
\end{overpic}
\begin{subfigure}[b]{\linewidth}
    \centering
	\begin{overpic}[trim=333 72 333 72 ,clip,width=\linewidth]{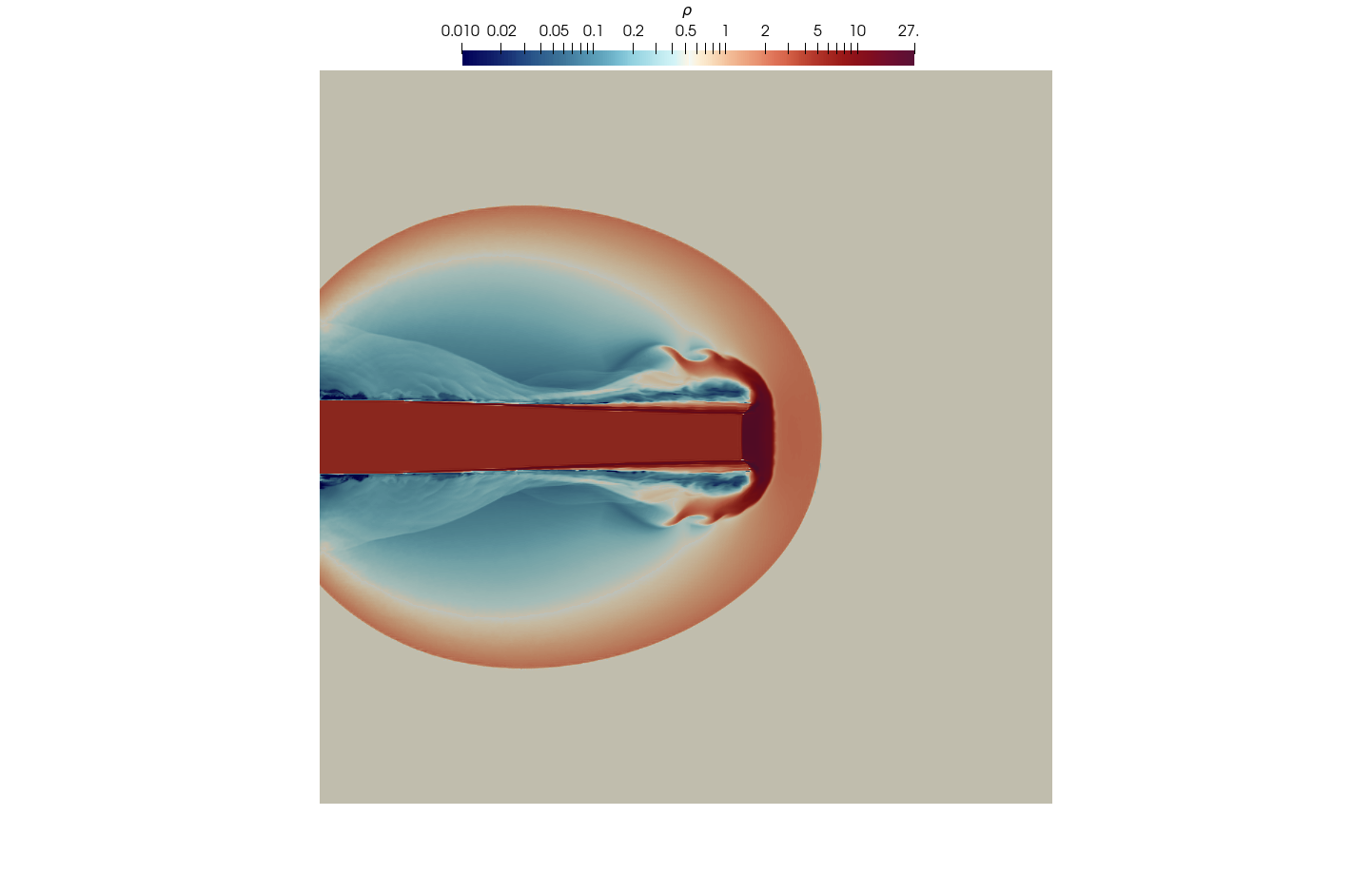}
	    \put (3,90) {CFL=0.09}
    \end{overpic}
	\caption{FCT/IDP}
\end{subfigure}
\end{minipage}

\caption{Density contours for the astrophysical jet simulations with MCL limiter \textbf{C} (density factor for all equations, sequential limiting, pressure and entropy limiters) and FCT/IDP (density and specific entropy) strategies.
DGSEM results with polynomial degree $N=3$ and $256\times 256$ elements.}
\label{fig:astrojet_cfl}
\end{figure}

We first compare the MCL limiter variant \textbf{C} from the previous section with the FCT/IDP method with local bar-state bounds for the density and specific entropy \eqref{eq:IDPbounds} at different CFL numbers.
Figure~\ref{fig:astrojet_cfl} shows that the amount of vortical structures in the density contours at $t=10^{-3}$ obtained with the FCT/IDP method is highly dependent on the CFL number, whereas the simulations that use the MCL limiter show a weaker dependence on the CFL number.
Table~\ref{tab:astrojet_timesteps} further shows that the total number of time steps needed to reach $t=10^{-3}$ depends linearly on the CFL number for MCL methods, but not for FCT/IDP methods.

The dependence of the spatial discretization on the time-step size for FCT/IDP methods causes the number of vortical structures to be highly dependent on the CFL number and the total number of time steps to be not inversely proportional to the CFL number.
In fact, the amount of dissipation is reduced for small CFL numbers, which leads to lower minimum densities and higher maximum pressures in FCT/IDP, as can be seen in Table~\ref{tab:astrojet_timesteps}.
Something like a feedback effect occurs when the lower densities and higher pressures increase the speed of sound in the medium, which in turn reduces the time-step size even more, which again reduces dissipation (indirectly, due to lowering the bounds of the \textit{a-posteriori} limiting approach in FCT/DIP). While having reduced dissipation is in general of course desirable, here it is more subtle as one buys the low dissipation with a strongly increased number of time steps, i.e., with a strongly increased CPU time. 

\begin{table}[]
    \centering
    \caption{Total number of time steps, density and pressure range at $t=10^{-3}$ for the MCL and FCT/IDP simulations of the astrophysical jet as a function of the CFL number}
    \begin{tabular}{c|cccc}
        Method & CFL & $\#$ time steps & $[\rho^{\min},\rho^{\max}]$ & $[p^{\min},p^{\max}]$\\
        \hline
\multirow{4}{*}{FCT/IDP}&0.9& $10758$ & $[4.43\times 10^{-3}, 32.15]$ & $[0.412, 172640]$ \\
&0.45& $25608$ & $[2.93\times 10^{-3}, 35.64]$ & $[0.409, 173935]$ \\
&0.225& $71577$ & $[1.56\times 10^{-3}, 32.62]$ & $[0.407, 217273]$ \\
&0.09& $241974$ & $[7.11\times 10^{-4}, 40.67]$ & $[0.405, 228165]$ \\
        \hline
\multirow{4}{*}{MCL}&0.9& $10284$ & $[1.07\times 10^{-2}, 24.71]$ & $[0.217, 174801]$ \\
&0.45& $20568$ & $[8.73\times 10^{-3}, 25.41]$ & $[0.279, 175133]$ \\
&0.225& $41138$ & $[7.40\times 10^{-3}, 23.85]$ & $[0.047, 175271]$ \\
&0.09& $102844$& $[1.33\times 10^{-2}, 24.64]$& $[0.126, 175304]$\\
        \hline
    \end{tabular}
    \label{tab:astrojet_timesteps}
\end{table}

\begin{figure}[h!]
\centering
\includegraphics[trim=450 846 450 0 ,clip,width=0.5\linewidth]{figures/astrojet/mcl_densityforall_seq_presExact_cfl_0.45.png}
\\
\includegraphics[trim=394 73 394 72 ,clip,width=0.24\linewidth]{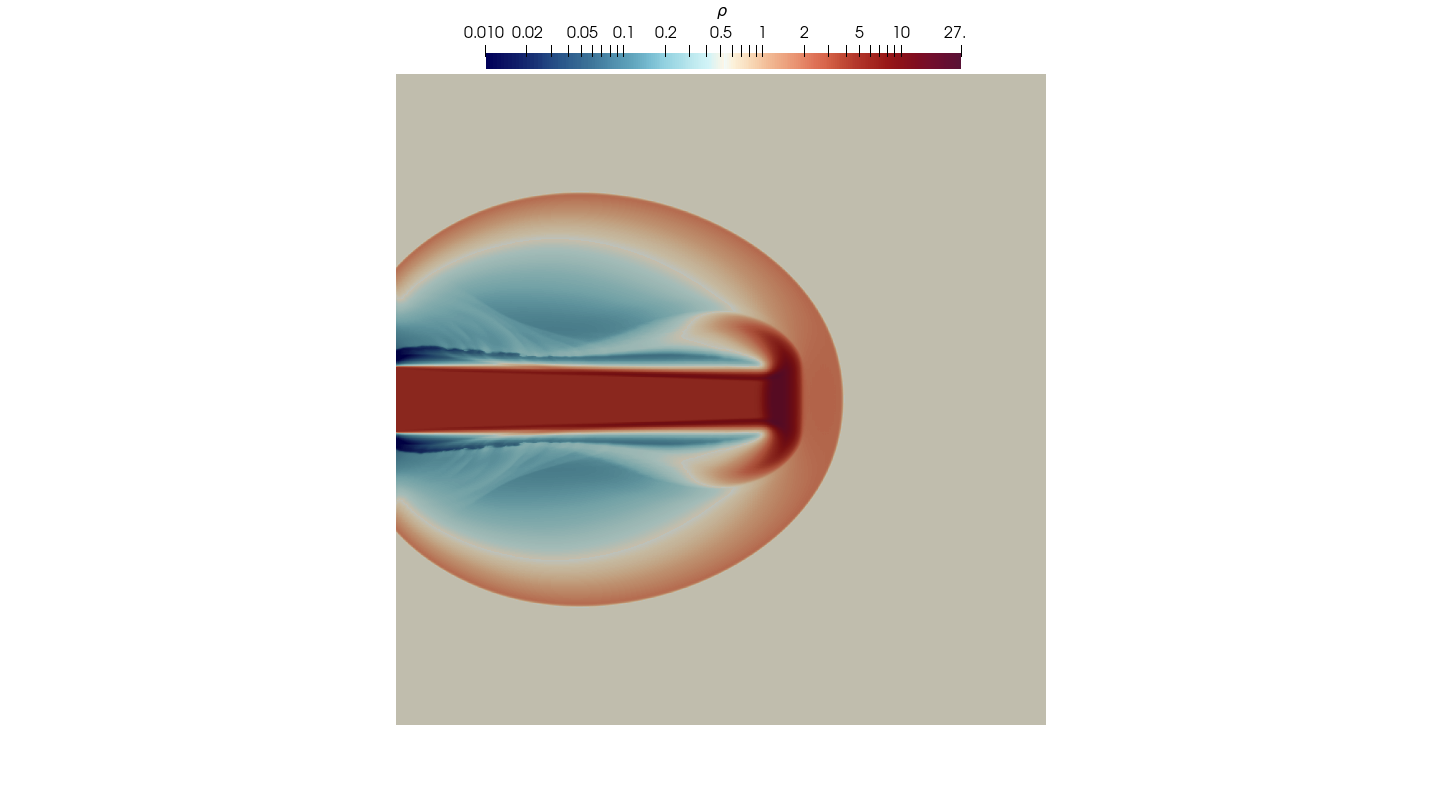}
\includegraphics[trim=336 73 336 72 ,clip,width=0.24\linewidth]{figures/astrojet/mcl_densityforall_seq_presExact_cfl_0.9.png}
\includegraphics[trim=394 73 394 72 ,clip,width=0.24\linewidth]{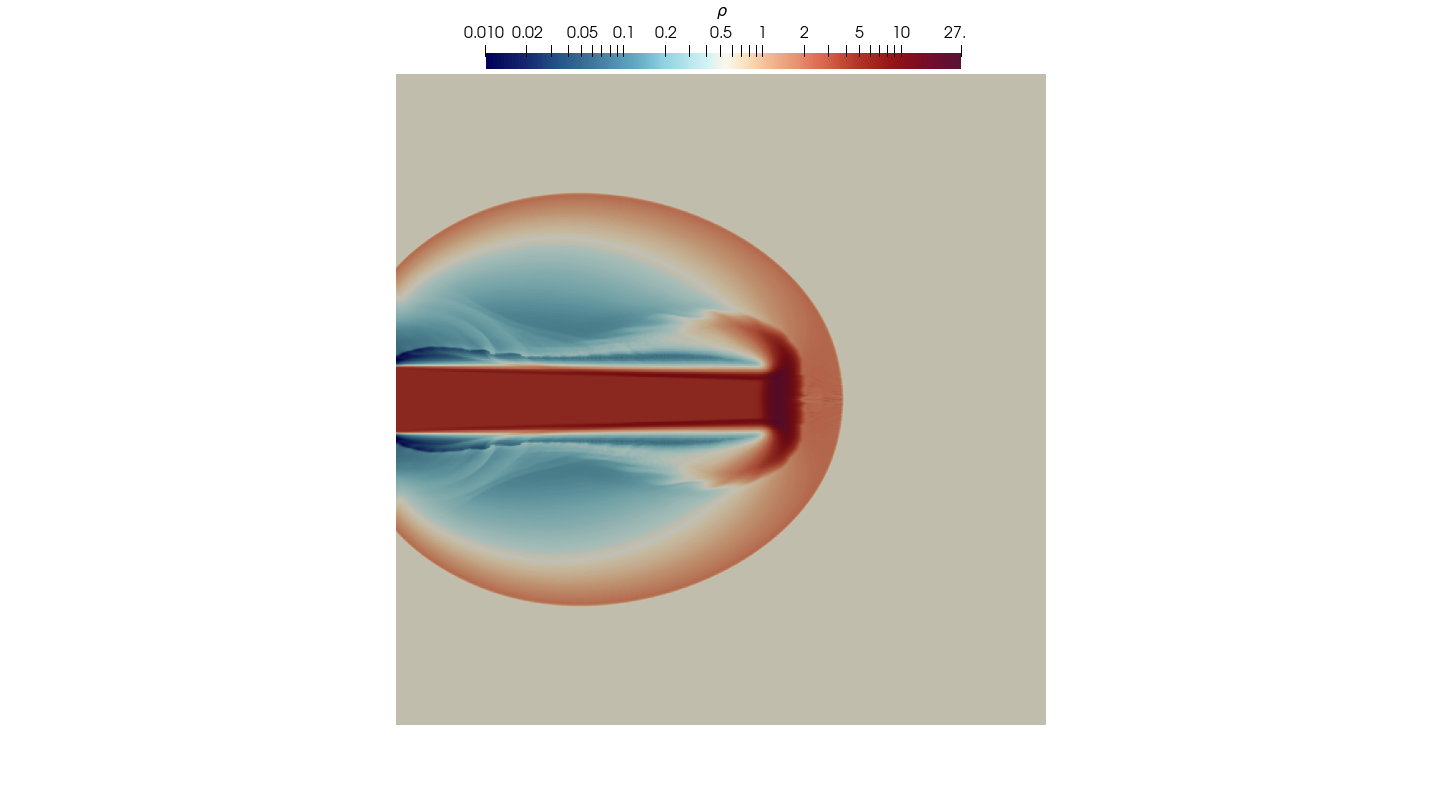}
\includegraphics[trim=394 73 394 72 ,clip,width=0.24\linewidth]{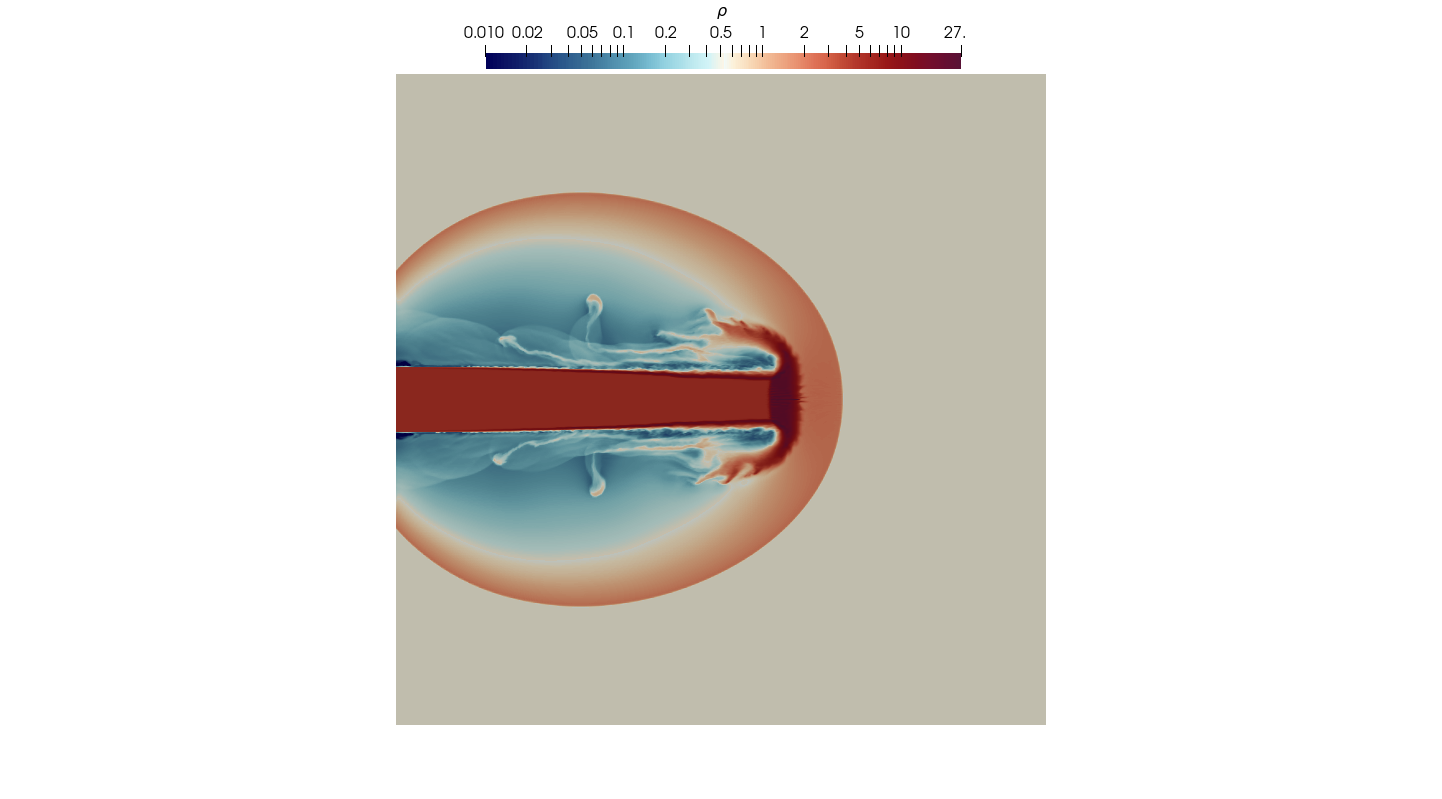}

\includegraphics[trim=463 726 463 0 ,clip,width=0.5\linewidth]{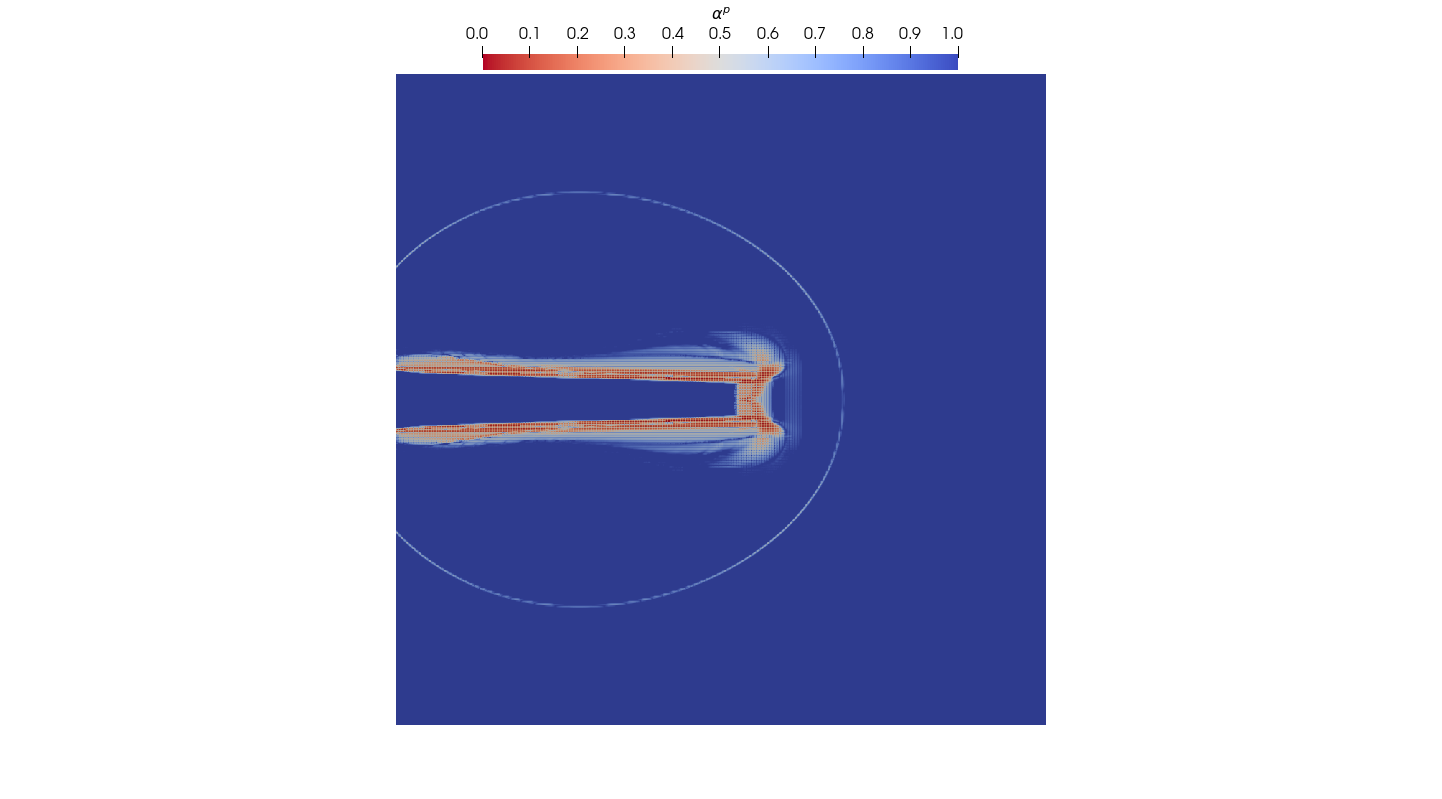}
\\
\begin{subfigure}[b]{0.24\linewidth}
    \centering
	
 	\includegraphics[trim=394 73 394 72 ,clip,width=\linewidth]{figures/astrojet/mcl_densityforall_seq_pres_cfl_0.9_alpha_pressure.png}
	\caption{MCL \textbf{C} cautious}
\end{subfigure}
\begin{subfigure}[b]{0.24\linewidth}
    \centering
 	\includegraphics[trim=394 73 394 72 ,clip,width=\linewidth]{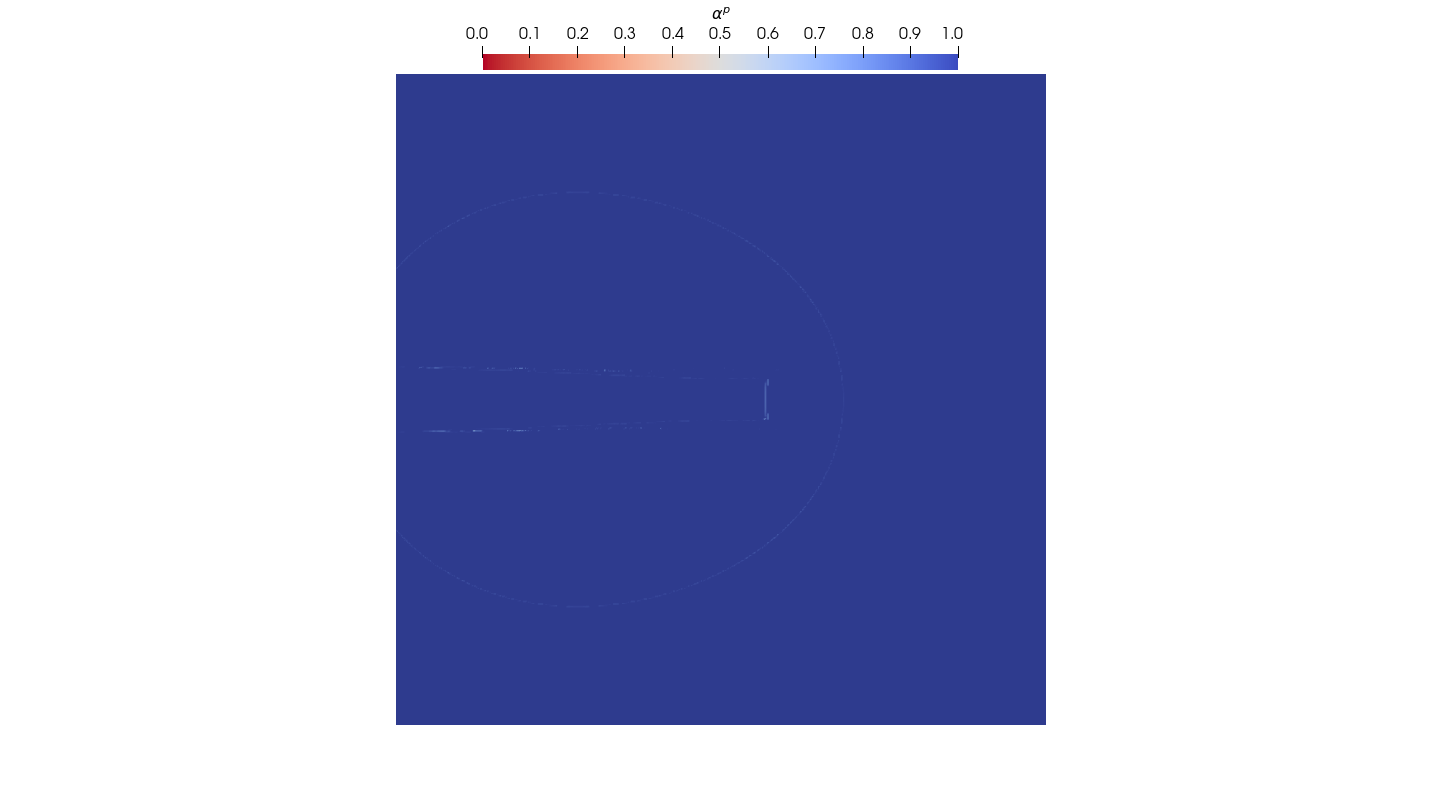}
	\caption{MCL \textbf{C} sharp}
\end{subfigure}
\begin{subfigure}[b]{0.24\linewidth}
    \centering
 	\includegraphics[trim=394 73 394 72 ,clip,width=\linewidth]{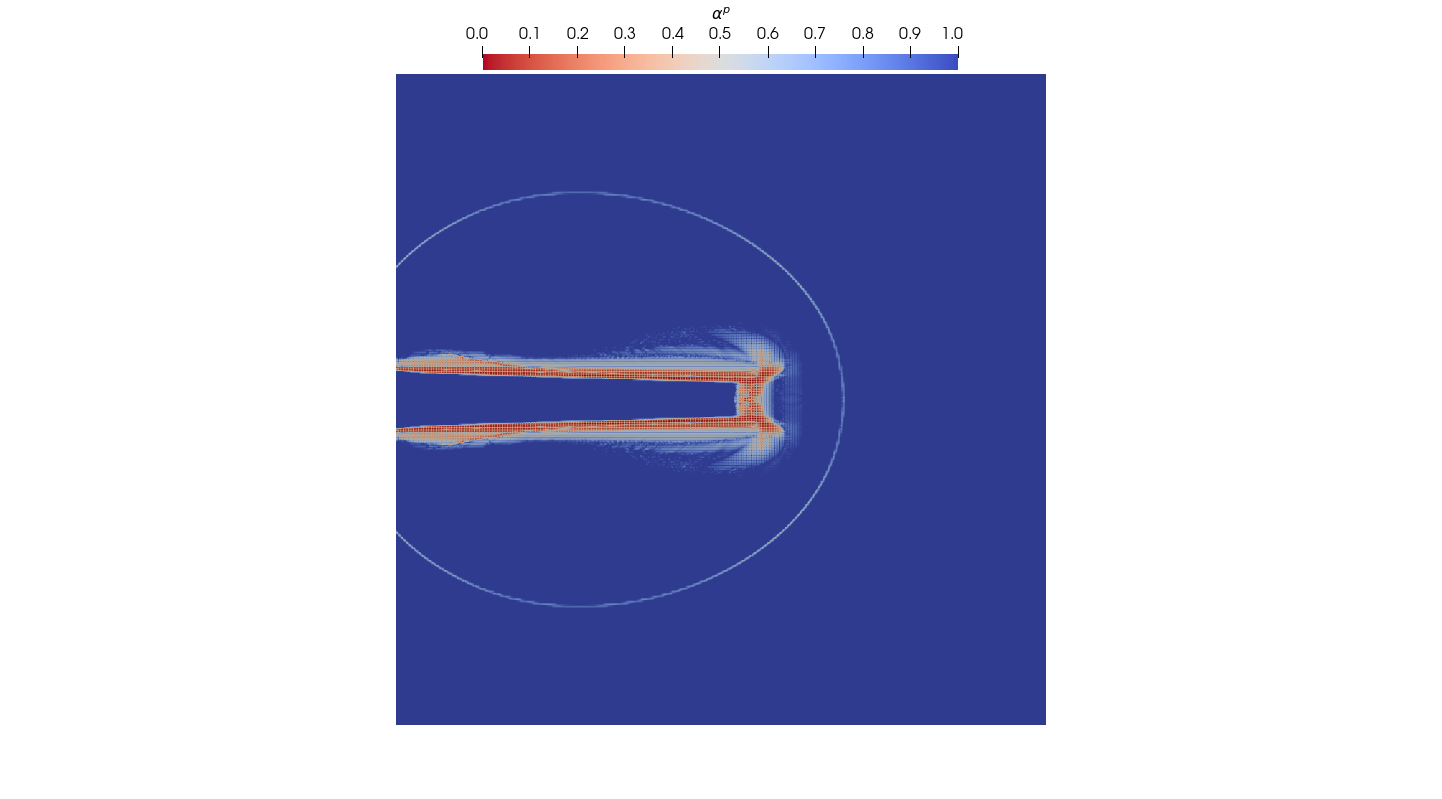}
	\caption{MCL \textbf{A} cautious}
\end{subfigure}
\begin{subfigure}[b]{0.24\linewidth}
    \centering
 	\includegraphics[trim=394 73 394 72 ,clip,width=\linewidth]{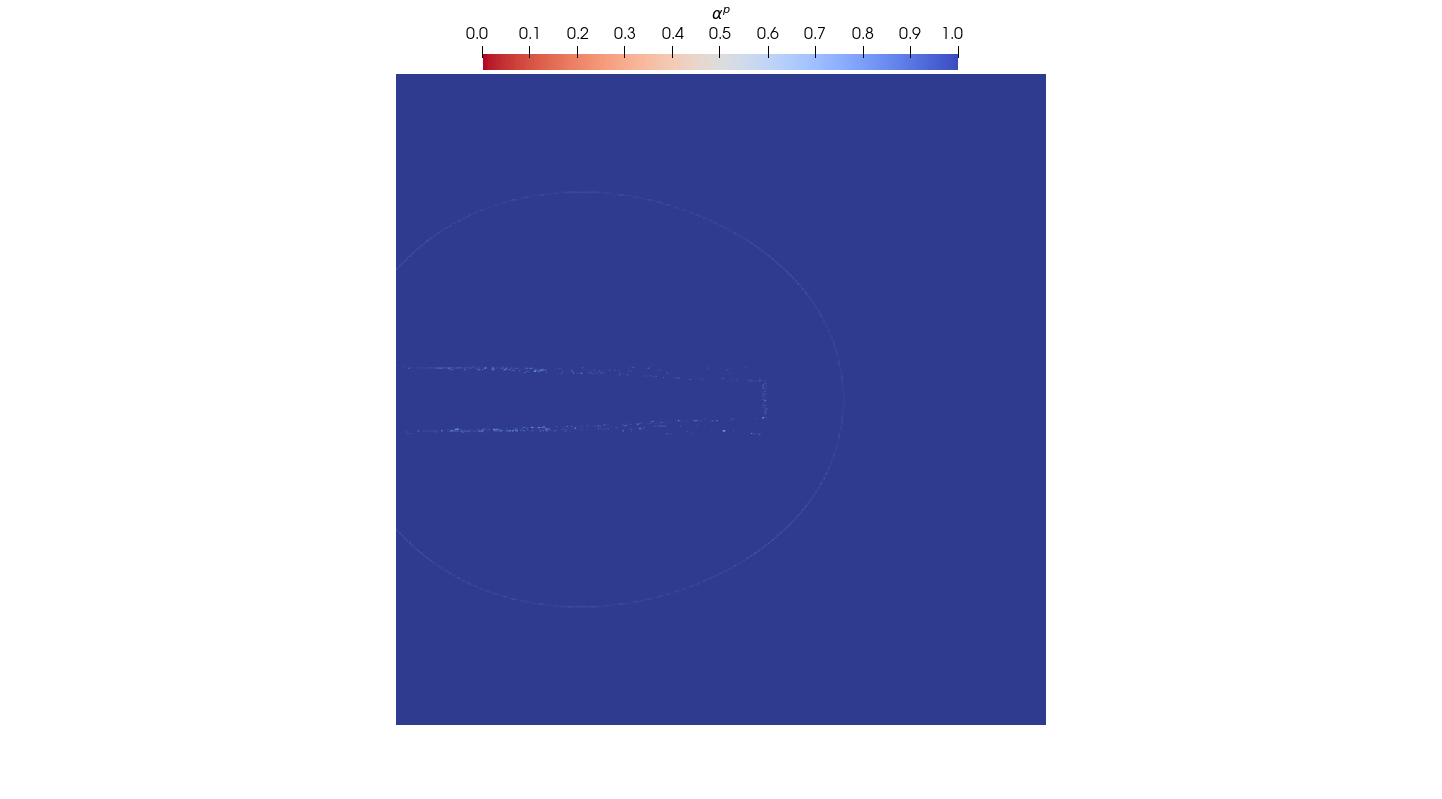}
	\caption{MCL \textbf{A} sharp}
\end{subfigure}


\caption{Density contours for the astrophysical jet simulations with MCL limiters \textbf{A} and \textbf{C} using the cautious and sharp pressure limiters.
DGSEM results with polynomial degree $N=3$ and $256\times 256$ elements. CFL=0.9.}
\label{fig:astrojet_simplfiedpressure}
\end{figure}

Finally, we compare the difference between the sharp pressure positivity limiter (\eqref{eq:pressureLimiter} with \eqref{eq:pressureLimiterExact}) and the cautious pressure positivity limiter (\eqref{eq:pressureLimiter} with \eqref{eq:pressureLimiterSimp}).
Figure~\ref{fig:astrojet_simplfiedpressure} shows the density contours and pressure limiter limiting factors at $t=10^{-3}$ obtained with the MCL limiters \textbf{A} and \textbf{C}, and the cautious and sharp pressure positivity limiters.
For this particular case, it is clear that the cautious pressure positivity limiter adds significant dissipation in the shear layer of the yet, which suppresses the appearance of vortical structures for both MCL limiters.
On the other hand, the sharp pressure positivity limiter adds just enough dissipation to maintain the pressure of all bar states non-negative, and hence allows the development of turbulence.

\section{Conclusion}\label{sec:conclusions}

In this paper, we have extended the monolithic convex limiting method to nodal discontinuous Galerkin methods (DGSEM) that use Legendre--Gauss--Lobatto (LGL) points.
We have shown that the collocated nature of the LGL-DGSEM approximation greatly simplifies the design of MCL limiters and the derivation of compatible low-order invariant domain preserving schemes, which are needed to apply the MCL strategy.
We have demonstrated the versatility of LGL-DGSEM/MCL methods to solve challenging simulation setups featuring supersonic, hypersonic and turbulent flow regimes.
We have compared the performance of MCL methods and predictor-corrector-type flux corrected transport (FCT) subcell limiting methods.
Unlike FCT methods, the amount of dissipation (and hence the spatial discretization) obtained with MCL techniques does not depend on the time-step size.
As a result, time convergence can be expected for MCL but not for FCT schemes in problems that need stabilization.

\backmatter

\bmhead{Supplementary information}.

The methods used in this paper were implemented in the open-source high-order DG code Trixi.jl \cite{ranocha2022adaptive,schlottkelakemper2021purely,schlottkelakemper2020trixi}.
We refer the interested reader to our reproducibility repository (\url{https://github.com/amrueda/paper_2023_MCL_LGL-DGSEM}), where we provide detailed instructions of how to reproduce the numerical results that we present here.

\bmhead{Acknowledgments}

The authors thank Dr. Hennes Hajduk for performing comparative studies with his code and giving a deeper insight into the subcell DG-MCL schemes he developed in \cite{hajduk2021monolithic} for high-order Bernstein finite elements. 
Moreover, the authors would like to thank Prof. Dr. Hendrik Ranocha for his support and advice during the implementation of FCT methods in Trixi.jl \cite{ranocha2022adaptive,schlottkelakemper2021purely,schlottkelakemper2020trixi}.

Gregor Gassner and Andrés M. Rueda-Ramírez acknowledge funding through the Klaus-Tschira Stiftung via the project ``HiFiLab''.
Gregor Gassner further acknowledges funding by the German Research Foundation (DFG) under the grant number DFG-FOR5409. Dmitri Kuzmin acknowledges DFG support under grant number KU 1530/23-3.

We furthermore thank the Regional Computing Center of the University of Cologne (RRZK) for providing computing time on the High Performance Computing (HPC) system ODIN, as well as for technical support.

\section*{Declarations}

The authors have no competing interests to declare that are relevant to the content of this article.

\begin{appendices}

\section{Euler Equations of Gas Dynamics}\label{app:euler}
The Euler equations describe the conservation of mass, momentum, and total energy per unit volume, $\stateL{u} = [\rho, \rho \vec{v}, \rho E]^T$.
The conservation law reads
\begin{equation} \label{eq:euler}
\bigpartialderiv{\stateL{u}}{t} 
+ \Nabla \cdot
\begin{pmatrix}
\rho \vec{v} \\
\rho \vec{v} \otimes \vec{v} + \mat{I} p \\
\vec{v} (\rho E + p)
\end{pmatrix}
= \stateL{0},
\end{equation}
where $\mat{I}$ is the identity matrix, the pressure is computed with the calorically perfect gas assumption,
\begin{equation} \label{eq:press}
p=(\gamma-1) \rho e,
\end{equation}
$\gamma$ is the heat capacity ratio, and $e = E - \norm{\vec{v}}^2/2$ is the internal energy.

With the physical assumption of positive density and pressure, $\rho,p>0$, a suitable, strictly convex entropy function for the compressible Euler equations is the thermodynamic entropy density divided by the constant $-(\gamma-1)$ \cite{barth1999numerical,Tadmor2003},
\begin{equation}
S(\stateL{u}) = - \frac{\rho s}{\gamma-1},
\label{entropy}
\end{equation}
where $S$ is the so-called mathematical entropy and $s = \ln\left(p \rho^{-\gamma}\right)$ is the thermodynamic entropy.
From the entropy function, we define the entropy variables,
\begin{equation} \label{eq:entrvars}
\entVar = \frac{\partial S}{\partial \stateL{u}} = \left(\frac{\gamma-s}{\gamma-1} - \beta \norm{\vec{v}}^2,~2\beta v_1,~2\beta v_2,~2\beta v_3,~-2\beta \right)^T,
\end{equation}
with $\beta = \frac{\rho}{2p}$, a quantity that is proportional to the inverse temperature.

If we contract \eqref{eq:euler} with the entropy variables, we obtain the entropy conservation law if the solution is smooth \cite{barth1999numerical,Tadmor2003},
\begin{equation}\label{eq:EntropyConservationLaw}
 \frac{\partial S}{\partial t} + \Nabla \cdot \vec{f}^{\,S} = 0,
\end{equation}
where $\vec{f}^{\,S} = \vec{v} S$ is the so-called entropy flux.

Furthermore, in the presence of discontinuities in the solution, the contraction of the compressible Euler equations with the entropy variables leads to an entropy inequality \cite{barth1999numerical,Tadmor2003},
\begin{equation}\label{eq:EntropyInequalityWeak}
\frac{\partial S}{\partial t} + \Nabla \cdot \vec{f}^{\,S} \le 0,
\end{equation}
where the total mathematical entropy within any physical domain, $\Omega$, can only increase over time if it is transported into the domain through its boundaries,	$\partial \Omega$.
Equation \eqref{eq:EntropyInequalityWeak} is the mathematical description of the second law of thermodynamics.

Finally, the entropy flux potential is defined as \cite{barth1999numerical,Tadmor2003}
\begin{equation}\label{eq:entPotential}
\vec{\Psi} := \entVar^T \blocktensor{f}^a -\vec{f}^S.
\end{equation}

\end{appendices}

\bibliography{sn-bibliography}


\end{document}